\documentclass[11pt]{article}

\usepackage[mathscr]{eucal}
\usepackage{color}
\usepackage{amsmath}
\usepackage{amssymb}
\usepackage{amsfonts}
\usepackage{amscd}
\usepackage{amsthm}
\usepackage{latexsym}
\usepackage{graphicx}
\usepackage{subfig}

\def\be{\begin{equation}}
\def\ee{\end{equation}}
\def\bse{\begin{subequations}}
\def\ese{\end{subequations}}

\usepackage{latexsym}

\let\er\eqref

\let\be\beta

\newcommand{\R}{{\mathbb R}}

\newtheorem{theorem}{Theorem}
\newtheorem{lemma}[theorem]{Lemma}
\newtheorem{proposition}[theorem]{Proposition}
\newtheorem{remark}[theorem]{Remark}

\newtheorem{definition}[theorem]{Definition}

\def\bse{\begin{subequations}}
\def\ese{\end{subequations}}

\usepackage{geometry}

\title{On the cauchy problem with degenerate diffusion and nonlocal nonlinear sources}

\author{ Shen Bian\footnote{Email: \texttt{bianshen66@163.com}. }
 \\
  \  
  \\
Beijing University of Chemical Technology, 100029, Beijing. }
\date{}

\begin{document}
\let\cleardoublepage\clearpage

\maketitle

\begin{abstract}
This paper is devoted to the analysis of non-negative solutions for a generalisation of the parabolic equation with porous medium like nonlinear diffusion and nonlinear nonlocal reaction. We investigate under which conditions equilibration between two competing effects, repulsion modelled by nonlinear diffusion and aggregation modelled by nonlinear reaction, occurs. Precisely, we exhibit that the qualitative behavior of solutions is decided by the nonlinear diffusion which is chosen in such a way that its scaling and the reaction term coincide, i.e. that there is a critical exponent $m+2/n$ for the reaction exponent $\alpha,$ solutions exist globally with uniformly upper bounds in the case of (i)$1\le \alpha<m+2/n$ for any initial data, (ii) $\alpha>m+2/n$ for small initial data and (iii) $\alpha=m+2/n$ for small mass capacity $M_0$. In the case of (ii) and (iii), the decay properties of the solution are also discussed. Moreover, numerical simulations are carried out to verify the theoretical analysis and explore other issues that lie beyond the scope of the analysis.
\end{abstract}

\section{Introduction}\label{sec1}
\def\theequation{1.\arabic{equation}}\makeatother
\setcounter{equation}{0}
\def\thetheorem{1.\arabic{theorem}}\makeatother
\setcounter{theorem}{0}

In this work, we analyse qualitative properties of non-negative solutions in dimension $n \ge 3$ for the degenerate equation of the type
\begin{align}\label{nkpp}
\left\{
  \begin{array}{ll}
    u_t=\Delta u^m+\chi u^\alpha\left(M_0- \int_{\R^n} u(x,t)dx\right),\quad   & x \in \R^n, t>0, \\
 u(x,0)=u_0(x) \geq 0 \in L^1 \cap L^\infty(\R^n),
  \end{array}
\right.
\end{align}
where $\chi>0, m>1, \alpha \ge 1$. \er{nkpp} is related to many equations arising from population dynamics \cite{f37,kpp}, $u$ is the density of the population. The purpose of nonlinear diffusion $\Delta u^m$ with $m>1$ is to model the local repulsion of population, this can be interpreted as taking into account anti-crowding effects \cite{cal06}. The reaction term presents a growth factor of logistic type defined in terms of the total mass of the population which is a competitive term limiting such growth, where the resources of the environment can be consumed nonlocally. $\chi u^\alpha$ can also be interpreted as the nonlinear source in the process of diffusion \cite{AD}, which is called heat source for $\chi>0$ and cold source for $\chi<0.$ In the case of $\alpha=1$, the coefficient $\chi M_0$ is sometimes called Malthusian parameter which induces an exponential growth for low density populations. The case $\alpha=2$, which is the motivation of this work, considers the addition of sexual reproduction to the model with the reproduction rate proportional to the square of the density \cite{VVpp}. Nonlocal type reaction terms can also describe Darwinian evolution of a structured population density or the behavior of cancer cells with therapy as well as polychemotherapy \cite{Lorz:2011hl,Lorz:2013vp}. The main feature of this class of equations is the interplay between the degeneracy in the principal part and the growth of the forcing term.

A fundamental property of the solutions to \er{nkpp} is the formal boundedness of the total mass of the system
\begin{align}
m(t)=\int_{\R^n} u(x,t) dx
\end{align}
which satisfies
\begin{align}
\frac{d}{dt}m(t)=\chi \left(M_0-m(t) \right) \int_{\R^n} u^\alpha dx.
\end{align}
If the initial mass $m_0:=\int_{\R^n} u_0 dx>M_0$, then $m(t)$ decreases in time and $M_0 \le  m_0$ for all $t>0.$ Thus we find that $u(x,t)$ is a subsolution of the porous medium equation $v_t=\Delta v^m$ which admits a global solution for any $m>1$ \cite{JLV03,JLV06}. By the comparison principle,  all solutions of \er{nkpp} exist globally. When the initial mass $m_0<M_0$, then $m(t)$ increases in time and
$ m_0 \le m(t) \le M_0.$ Therefore, we assume that the initial mass satisfies
$$ m_0=\int_{\R^n} u_0 dx<M_0 $$
throughout this paper. In this sense, $M_0$ can be considered as the carrying capacity \cite{bp07}.

In any dimension $n \ge 3,$ we will concentrate on a particular choice of the nonlinear reaction exponent
\begin{align*}
\alpha=m+2/n
\end{align*}
which produces a balance in the mass-invariant scaling of diffusion and reaction. Indeed, let $u_\lambda(x)=\lambda^n u(\lambda x,t)$ of same mass as $u$, the diffusion term $\lambda^{2+nm} \Delta u_\lambda^m$ has the same scaling as the reaction term $\lambda^{n \alpha} u_\lambda^\alpha \left( M_0-\int_{\R^n}  u_\lambda dx \right)$ if and only if $nm+2=n\alpha$ or equivalently $\alpha=m+2/n.$ The dynamics in \er{nkpp} are governed by the interaction between nonlinear diffusion and reaction. We are interested in three cases: critical case $\alpha=m+2/n$, subcritical case $\alpha<m+2/n$, supercritical case $\alpha>m+2/n$. Under what conditions which of the two competing items dominates will be explored in this paper.

\subsection{Comments on the non-degenerate case $m=1$}

In the non-degenerate case $m=1,$ there is a vast body of literature on the semi-linear reaction-diffusion equations $u_t=\Delta u+F(u)$ in bounded domain, cf. the papers \cite{BCL17,hy95,lcl09,Lorz:2011hl,Lorz:2013vp,vol1,vol2,VVpp,ww11,ww96}. For instance, the reaction term is $F(u)=\int_{\Omega} u^p dx-\beta u^q$ \cite{ww96} where the competitive effect of the local term $u^q$ becomes more influential as population grows such that the equation possesses the comparison principle which helps in proving the existence of global solutions by virtue of the boundedness of $\Omega$, or $F(u)=u^p-\frac{1}{|\Omega|} \int_{\Omega} u^p dx$ \cite{hy95} where the equation is equipped with a decreasing Lyapunov functional. Yet there are few results on this type of equations in the whole space. This is partially due to the apparent lack of a good Lyapunov functional and the unboundedness of the domain. More importantly, the comparison principle (which has been used in the equation with local reaction term for the global existence \cite{WMZ14}) is no longer applicable to our model with nonlocal reaction term.

Before turning to the nonlocal term in the whole space, we firstly mention the following fundamental work of Fujita \cite{fuj1}, he proved that for $1<\alpha<1+2/n,$ the local classical solution blows up in finite time (the same is true for $\alpha=1+2/n$ \cite{kh73,tana77}). The natural guess is that if $M_0-\int_{\R^n} u dx$ remains positive, our model has similar structure to Fujita equation. However, the result obtained in our previous paper \cite{BCL15} gives an opposite consequence. That's, for $1<\alpha<1+2/n$, our model admits a global classical solution as long as the initial value $u_0(x)$ is nontrivial.

\subsection{Our results for the degenerate case $m>1$}
In the degenerate case $m>1,$ before proceeding further, let us state the notion of weak solutions we will deal throughout this paper with:
\begin{definition}\label{weakdefine}(Weak solution)
Let $u_0$ be an initial condition satisfying
\begin{align}\label{initialdata}
u_0 \in L^1 (\R^n;(1+|x|^2)dx ) \cap L^\infty(\R^n),~~\nabla u_0^m \in L^2(\R^n),~~u_0 \ge 0,~~\int_{\R^n} u_0 dx <M_0
\end{align}
and $T \in (0,\infty]$. The non-negative functions defined in $\R^n \times [0,T)$ is called a weak solution of \er{nkpp} on $[0,T)$ if
\begin{enumerate}
\item[\textbf{(i)}] Regularity:
\begin{align}
u \in L^\infty(0,T;L^1 \cap L^\infty(\R^n)),~u^m \in L^2(0,T;H^1(\R^n)).
\end{align}
\item[\textbf{(ii)}]
$u$ satisfies the equation in the sense of distribution, i.e. that
\begin{align}
& \int_0^T \int_{\R^n} \nabla u^m \cdot \nabla \varphi dxdt-\chi \int_0^T \int_{\R^n} u^\alpha \varphi dx \left( M_0-\int_{\R^n} u dx  \right) dt   \nonumber \\
& =\int_{\R^n} u_0(x) \varphi(x,0)dx+\int_0^T \int_{\R^n} u \varphi_t dxdt \label{faiweak}
\end{align}
for any continuously differentiable function $\varphi$ with compact support in $\R^n \times [0,T).$
\end{enumerate}
\end{definition}

The main results of this work can be listed as follows. For the subcritical case $1\le \alpha <m+2/n$, the following theorem gives the existence of a time global weak solution.
\begin{theorem}\label{sub}(Uniform boundedness in the subcritical case $1\le \alpha <m+2/n$)
Let $n \ge 3, m>1, 1 \le \alpha <m+2/n.$ For any $T>0,$ under assumption \er{initialdata}, there exists a weak solution $u$ to \er{nkpp} on $[0,T)$. Moreover, $u$ is uniformly bounded, i.e. that there exists a constant $C=C\left( \|u_0\|_{L^1(\R^n)},\|u_0\|_{L^\infty(\R^n)},m,\alpha,n,\chi,m_0,M_0 \right)$ such that for all $1\le k \le \infty$
\begin{align*}
\displaystyle \sup_{0<t<T} \|u(\cdot,t)\|_{L^k(\R^n)} \le C.
\end{align*}
\end{theorem}
\begin{remark}
For the subcritical case, using the mass invariant scaling, the reaction with nonlocal term dominates the diffusion for low density and prevent spreading. While for high density, the diffusion dominates the reaction, thus blow-up is precluded.
\end{remark}
\begin{remark}
As pointed out in \er{nkpp} with the absence of the nonlocal term $\chi u^\alpha \int_{\R^n} u dx$, global solutions cannot exist for $1<\alpha<m+2/n,$ see \cite{AD}. While Theorem \ref{sub} shows that the solution of \er{nkpp} exists globally without any restriction on the size of the initial data.
\end{remark}
\begin{remark}
For $\alpha=1,$ consider a solution $u$ of \er{nkpp} and define the rescaled function $v$ by:
\begin{align*}
u(x,t)=\frac{1}{R^n(t)} v\left( \frac{x}{R(t)},\tau(t) \right)=\frac{1}{R^n(t)} v(y,\tau)
\end{align*}
with
\begin{align*}
R(t)=\left( 1+\mu t \right)^{\frac{1}{\mu}},~~\tau(t)=\log R(t),
\end{align*}
where $\mu=nm-n+2$. The rescaled system is
\begin{align}\label{20210621}
\left\{
  \begin{array}{ll}
    \frac{\partial v}{\partial \tau}=\Delta v^m +\nabla \cdot (vy)+\chi v \left( M_0-\int_{\R^n} v dy  \right) e^{\mu \tau}, & y \in \R^n, ~\tau>0, \\
    v(\cdot,\tau=0)=u_0 \ge 0,  & y \in \R^n.
  \end{array}
\right.
\end{align}
Integrating \er{20210621} over $\R^n$ we obtain
\begin{align}
\left\{
  \begin{array}{ll}
    \frac{d m(\tau)}{d\tau}=\chi e^{\mu \tau} m(\tau)\left( M_0-m(\tau) \right), \\
 m(0)=m_0.
  \end{array}
\right.
\end{align}
As a consequence,
\begin{align}
M_0-m(\tau)=\frac{M_0}{\frac{m_0}{M_0-m_0}e^{-\frac{\chi M_0}{\mu}}e^{\frac{\chi M_0}{\mu}e^{\mu \tau}}+1}
\end{align}
which tells us that
\begin{align}\label{20210622}
\frac{\partial v}{\partial \tau}=\Delta v^m +\nabla \cdot (v y)+\chi v \frac{M_0 e^{\mu \tau} }{C_1 e^{C_2 e^{\mu \tau}}+1},
\end{align}
where $C_2=\frac{\chi M_0}{\mu}, C_1=\frac{m_0}{M_0-m_0}e^{-C_2}.$ For all $m>1,$ the equation \er{20210622} has a unique integrable stationary solution. Computing it we get
\begin{align}
v_{\infty,M_0}=\left( C_{M_0}-\frac{m-1}{2m}|y|^2  \right)^{\frac{1}{m-1}}.
\end{align}
Here the mass $M_0$ of the steady state $v_{\infty,M_0}$ fixes $C_{M_0}$, i.e. calculating
\begin{align*}
M_0=\int_{\R^n} v_{\infty,M_0}(y) dy
\end{align*}
one finds that
\begin{align*}
C_{M_0}=\left( \frac{m-1}{2m} \right)^{\frac{n(m-1)}{\mu}} \left( \frac{n \alpha_n}{2}\mathcal{B} \left( \frac{n}{2},\frac{m}{m-1} \right)  \right)^{-\frac{2(m-1)}{\mu}} M_0^{\frac{2(m-1)}{\mu}},
\end{align*}
where $\alpha_n=\frac{\pi^{n/2}}{\Gamma (n/2+1)}$. The behavior of the solution $u$ can be described for large $t$ by
\begin{align*}
u_{\infty,M_0}=\left( \frac{ C_{M_0}(1+\mu t)^{\frac{2}{\mu}}-\frac{m-1}{2m} |x|^2}{1+\mu t} \right)^{\frac{1}{m-1}}
\end{align*}
on expanding sets of the form $|x|<\sqrt{\frac{2mC_{M_0}}{m-1}} \left( 1+\mu t \right)^{\frac{1}{\mu}}.$
\end{remark}

Let us now discuss the critical case $\alpha=m+2/n$ in which the weak solution exists globally in time for small capacity of the total mass $M_0.$
\begin{theorem}\label{critical}(Decay properties in the critical case $\alpha=m+2/n$)
Let $n \ge 3, m>1, \alpha=m+2/n.$ Let $u_0$ be an initial data satisfying \er{initialdata} and the capacity $M_0$ satisfies
\begin{align*}
M_0 \le M_\ast
\end{align*}
where $M_\ast$ is expressed as
\begin{align}\label{Mstar}
M_\ast=\left( \frac{S_n (\alpha-m)}{\chi}\right)^{\frac{1}{\alpha-m+1}} \frac{\alpha-m+1}{\alpha-m}
\end{align}
with
\begin{align}\label{Sn}
    S_n=\frac{n(n-2)}{4}2^{\frac{2}{n}}\pi^{1+\frac{1}{n}}\Gamma\left(\frac{n+1}{2}\right)^{-\frac{2}{n}}.
\end{align}
Then there exists a weak solution $u$ to \er{nkpp} with the following decay property that for any $t>0$
\begin{align*}
\|u(\cdot,t)\|_{L^k(\R^n)} \le C_1(1+t)^{-\frac{k-1}{k(m+2/n-1)}},~~1<k<\infty
\end{align*}
and the uniform estimate
\begin{align*}
\|u(\cdot,t)\|_{L^\infty(\R^n)} \le C_2,
\end{align*}
where $C_1, C_2$ are constants depending only on $\|u_0\|_{L^1 \cap L^\infty(\R^n)}$, $m$, $n$, $\chi, m_0,M_0$.
\end{theorem}

For the supercritical case $\alpha>m+2/n,$ we present the decay property of the weak solution to \er{nkpp} under the smallness assumption on $\|u_0\|_{L^{\frac{n(\alpha-m)}{2}}(\R^n)}$. Throughout this paper, we define a constant which is related to the initial condition for the existence results:
\begin{align}\label{cp0}
  C_{p_0}=\left(\frac{n+2}{2 }\right)^{\frac{1}{p_0}}\left( \frac{4S_n m(p_0-1)(n+2)}{n\chi (p_0+m-1)^2} \right)^{\frac{1}{\alpha-m}} \frac{m_0^{\frac{1}{p_0}}}{M_0^{\frac{1}{p_0}+\frac{1}{\alpha-m}}},
\end{align}
where $S_n$ is defined by \er{Sn} and
\begin{align}\label{p0}
p_0=\frac{n(\alpha-m)}{2}.
\end{align}
\begin{theorem}\label{super}(Decay property in the supercritical case $\alpha>m+2/n$)
Let $n \ge 3, m>1, \alpha>m+2/n.$ Suppose that $u_0$ has the property \er{initialdata} satisfying
\begin{align*}
\|u_0\|_{L^{\frac{n(\alpha-m)}{2}}(\R^n)} < C_{p_0},
\end{align*}
then \er{nkpp} has a global weak solution and it holds that for any $t>0$
\begin{align*}
&\|u(\cdot,t)\|_{L^k(\R^n)} \le C_1(1+t)^{-\frac{k-1}{k(\alpha-1)}},~~1<k<\infty, \\
&\|u(\cdot,t)\|_{L^\infty(\R^n)} \le C_2,
\end{align*}
where $C_1,C_2$ are constants depending on $\|u_0\|_{L^1 \cap L^\infty(\R^n)}, n, \alpha, m,\chi, m_0,M_0$.
\end{theorem}
\begin{remark}
For the supercritical case $\alpha>m+2/n,$ from the perspective of scaling analysis, the diffusion becomes much more influential than the reaction for low density and the density has infinite-time spreading. An interesting conclusion is that the higher the norm under consideration, the faster is the time decay.
\end{remark}

Let us mention that this model shares many common features with the nonlinear Schr\"{o}dinger equation and the unstable thin film equation, such as the competition between the attractive and repulsive terms. As in our scaling analysis, the balance between reaction(attraction) and diffusion(repulsive) happens precisely for our chosen exponent $\alpha=m+2/n$. In the nonlinear Schr\"{o}dinger equation, Weinstein \cite{wein83} proposed the existence of the critical exponent $\sigma=2/n$ that would separate those equations that only have local solutions from those that do not, see \cite{Mer04}. In the unstable thin film equation, $m=n+2$ is the critical exponent separating equations with possible finite-time blow-up from problems where the solutions are always bounded \cite{hoch93}, a comprehensive discussion of how scaling properties of the equations relate to infinite-time diffusive spreading and finite-time blow-up can be found in \cite{witel04}, see \cite{ber98,ber20} for the subcritical case $m<n+2$ where blow-up is impossible and for the supercritical case $m>n+2$ where the existence of a solution that blows up in finite time.

\subsection{Structure of the paper}

This paper is organized as follows. Section \ref{sec2} prepares some preliminary lemmas. The following sections are devoted to the detailed proof of the existence of weak solutions for the three cases. The point here is to establish the result with
all necessary details: regularized problem, uniform estimates, passing to the limit in the
regularization parameter.

Precisely, in Section \ref{sec3}, a regularized equation is constructed for which the global strong solution exists. Firstly, a key maximal time of local existence criterion for solutions of the regularized problem is established, see Proposition \ref{Prop1}. Then a priori estimates of the local solution have been derived, see Proposition \ref{Prop2} for a detailed study of the regularity properties of the solutions. Furthermore, by Moser iterative method, we prove that the solution is uniformly bounded in $L^\infty$ space for almost any positive $t$, see Proposition \ref{Prop3}.

Section \ref{sec4} displays the global existence of a weak solution to \er{nkpp} by passing the regularized parameter to zero. The main difficulty comes from the nonlocal term $\int_{\R^n} u dx$ and we prove the three cases in which we use the standard arguments relying on the evolution of the second moment of solutions. We also derive the decay rate of global solutions in the critical case and the supercritical case.

Finally, in Section \ref{sec5}, series of numerical experiments are carried out to verify the results of the earlier sections and explore other issues that lie beyond the scope of the analysis. Further problems and open questions for the nonlinear dynamics of \er{nkpp} are also addressed using numerical simulations.

\section{Preliminaries}\label{sec2}
\def\theequation{2.\arabic{equation}}\makeatother
\setcounter{equation}{0}
\def\thetheorem{2.\arabic{theorem}}\makeatother
\setcounter{theorem}{0}
Before showing the global existence, we shall prepare several lemmas which will be used often in the next sections.
\begin{lemma}[\cite{lieb202}]\label{sobolev}
Let $n\ge 3$. Suppose $u\in H^1(\R^n)$. Then $u\in L^{\frac{2n}{n-2}}(\R^n)$ and the following holds:
\begin{align}
     S_n \|u\|_{L^{\frac{2n}{n-2}}(\R^n)}^2\leq \|\nabla u \|_{L^2(\R^n)}^2,
\end{align}
where $S_n$ is defined in \er{Sn}.
\end{lemma}
 \begin{lemma}\label{GNS}
 Let $n\ge 3$, $1<\frac{b}{a}<\frac{2n}{a(n-2)}$ and $\frac{b}{a}<\frac{2}{a}+\frac{2}{n}$. Assume $w\in L_+^1(\R^n)$ and $w^{\frac{1}{a}}\in H^1(\R^n)$ with $a>0$, then
\begin{align}\label{gns}
    \|w\|_{L^\frac{b}{a}(\R^n)}^{\frac{b}{a}} \le S_n^{\frac{-\lambda b}{2}}\|w\|_{L^1(\R^n)}^{\frac{b}{a}(1-\lambda)}\|\nabla w^{\frac{1}{a}}\|_{L^2(\R^n)}^{b\lambda},
\end{align}
where $\lambda=\frac{1/a-1/b}{1/a-\frac{n-2}{2n}}$.
\end{lemma}
\noindent{\it Proof.} We take $u=w^{\frac{1}{a}}$ in Lemma \ref{sobolev} and employ H\"{o}lder inequality with $1<\frac{b}{a}<\frac{2n}{a(n-2)}$ yield
\begin{align*}
    \|w\|_{L^{\frac{b}{a}}(\R^n)} \le \|w\|_{L^1(\R^n)}^{1-\lambda}\|w^{\frac{1}{a}}\|_{L^{\frac{2n}{n-2}}(\R^n)}^{\lambda a} \le S_n^{-\frac{\lambda a}{2}}\|w\|_{L^1(\R^n)}^{1-\lambda}\|\nabla w^{\frac{1}{a}}\|_{L^2(\R^n)}^{\lambda a},
\end{align*}
where $\lambda=\frac{1/a-1/b}{1/a-\frac{n-2}{2n}}$. \quad$\Box$

The following lemma which have been proved in \cite{BCL15} will play an important role in the proof of global existence of solutions to equation \er{nkpp}.
 \begin{lemma}[\cite{BCL15}]\label{2018dcds} (Gagliardo-Nirenberg-Sobolev inequality) Let $n\ge 3$, $p=\frac{2n}{n-2}$, $1\le r<q<p$ and $\frac{q}{r}<\frac{2}{r}+1-\frac{2}{p}$, then for any $w\in H^1(\R^n)$ and $w\in L^r(\R^n)$, it holds
 \begin{align}
     \|w\|_{L^q(\R^n)}^q\le C_0\|\nabla w\|_{L^2(\R^n)}^2+\left(1-\frac{\lambda q}{2}\right)\left( \frac{2 S_n C_0}{\lambda q} \right)^{-\frac{\lambda q}{2-\lambda q}} \|w\|_{L^r(\R^n)}^{\frac{2(1-\lambda)q}{2-\lambda q}},~~n \ge 3,
 \end{align}
where $\lambda=\frac{\frac{1}{r}-\frac{1}{q}}{\frac{1}{r}-\frac{1}{p}}\in(0,1)$ and $S_n$ is given by \er{Sn}.
\end{lemma}

\begin{lemma}[\cite{BL14}]\label{ode}
Assume $y(t)\ge 0$ is a $C^1$ function for $t>0$ satisfying
\begin{align*}
    y'(t)\le \eta-\beta y(t)^a
\end{align*}
for $\eta>0$, $\beta>0$, then
\begin{itemize}
\item[(i)]~For $a>1$, $y(t)$ has the following hyper-contractive property
\begin{align}
   y(t)\le (\eta/\beta)^{\frac{1}{a}}+\left( \frac{1}{\beta(a-1)t} \right)^{\frac{1}{a-1}}, ~~ \mbox{for any} ~~ t>0.
\end{align}
In addition, if $y(0)$ is bounded, then
\begin{align}
    y(t)\le \max\left( y(0),(\eta/\beta)^{\frac{1}{a}} \right).
\end{align}
\item[(ii)]~For $a=1,$ $y(t)$ decays exponentially
\begin{align}
  y(t) \le \eta/\beta+y(0)e^{-\beta t}.
\end{align}
\end{itemize}
\end{lemma}
More generally, we have
 \begin{lemma}[\cite{BL13}]\label{ODE}
Assume $f(t)\geq 0$ is a non-increasing function for $t>0$. $y(t)\geq0$ is a $C^1$ function and satisfies
$$y'(t)\leq f(t)-\beta y(t)^a$$
for $a>1, \beta>0$. Then for any $t_0>0$ one has
\begin{align}
    y(t)\leq \left( \frac{f(t_0)}{\beta} \right)^{1/a}+\left(\frac{1}{\beta(a-1)(t-t_0)}\right)^{\frac{1}{a-1}},~~\mbox{for any}~~t>t_0.
\end{align}
\end{lemma}

\section{Regularized problem}\label{sec3}
\def\theequation{3.\arabic{equation}}\makeatother
\setcounter{equation}{0}
\def\thetheorem{3.\arabic{theorem}}\makeatother
\setcounter{theorem}{0}

In order to justify the formal arguments of the priori estimates (which will be given in Proposition \ref{Prop2}), we consider the regularized problem
\begin{align}\label{fkppepsilon}
\left\{
  \begin{array}{ll}
  \frac{\partial u_\varepsilon(x,t)}{\partial t}=\Delta u_\varepsilon^m+\varepsilon \Delta u_\varepsilon   + \chi u_\varepsilon^{\alpha} \left( M_0-\int_{\R^n} u_\varepsilon dx \right), ~~& x\in \R^n,~t>0, \\
  u_\varepsilon(x,0)=u_{0\varepsilon}(x),~~& x \in \R^n.
  \end{array}
\right.
\end{align}
Here we define the convolution $u_{0\varepsilon}=J_{\varepsilon}\ast u_0$ where the regularizing kernel $J_\varepsilon=\frac{1}{\varepsilon^n}J\left( \frac{x}{\varepsilon} \right)$ with $J \in C_0^\infty(\R^n)$ and $\int_{\R^n}J dx=1$ so that $\int_{\R^n} J_\varepsilon dx=1.$ $u_{0\varepsilon}$ satisfies $\|u_{0\varepsilon}\|_{L^1(\R^n)} <M_0$ and there exists $\delta>0$ such that for all $0<\varepsilon<\delta$
\begin{align}\label{u0epsilon}
\left\{
  \begin{array}{ll}
   (i)~u_{0\varepsilon}\in L^q(\R^n) ~\mbox{and}~ \|u_{0\varepsilon}\|_{L^q(\R^n)} \le \|u_0\|_{L^q(\R^n)}~ \mbox{for all}~ 1 \le q \le \infty, \\
   (ii)~0\le u_{0\varepsilon}\in L^1 \cap W^{2,q}(\R^n)~ \mbox{for all}~q \in [1,n+3], \\
   (iii)~u_{0\varepsilon} \to u_0~\mbox{strongly in}~L^q(\R^n)~\mbox{as}~\varepsilon \to 0,~\mbox{for some}~q\in [1,\infty), \\
   (iv)~\left\|\nabla u_{0\varepsilon}^m \right\|_{L^2(\R^n)} \le \left\|\nabla u_0^m \right\|_{L^2(\R^n)}, \\[1mm]
   (v)\int_{\R^n} |x|^2 u_{0\varepsilon} dx \to \int_{\R^n} |x|^2 u_{0} dx ~\mbox{as}~\varepsilon \to 0. \\
  \end{array}
\right.
\end{align}
We denote $Q_T=\R^n \times [0,T)$ and
\begin{align}\label{WT}
&W^{2,1}_q(Q_T):=\left\{ u \in L^q(0,T;W^{2,q}(\R^n))~\mbox{and}~u_t \in L^q(0,T;L^q(\R^n))\right\}, \\
&W(Q_T)=W^{2,1}_{\frac{n}{n-1}} \cap W^{2,1}_{n+3}(Q_T).
\end{align}
This section aims to prove the time global solution of \er{fkppepsilon} which reads:

\begin{theorem}\label{global}(Time global strong solution)
Let $n \ge 3, \alpha\ge 1, m>1.$ Suppose that $u_{0\varepsilon}$ satisfies \er{u0epsilon}, then \er{fkppepsilon} has the unique strong solution in $W(Q_T)$ for all $T>0.$
\end{theorem}

For the proof of Theorem \ref{global}, it suffices to show the following three propositions: Proposition \ref{Prop1}, Proposition \ref{Prop2}, Proposition \ref{Prop3}. We first establish the local existence and blow-up criteria, then we show that the local solution admits the uniformly boundedness for extension in time.

\begin{proposition}\label{Prop1}(Time local existence and blow-up criteria)  Let $n \ge 3, \alpha \ge 1, m>1.$ Suppose that $u_{0\varepsilon}$ satisfies \er{u0epsilon}, then there exists a number $T_{\max}=T\left( \|u_{0\varepsilon}\|_{W^{2,n+2}(\R^n)},m,\alpha,n,\chi\right)>0$ such that $u_\varepsilon(x,t) \in W^{2,1}_{n+2}(Q_T) \cap L^\infty(0,T;L^1 \cap L^\infty(\R^n))$ is a non-negative strong solution of \er{fkppepsilon}. Furthermore, if  $T_{\max}<\infty,$ then we have $\displaystyle \limsup_{t \to T_{\max}} \|u_\varepsilon(\cdot,t)\|_{L^\infty(\R^n)}=\infty.$
\end{proposition}
\noindent{\it Proof.} The proof will be carried out as follows. The first route we shall follow here is to consider the problem
\begin{align}\label{fkpp1}
\left\{
  \begin{array}{ll}
  \frac{\partial u_\varepsilon(x,t)}{\partial t}=\nabla \cdot \left( (mu_\varepsilon^{m-1}+\varepsilon)\nabla u_\varepsilon  \right)+\chi u_\varepsilon^\alpha \left( M_0-\int_{\R^n} h dx  \right), ~~& x\in \R^n,~t>0, \\
  u_\varepsilon(x,0)=u_{0\varepsilon}(x),~~& x \in \R^n,
  \end{array}
\right.
\end{align}
where $h \in L^\infty(0,T;L^1(\R^n))$ is a non-negative function. We show the local existence of a strong solution of \er{fkpp1}. Secondly, by the fixed point theorem we prove the local existence of the strong solution of \er{fkppepsilon}. Finally, we state that the local solution satisfies a blow-up criterion and thus close the proof. In the following precise discussions, we will deal with the nonlinear reaction $\alpha>1$ and the linear  reaction $\alpha=1$ respectively.

\noindent{\it\textbf{Step 1}} (Local existence of the non-negative solution of \er{fkpp1})\quad In this step, in order to prove the existence of a strong solution $u_\varepsilon$ in \er{fkpp1}, we observe the equation:
\begin{align}\label{star}
(u_\varepsilon)_t &=\nabla \cdot \left((mf^{m-1}+\varepsilon)\nabla u_\varepsilon\right)+\chi f^{\alpha-1}u_\varepsilon\left(M_0-\int_{\R^n}hdx  \right) \nonumber \\
&=(mf^{m-1}+\varepsilon)\Delta u_\varepsilon +m \nabla f^{m-1}\cdot \nabla u_\varepsilon+\chi f^{\alpha-1}u_\varepsilon \left(M_0-\int_{\R^n}hdx  \right).
\end{align}
Here $\alpha>1.$ The proof is refined in the spirit of \cite{bdp06,2006Sugiyama}. We shall use the notation
\begin{align}\label{XT}
X_T:=\{ & f\in L^\infty(0,T;W^{2,n+2}(\R^n)\cap L^\infty(0,T;L^1 \cap L^\infty(\R^n)), f_t \in L^{n+2}(Q_T): \nonumber \\
& f \ge 0, ~\|f\|_{L^\infty(0,T;W^{2,n+2}(\R^n)}+\|f\|_{L^\infty(0,T;L^1 \cap L^\infty(\R^n))}+\|f_t\|_{L^{n+2}(Q_T)} \nonumber \\
& \qquad~~\le c_1 \|u_{0\varepsilon}\|_{W^{2,n+2}(\R^n)}+c_2 \|u_{0\varepsilon}\|_{L^1 \cap L^\infty (\R^n)}+c_3 \}
\end{align}
for some constants $c_1,c_2,c_3$ only depending on $m,\alpha,n,\chi,M_0.$ By Theorem 9.1 of \cite{LSU} with $f \in X_T$ and $h \in L^\infty(0,T;L^1(\R^n))$, it follows that \er{star} corresponding to the initial data $u_{0\varepsilon}$ has the unique strong solution $u^f_\varepsilon\in W(Q_T).$ Hence we can define a mapping $\Phi$ by
\begin{align}
\Phi:~f \in X_T \mapsto u^f_\varepsilon \in W(Q_T).
\end{align}

Now we claim that $u^f_\varepsilon \ge 0 \in L^\infty(0,T;L^1 \cap L^\infty(\R^n)).$ In the following, we use $u_\varepsilon$ instead of $u^f_\varepsilon$ for simplicity. Multiplying \er{star} by $|u_\varepsilon|^{k-2}u_\varepsilon (k>1)$ yields
\begin{align*}
\frac{1}{k} \frac{d}{dt} \int_{\R^n} |u_\varepsilon|^k dx&=-(k-1)\int_{\R^n} (mf^{m-1}+\varepsilon)|u_\varepsilon|^{k-2}|\nabla u_\varepsilon|^2dx+\chi\int_{\R^n} f^{\alpha-1}|u_\varepsilon|^k dx \left( M_0-\int_{\R^n}h dx \right) \\
&\le \chi M_0 \|f\|_{L^\infty(Q_T)}^{\alpha-1} \int_{\R^n} |u_\varepsilon|^k dx
\end{align*}
which follows that for $k>1$
\begin{align}\label{estimate1}
\frac{d}{dt} \|u_\varepsilon\|_{L^k(\R^n)} \le\chi M_0  \|f\|_{L^\infty(Q_T)}^{\alpha-1}\|u_\varepsilon\|_{L^k(\R^n)}.
\end{align}
Thus we have
\begin{align*}
\|u_\varepsilon\|_{L^k(\R^n)} \le \|u_{0\varepsilon}\|_{L^k(\R^n)}+\chi M_0 \|f\|_{L^\infty(Q_T)}^{\alpha-1} \int_0^t \|u_\varepsilon\|_{L^k(\R^n)} ds.
\end{align*}
Taking $k \to \infty$ and using Gronwall inequality assure that
\begin{align}
\displaystyle \sup_{0<t<T} \|u_\varepsilon(t)\|_{L^\infty(\R^n)} \le \|u_{0\varepsilon}\|_{L^\infty(\R^n)}e^{\chi M_0\|f\|_{L^\infty(Q_T)}^{\alpha-1}~T}.
\end{align}
The nonnegativity of $u_\varepsilon$ can be obtained by multiplying \er{star} with $u_\varepsilon^-:=-\min(u_\varepsilon,0)$ that
\begin{align*}
\frac{1}{2}\frac{d}{dt}\int_{\R^n} |u_\varepsilon^-|^2 dx &=-\int_{\R^n} (mf^{m-1}+\varepsilon)|\nabla u_\varepsilon^-|^2 dx+\chi\int_{\R^n}f^{\alpha-1} |u_\varepsilon^-|^2 dx \left(M_0-\int_{\R^n} h dx \right) \\
& \le \chi M_0\|f\|_{L^\infty(Q_T)}^{\alpha-1}\int_{\R^n} |u_\varepsilon^-|^2 dx.
\end{align*}
It follows
\begin{align*}
\displaystyle \sup_{0<t<T} \|u_\varepsilon^-(\cdot,t)\|_{L^2(\R^n)} \le e^{\chi M_0  \|f\|_{L^\infty(Q_T)}^{\alpha-1} T }~\|u_{0\varepsilon}^-(\cdot,0)\|_{L^2(\R^n)}=0
\end{align*}
which guarantees that for all $0 \le t <T$
\begin{align}
u_\varepsilon(x,t) \ge 0,~~a.e.~~x \in \R^n.
\end{align}
This allows us to integrate \er{star} over $\R^n$
\begin{align*}
\frac{d}{dt}\int_{\R^n} u_\varepsilon dx=\chi\int_{\R^n}f^{\alpha-1}u_\varepsilon dx \left(M_0-\int_{\R^n} h dx \right)
\end{align*}
to obtain
\begin{align}
\displaystyle \sup_{0<t<T} \|u_\varepsilon(t)\|_{L^1(\R^n)} \le \|u_{0\varepsilon}\|_{L^1(\R^n)}e^{\chi M_0 \|f\|_{L^\infty(Q_T)}^{\alpha-1}~T}.
\end{align}
Now we can see that there exists $T_\ast=T_\ast(\varepsilon,\|h\|_{L^\infty(0,T;L^1(\R^n))},\|u_{0\varepsilon}\|_{W^{2,n+2}(\R^n)},\|u_{0\varepsilon}\|_{L^1 \cap L^\infty(\R^n)}, T)$ such that $\Phi$ maps $X_{T_\ast}$ into itself.

Considering the complete metric space $(X_T,d)$ where $d$ is defined by $d(f_1-f_2)=\|f_1-f_2\|_{L^\infty(0,T;L^{n+2}(\R^n))}$, we denote
\begin{align*}
u_1=u_\varepsilon^{f_1},~~u_2=u_\varepsilon^{f_2},~~w=u_1-u_2,
\end{align*}
from \er{star} one has
\begin{align}
(u_1-u_2)_t
= &\nabla \cdot \left( m( f_1^{m-1}-f_2^{m-1} )\nabla u_1+(mf_2^{m-1}+\varepsilon)\nabla(u_1-u_2)\right) \nonumber \\
 & +  \chi\left(( f_1^{\alpha-1}-f_2^{\alpha-1} )u_1+f_2^{\alpha-1}(u_1-u_2) \right)\left( M_0-\int_{\R^n} h dx \right).\label{stars}
\end{align}
The multiplication \er{stars} by $|w|^n w$ gives rise to
\begin{align}\label{starstar}
& \frac{1}{n+2}\frac{d}{dt}\int_{\R^n}|w|^{n+2} dx=-(n+1)\int_{\R^n} (mf_2^{m-1}+\varepsilon)|\nabla (u_1-u_2)|^2 |u_1-u_2|^n dx \nonumber \\
& -(n+1)m \int_{\R^n} (f_1^{m-1}-f_2^{m-1})|u_1-u_2|^n \nabla u_1 \cdot \nabla(u_1-u_2) dx \nonumber \\
& +\chi \int_{\R^n} (f_1^{\alpha-1}-f_2^{\alpha-1})u_1(u_1-u_2)|u_1-u_2|^n dx \left(M_0-\int_{\R^n} h dx\right) \nonumber \\
& +\chi\int_{\R^n} f_2^{\alpha-1}|u_1-u_2|^{n+2} dx \left(M_0-\int_{\R^n} h dx  \right) \nonumber \\
& :=I_1+I_2+I_3+I_4.
\end{align}
By Young's inequality we learn that
\begin{align}\label{I2}
I_2 &\le C(m,n,\|f\|_{L^\infty(Q_T)}) \int_{\R^n} |f_1-f_2|^{\min(m-1,1)} |u_1-u_2|^n |\nabla u_1| ~|\nabla (u_1-u_2)| dx  \nonumber \\
&\le m(n+1)\varepsilon \int_{\R^n} |\nabla (u_1-u_2)|^2 |u_1-u_2|^n dx \nonumber \\
&~~ + C\left(\frac{1}{\varepsilon}, m,\|f\|_{L^\infty(Q_T)}\right) \int_{\R^n} |f_1-f_2|^{2\min(m-1,1)} |\nabla u_1|^2 |u_1-u_2|^n dx \nonumber\\
&\le m(n+1)\varepsilon\int_{\R^n} |\nabla(u_1-u_2)|^2|u_1-u_2|^n dx
\nonumber \\
&~~ +C\left(\frac{1}{\varepsilon},\|f\|_{L^\infty(Q_T)},\|u\|_{L^\infty(0,T;W^{2,n+2}(\R^n))},m\right)\|u_1-u_2\|_{L^{n+2}(\R^n)}^n\|f_1-f_2\|_{L^{n+2}(\R^n)}^2,
\end{align}
where the last inequality is given by H\"{o}lder inequality. Again by virtue of H\"{o}lder inequality we observe
\begin{align}\label{I3}
& I_3
 \le \chi(M_0+\|h\|_{L^\infty(0,T;L^1(\R^n))}) C(\alpha,\|f\|_{L^\infty(Q_T)}) \int_{\R^n} |f_1-f_2|^{\min(\alpha-1,1)} u_1 |u_1-u_2|^{n+1} dx \nonumber\\
& \le \chi(M_0+\|h\|_{L^\infty(0,T;L^1(\R^n))}) C\left(\alpha,\|f\|_{L^\infty(Q_T)},\|u\|_{L^\infty(Q_T)}\right) \|f_1-f_2\|_{L^{n+2}(\R^n)} \|u_1-u_2\|_{L^{n+2}(\R^n)}^{n+1}.
\end{align}
Substituting \er{I2} and \er{I3} into \er{starstar} arrives at
\begin{align*}
\frac{d}{dt} \|w\|_{L^{n+2}(\R^n)}^2 \le c_1 \|f_1-f_2\|_{L^{n+2}(\R^n)}^2+c_1 \|u_1-u_2\|_{L^{n+2}(\R^n)}^2,
\end{align*}
where $c_1$ are constants depending on $\varepsilon, \alpha,m,n, \chi, \|u_{0\varepsilon}\|_{W^{2,n+2}(\R^n)}, \|u_{0\varepsilon}\|_{L^1\cap L^\infty(\R^n)},\|h\|_{L^\infty(0,T;L^1(\R^n))}$. By Gronwall inequality it holds that
\begin{align*}
\displaystyle \sup_{0<t<T_\ast} \|w\|_{L^{n+2}(\R^n)}^2 \le c_1 \|f_1-f_2\|_{L^2(0,T_\ast;L^{n+2}(\R^n))}^2 e^{c_1 T_\ast}.
\end{align*}
Therefore, there exists $T_1=T_1(c_1) \le T_\ast$ such that
\begin{align*}
\displaystyle \sup_{0<t<T_1} \|w\|_{L^{n+2}(\R^n)}^2 \le \frac{1}{2} \|f_1-f_2\|_{L^\infty(0,T_1;L^{n+2}(\R^n))}.
\end{align*}
We find that $\Phi$ becomes a contraction from $X_{T_1}$ into $X_{T_1}$ which is achieved by Banach fixed point theorem. Consequently, $\Phi$ has a fixed point $f=\Phi(f)=u_\varepsilon^f \in X_{T_1}$. Hence, there exists $T_1=T_1\left(\varepsilon, \alpha,m,n,\chi,\|u_{0\varepsilon}\|_{W^{2,n+2}(\R^n)}, \|u_{0\varepsilon}\|_{L^1\cap L^\infty(\R^n)},\|h\|_{L^\infty(0,T;L^1(\R^n))} \right)$ such that there is a desired strong solution $u_\varepsilon^f$ of \er{fkpp1} on $[0,T_1]$ corresponding to the initial data $u_{0\varepsilon}$.

The proof for the case $\alpha=1$ is a word for word translation of the proof for $\alpha>1$ except $\alpha-1$ is replaced by zero.

\noindent{\it\textbf{Step 2}} (Local existence of the non-negative solution of \er{fkppepsilon}) \quad We firstly claim that the solution $u_\varepsilon$ is bounded in $L^\infty(0,T_2;L^\infty(\R^n))$ ($T_2$ is to be determined) as a consequence of the following computations:
\begin{align} \label{AB}
& \frac{1}{k} \frac{d}{dt} \int_{\R^n} u_\varepsilon^k dx=-(k-1)\int_{\R^n} (m u_\varepsilon^{m-1}+\varepsilon) u_\varepsilon^{k-2}|\nabla u_\varepsilon|^2 dx + \chi\int_{\R^n}  u_\varepsilon^{k+\alpha-1} dx \left(M_0-\int_{\R^n} h dx \right) \nonumber \\
\le &-\frac{4m(k-1)}{(m+k-1)^2} \int_{\R^n} |\nabla u^{\frac{m+k-1}{2}}|^2 dx + \chi\left( M_0+\|h\|_{L^\infty(0,T_1;L^1(\R^n))} \right)  \int_{\R^n} u_\varepsilon^{k+\alpha-1} dx \nonumber \\
= & -\frac{4m(k-1)}{(m+k-1)^2} \int_{\R^n} |\nabla u^{\frac{m+k-1}{2}}|^2 dx + \bar{A} \int_{\R^n} u_\varepsilon^{k+\alpha-1} dx,
\end{align}
where $\bar{A}=\chi\left( M_0+\|h\|_{L^\infty(0,T_1;L^1(\R^n))} \right)$.

For $\alpha>1$, we apply
\begin{align*}
w=u_\varepsilon^{\frac{m+k-1}{2}},~q=\frac{2(k+\alpha-1)}{k+m-1},~r=\frac{2k}{k+m-1},~C_0=\frac{2m(k-1)}{\bar{A} (k+m-1)^2}
\end{align*}
in Lemma \ref{2018dcds} for $k>\max\left(\frac{n(\alpha-m)}{2},1\right)$ such that
\begin{align*}
 &\bar{A} \|u_\varepsilon\|_{L^{k+\alpha-1}(\R^n)}^{k+\alpha-1} \\
 \le & \frac{2m(k-1)}{(k+m-1)^2} \int_{\R^n} |\nabla u_\varepsilon^{\frac{m+k-1}{2}}|^2 dx+C(n,m,\bar{A}) \left( \frac{(k+m-1)^2}{k-1} \right)^{\frac{\lambda q}{2-\lambda q}} \|u_\varepsilon\|_{L^k(\R^n)}^{\frac{(m+k-1)(1-\lambda)q}{2-\lambda q}} \\
= &\frac{2m(k-1)}{(k+m-1)^2} \int_{\R^n} |\nabla u_\varepsilon^{\frac{m+k-1}{2}}|^2 dx+C(n,m,\bar{A}) \left( \frac{(k+m-1)^2}{k-1} \right)^{\frac{n(\alpha-1)}{2k+n(m-\alpha)}} \|u_\varepsilon\|_{L^k(\R^n)}^{\frac{k+\alpha-1+\frac{n(m-\alpha)}{2}}{\frac{n(m-\alpha)}{2k}+1}},
\end{align*}
where $\lambda=\frac{\frac{1}{r}-\frac{1}{q}}{\frac{1}{r}-\frac{n-2}{2n}}$. Plugging it into \er{AB} we compute
\begin{align*}
\|u_\varepsilon\|_{L^k(\R^n)} \le \|u_{0\varepsilon}\|_{L^k(\R^n)}+C(n,m,\bar{A}) \left( \frac{(k+m-1)^2}{k-1} \right)^{\frac{n(\alpha-1)}{2k+n(m-\alpha)}} \int_0^t \|u_\varepsilon(s)\|_{L^k(\R^n)}^{\frac{2(\alpha-1)}{2+\frac{n(m-\alpha)}{k}}+1}ds.
\end{align*}
Taking $k \to \infty$ gives
\begin{align*}
\|u_\varepsilon\|_{L^\infty(\R^n)} \le \|u_{0\varepsilon}\|_{L^\infty(\R^n)}+C(n,m,\bar{A}) \int_0^t \|u_\varepsilon(s)\|_{L^\infty(\R^n)}^{\alpha} ds
\end{align*}
which provides the estimate
\begin{align}\label{A222}
\|u_\varepsilon\|_{L^\infty(\R^n)} \le \left( \frac{1}{C\left(n,m,\bar{A}\right)(\alpha-1)(\overline{T}-t)}    \right)^{\frac{1}{\alpha-1}},\quad \overline{T}=\frac{\|u_{0\varepsilon}\|_{L^\infty(\R^n)}^{1-\alpha}}{C\left(n,m,\bar{A}\right)(\alpha-1)}.
\end{align}
Hence $u_\varepsilon$ is bounded in $L^\infty(0,T_2;L^\infty(\R^n))$ where $T_2=\frac{\overline{T}}{2}$.

For $\alpha=1,$ \er{AB} directly follows that for any $t>0$ and $k>1$
\begin{align*}
\|u_\varepsilon\|_{L^k(\R^n)} \le \|u_{0\varepsilon}\|_{L^k(\R^n)} e^{\bar{A} t}.
\end{align*}
Letting $k \to \infty$ we get
\begin{align}\label{A111}
\|u_\varepsilon\|_{L^\infty(\R^n)} \le \|u_{0\varepsilon}\|_{L^\infty(\R^n)} e^{\bar{A} t}.
\end{align}

On the other hand, integrating \er{fkpp1} over $\R^n$ we obtain that for $\alpha \ge 1$
\begin{align*}
\frac{d}{dt}\int_{\R^n} u_\varepsilon dx  =\chi \int_{\R^n}  u_\varepsilon^\alpha dx \left( M_0-\int_{\R^n} h dx \right)
 \le \chi M_0 \|u_\varepsilon\|_{L^\infty(\R^n)}^{\alpha-1} \int_{\R^n}  u_\varepsilon dx.
\end{align*}
By \er{A222} and \er{A111} we prove that $u_\varepsilon \in L^\infty(0,T_2;L^1(\R^n)).$

Next, for $h \in X_{T_1},$ consider the map
\begin{align*}
F:~h \in L^{\infty}(0,T_1;L^1(\R^n)) \mapsto u_\varepsilon^h \in L^{\infty}(0,T_1;L^1(\R^n)).
\end{align*}
Analogous to Step 1, it can be seen that there exists $T_3=T_3(M_0,\alpha,m,n,\chi,\|u_{0\varepsilon}\|_{L^1 \cap L^\infty(\R^n)}) \le \min(T_2,T_1)$ such that
$F$ is a contraction from $L^{\infty}(0,T_3;L^1(\R^n))$ to $L^{\infty}(0,T_3;L^1(\R^n))$ by making use of Banach's fixed point theorem. Thus $F$ has a fixed point $h=F(h)=u_\varepsilon^h \in L^{\infty}(0,T_3;L^1(\R^n))$ and we prove the existence of a solution of \er{fkppepsilon} on the time interval $[0,T_3].$ This is exactly the anticipated result and we complete the proof of Proposition \ref{Prop1}. \quad$\Box$

\begin{proposition}\label{Prop2}
(Priori estimates in $L^k$ for $1<k<\infty$) \quad Let the same assumption as that in Proposition \ref{Prop1} hold. Suppose that $u_\varepsilon$ is the non-negative strong solution of \er{fkppepsilon}, $C$ is a positive constant depending on $\|u_{0\varepsilon}\|_{L^k(\R^n)},k,m,\alpha,n,\chi,m_0,M_0$ but not on $\varepsilon,$ then $u_\varepsilon$ satisfies the following estimates:
\begin{enumerate}
\item[\textbf{(i)}] For $1\le \alpha<m+2/n,$ the following holds true that for any $t>0$,
\begin{align}
  \|u_\varepsilon(\cdot,t) \|_{L^k(\R^n)} \le   C,~~\mbox{for all}~k \in (1,\infty).
\end{align}
\item[\textbf{(ii)}] For $\alpha=m+2/n,$ if
\begin{align*}
   M_0 \le\left( \frac{S_n(\alpha-m)}{\chi } \right)^{\frac{1}{\alpha-m+1}}\frac{\alpha-m+1}{\alpha-m},
\end{align*}
then $u_\varepsilon$ satisfies that for any $t>0$
\begin{align}
     \|u_\varepsilon(\cdot,t)\|_{L^k(\R^n)} \le C\left( 1+ t \right)^{-\frac{k-1}{k(m+2/n-1)}},~~\mbox{for all}~~k\in (1,\infty),
\end{align}
where $S_n$ is given by \er{Sn}.
\item[\textbf{(iii)}] For $\alpha>m+2/n$, $p_0=\frac{n(\alpha-m)}{2}$. Now we assume
\begin{align}\label{p0}
\|u_{0\varepsilon}\|_{L^{\frac{n(\alpha-m)}{2}}(\R^n)}<C_{p_0},
\end{align}
where $C_{p_0}$ is defined as \er{cp0}. Then $u_\varepsilon$ has the following decay property:
 \begin{align}
     \|u_\varepsilon(\cdot,t)\|_{L^k(\R^n)} \le C\left( 1+ t \right)^{-\frac{k-1}{k(\alpha-1)}},~~\mbox{for all}~~k\in (1,\infty).
\end{align}
\end{enumerate}
Furthermore, the following regularities hold true for any $T>0$
\begin{align*}
& u_\varepsilon \in L^\infty \left(0,T; L^k(\R^n) \right), \\
& u_\varepsilon \in L^{k+\alpha-1}\left(0,T;L^{k+\alpha-1}(\R^n)\right), \\
& \nabla u_\varepsilon^{\frac{m+k-1}{2}} \in L^2\left(0,T;L^2(\R^n) \right).
\end{align*}
\end{proposition}

For the proof of proposition \ref{Prop2}, it suffices to show the following three lemmas. For simplicity in presentation, throughout this section, we omit all the $\varepsilon$ dependents and use $u$ instead of $u_\varepsilon$. We denote $C$ by a positive constant depending not only on $m,n,\alpha,\chi$, but also on other associated quantities (we will show them clearly at different occurrences in $C(\cdot)$). Most of the prior estimates are based on the following arguments. Multiplying \er{nkpp} by $ku^{k-1}(k>1)$ we obtain
\begin{align}\label{fsk}
    & \frac{d}{dt}\int_{\R^n} u^k dx +\frac{4mk(k-1)}{(k+m-1)^2}\int_{\R^n} |\nabla u^{\frac{k+m-1}{2}}|^2 dx +k\chi\int_{\R^n} udx \int_{\R^n} u^{k+\alpha-1} dx\nonumber\\
    = & k\chi M_0\int_{\R^n} u^{k+\alpha-1}dx.
\end{align}
\begin{lemma}\label{lemma1}
(Case of $1\le\alpha<m+2/n$) Let the same assumptions as that in proposition \ref{Prop2} hold. Then there exists a positive constant $C$ depending on $\|u_{0\varepsilon}\|_{L^k(\R^n)},k,$ $m$, $\alpha$, $n$, $\chi, m_0,M_0$ but not on $\varepsilon$.  The following holds true that for any $t>0$
\begin{align}
    \|u_\varepsilon\|_{L^k(\R^n)} \le  C,~~\mbox{for all}~~k\in (1,\infty).
\end{align}
\end{lemma}
\noindent{\it Proof.}
Using
\begin{align*}
    w^{\frac{1}{a}}=u^{\frac{k+m-1}{2}},~~ b=\frac{2(k+\alpha-1)}{k+m-1},~~ a=\frac{2k'}{k+m-1}
\end{align*}
in Lemma \ref{GNS} for $k>\max\left\{\frac{n(\alpha-m)}{2}-(\alpha-1),1 \right\}$ and $\max\left\{\frac{n(\alpha-m)}{2},1 \right\}<k'<k+\alpha-1$, it holds that
\begin{align}\label{gj1}
    \|u\|_{L^{k+\alpha-1}(\R^n)}^{k+\alpha-1} \le S_n^{-\frac{\lambda(k+\alpha-1)}{k+m-1}} \|\nabla u^{\frac{k+m-1}{2}}\|_{L^2(\R^n)}^{\frac{2\lambda(k+\alpha-1)}{k+m-1}} \|u^{\frac{k+m-1}{2}}\|_{L^{\frac{2k'}{k+m-1}}(\R^n)}^{\frac{2(1-\lambda)(k+\alpha-1)}{k+m-1}},
\end{align}
where $\lambda=\frac{\frac{k+m-1}{2k'}-\frac{k+m-1}{2(k+\alpha-1)}}{\frac{k+m-1}{2k'}-\frac{n-2}{2n}}$. We further use Young's inequality to get
\begin{align}\label{*}
    \|u\|_{L^{k+\alpha-1}(\R^n)}^{k+\alpha-1} \le \frac{m(k-1)}{\chi M_0(k+m-1)^2}\|\nabla u^{\frac{k+m-1}{2}}\|_{L^2(\R^n)}^2+C(k,M_0)\|u\|_{L^{k'}(\R^n)}^{\frac{(1-\lambda)(k+\alpha-1)}{1-\frac{\lambda(k+\alpha-1)}{k+m-1}}}
\end{align}
since the choices of $k,k'$ ensure
\begin{align*}
    \frac{2\lambda(k+\alpha-1)}{k+m-1}<2.
\end{align*}
Therefore, by H\"{o}lder inequality with $1<k'<k+\alpha-1$  and plugging \er{*} into \er{fsk} one has
\begin{align}
     & \frac{d}{dt}\int_{\R^n} u^k dx  +\frac{3mk(k-1)}{(k+m-1)^2}\int_{\R^n} |\nabla u^{\frac{k+m-1}{2}}|^2 dx +k\chi m_0 \int_{\R^n} u^{k+\alpha-1} dx \nonumber\\
     \le  & C(k,M_0)\|u\|_{L^{k'}(\R^n)}^r \nonumber\\
      \le & C(k,M_0) \|u\|_{L^{k+\alpha-1}(\R^n)}^{r\theta}\|u\|_{L^1(\R^n)}^{r(1-\theta)} \nonumber\\
     \le & C(k,M_0)\|u\|_{L^{k+\alpha-1}(\R^n)}^{r\theta},\label{gj2}
\end{align}
where $r=\frac{(1-\lambda)(k+\alpha-1)}{1-\frac{\lambda(k+\alpha-1)}{k+m-1}}$, $\theta=\frac{(k'-1)(k+\alpha-1)}{k'(k+\alpha-2)}$. Moreover, a tedious calculation assures
\begin{align}
r\theta<k+\alpha-1
\end{align}
if and only if
\begin{align}
1\le \alpha<m+2/n.
\end{align}
An immediate application of Young's inequality in \er{gj2} leads to
\begin{align}\label{gj3}
    \frac{d}{dt}\int_{\R^n} u^k dx  +\frac{3mk(k-1)}{(k+m-1)^2}\int_{\R^n} |\nabla u^{\frac{k+m-1}{2}}|^2 dx +\frac{km_0\chi}{2} \int_{\R^n} u^{k+\alpha-1} dx \le C(k,m_0,M_0).
\end{align}
Recalling the fact that $m(t) \le M_0$, by H\"{o}lder inequality we have
\begin{align}\label{interpolation}
    \left( \|u\|_{L^k(\R^n)}^k\right)^{\frac{k+\alpha-2}{k-1}}  \le \|u\|_{L^{k+\alpha-1}(\R^n)}^{k+\alpha-1}\|u\|_{L^1(\R^n)}^{\frac{\alpha-1}{k-1}} \le \|u\|_{L^{k+\alpha-1}(\R^n)}^{k+\alpha-1} M_0^{\frac{\alpha-1}{k-1}}.
\end{align}
Hence \er{gj3} is equivalent to
\begin{align}
    \frac{d}{dt}\int_{\R^n} u^kdx+\frac{km_0\chi}{2M_0^{\frac{\alpha-1}{k-1}}}\left( \|u\|_{L^k(\R^n)}^k \right)^{\frac{k+\alpha-2}{k-1}} \le  C(k,m_0,M_0).
\end{align}
Setting
\begin{align}
    y(t)=\int_{\R^n} u^kdx,~~ a=\frac{k+\alpha-2}{k-1},~~ \eta=C(k,m_0,M_0)
\end{align}
in Lemma \ref{ode} one has that for any $1<k<\infty$
\begin{align}
    \|u(\cdot,t)\|_{L^k(\R^n)} \le C, ~~\mbox{for any } t>0,
\end{align}
where $C$ is a constant depending on $\|u_{0\varepsilon}\|_{L^k(\R^n)}$, $k, m_0, M_0, \alpha,m,n,\chi$.

On the other hand, we integrate \er{gj3} from $0$ to $T$ in time to obtain that for any $T>0$
\begin{align}
    &\int_{\R^n} u^k(T)dx +\frac{3mk(k-1)}{(k+m-1)^2}\int_0^T\int_{\R^n} |\nabla u^{\frac{k+m-1}{2}}|^2dxdt +\frac{k\chi m_0}{2}\int_0^T\int_{\R^n} u^{k+\alpha-1} dxdt \nonumber \\
    \le & \int_{\R^n} u_{0\varepsilon}^kdx +C(k,m_0,M_0)T,
\end{align}
from which we derive that for any $T>0$ and $1<k<\infty$
\begin{align}
   u\in L^{k+\alpha-1}\left(0,T;L^{k+\alpha-1}(\R^n)\right),~ \nabla u^{\frac{k+m-1}{2}} \in L^2\left(0,T;L^2(\R^n)\right).
\end{align}
Thus we complete the proof of Lemma \ref{lemma1}.\quad$\Box$

The next goal is to consider the critical exponent case $\alpha=m+2/n.$

\begin{lemma}\label{lemma2}
(Case of $\alpha=m+2/n$) Let the same assumptions as that in proposition \ref{Prop2} hold. If the total mass $M_0$ satisfies
\begin{align*}
    M_0 \le \left( \frac{S_n(\alpha-m)}{\chi} \right)^{\frac{1}{\alpha-m+1}}\frac{\alpha-m+1}{\alpha-m},
\end{align*}
then there exists a positive constant $C$ depending on $\|u_{0\varepsilon}\|_{L^k(\R^n)}, k$, $m$, $\alpha$, $n$, $\chi, m_0,M_0$ but not on $\varepsilon$ such that $u_\varepsilon$ satisfies
\begin{align}
     \|u_\varepsilon\|_{L^k(\R^n)} \le C\left( 1+ t \right)^{-\frac{k-1}{k(m+2/n-1)}},~~\mbox{for all}~~k\in (1,\infty).
\end{align}
\end{lemma}
\noindent{\it Proof.} For any $k>\frac{n(\alpha-m)}{2}-(\alpha-1)$, keeping the fact $\alpha=m+2/n$ in mind we get the following estimate
 \begin{align}
     ||u||_{L^{k+\alpha-1}(\R^n)}^{k+\alpha-1} &\leq  \frac{1}{S_n} {||\nabla u^{\frac{m+k-1}{2}}||}_{L^2{(\R^n)}}^2{||u||}_{L^1(\R^n)}^{\alpha-m}\label{Inequa}
\end{align}
by Lemma \ref{GNS}. Combining \er{fsk} with \er{Inequa} we obtain
\begin{align}\label{339}
     \frac{d}{dt}\int_{\R^n} u^kdx + k\left( \frac{4mS_n(k-1)}{(k+m-1)^2}\|u\|_{L^1(\R^n)}^{-(\alpha-m)}+ \chi\|u\|_{L^1(\R^n)}-\chi M_0 \right) \int_{\R^n} u^{k+\alpha-1} dx \le 0.
 \end{align}
Denote
\begin{align}\label{denote}
     y=\|u\|_{L^1(\R^n)},~~
    \gamma(k)=\frac{4mS_n(k-1)}{(k+m-1)^2},~~f(y)=\gamma(k)y^{-(\alpha-m)}+\chi y-\chi M_0.
\end{align}
A straightforward calculation shows that at
\begin{align}
    y_0=\left( \frac{\gamma(k)(\alpha-m)}{\chi} \right)^{\frac{1}{\alpha-m+1}},
\end{align}
$f(y)$ attains its minimum
\begin{align}
    f(y_0)=\chi^{\frac{\alpha-m}{\alpha-m+1}}\gamma(k)^{\frac{1}{\alpha-m+1}}(\alpha-m)^{\frac{1}{\alpha-m+1}}\frac{\alpha-m+1}{\alpha-m}-\chi M_0.
\end{align}
It's easy to verify
 \begin{align}
     f(y)>f(y_0)\ge 0,~~\mbox{for any}~~ y>0,
 \end{align}
whenever
\begin{align}\label{chi}
   M_0\le  \left( \frac{\gamma(k)(\alpha-m)}{\chi } \right)^{\frac{1}{\alpha-m+1}}\frac{\alpha-m+1}{\alpha-m}.
\end{align}
We point out that $\gamma(k)$ attains its maximum when $k=m+1$. Therefore, we first treat the case of $k=m+1$.\\
\noindent{{\it\textbf{Step 1}}} (Decay estimate in $\|u\|_{L^{m+1}(\R^n)}$) \quad Taking $k=m+1$ in \er{339} we have
 \begin{align}\label{gjm+1}
     \frac{d}{dt}\int_{\R^n} u^{m+1}dx +(m+1)\left( \chi\|u\|_{L^1(\R^n)}-\chi M_0+S_n\|u\|_{L^1(\R^n)}^{-(\alpha-m)} \right)\|u\|_{L^{\alpha+m}(\R^n)}^{\alpha+m} \le 0.
 \end{align}
Similar arguments from \er{denote} to \er{chi} with $k=m+1$ yield that
\begin{align}
     M_0\le  \left( \frac{S_n(\alpha-m)}{\chi } \right)^{\frac{1}{\alpha-m+1}}\frac{\alpha-m+1}{\alpha-m},
\end{align}
results in
\begin{align*}
    f(\|u\|_{L^1(\R^n)})>\eta \ge 0,
\end{align*}
where
\begin{align*}
    \eta=f\left(\left( \frac{S_n(\alpha-m)}{\chi}\right)^{\frac{1}{\alpha-m+1}}\right)=\chi^{\frac{\alpha-m}{\alpha-m+1}}(S_n(\alpha-m))^{\frac{1}{\alpha-m+1}}\frac{\alpha-m+1}{\alpha-m}-\chi M_0.
\end{align*}
Hence we recover the inequality
\begin{align}\label{re}
    \frac{d}{dt}\int_{\R^n} u^{m+1} dx +(m+1)\eta \|u\|_{L^{\alpha+m}(\R^n)}^{\alpha+m} \le 0.
\end{align}
Using H\"{o}lder inequality with $1<m+1<m+\alpha$ we obtain
\begin{align}\label{holder}
    \left( \|u\|_{L^{m+1}(\R^n)}^{m+1} \right)^{\frac{\alpha+m-1}{m}} \le \|u\|_{L^{\alpha+m}(\R^n)}^{\alpha+m}\|u\|_{L^1(\R^n)}^{\frac{\alpha-1}{m}} \le \|u\|_{L^{\alpha+m}(\R^n)}^{\alpha+m} M_0^{\frac{\alpha-1}{m}}.
\end{align}
Thus taking \er{gjm+1} and \er{holder} together gives
\begin{align}
    \frac{d}{dt}\int_{\R^n} u^{m+1}dx +(m+1)\eta M_0^{-\frac{\alpha-1}{m}}\left( \int_{\R^n} u^{m+1}dx \right)^{\frac{\alpha+m-1}{m}} \le 0
\end{align}
which follows
\begin{align}
     \|u\|_{L^{m+1}(\R^n)}
     & \le \left( \frac{1}{\frac{(m+1)(\alpha-1)\eta}{m}M_0^{-\frac{\alpha-1}{m}}t +\|u_{0\varepsilon}\|_{L^{m+1}(\R^n)}^{-\frac{(m+1)(\alpha-1)}{m}}} \right)^{\frac{m}{(m+1)(\alpha-1)}} \le C(1+t)^{-\frac{m}{(m+1)(\alpha-1)}},\label{estimate2.1}
\end{align}
where C is a constant depending on $m$, $n, \chi$, $M_0$, $\|u_{0\varepsilon}\|_{L^{m+1}(\R^n)}$.

In addition, integrating \er{re} from $0$ to $T$ in time we obtain that for any $T>0$
\begin{align*}
    \int_{\R^n} u^{m+1}(T)dx +(m+1)\eta \int_0^T\int_{\R^n} u^{\alpha+m} dxdt \le \int_{\R^n} u_{0\varepsilon}^{m+1}dx
\end{align*}
which assures that
\begin{align}\label{regular}
    \int_0^T\int_{\R^n} u^{\alpha+m}dxdt  \le C \left( \|u_{0\varepsilon}\|_{L^{m+1}(\R^n)}^{m+1}, M_0 \right).
\end{align}
Therefore, integrating the following equality
\begin{align}\label{m+1}
  &\frac{d}{dt}\int_{\R^n} u^{m+1}dx +(m+1)\int_{\R^n} |\nabla u^m|^2dx +(m+1)\chi\int_{\R^n} udx\int_{\R^n} u^{\alpha+m} dx\nonumber\\
  = & (m+1)\chi M_0\int_{\R^n} u^{\alpha+m}dx.
\end{align}
from $0$ to $T$ in time we also obtain
\begin{align}
   \nabla u^{m} \in L^2(0,T;L^2(\R^n)).
\end{align}

\noindent{\it \textbf{Step 2}} (Decay estimates in $\|u\|_{L^k(\R^n)}$ for $1<k<\infty$) \quad For $1<k<m+1$, by \er{estimate2.1} with H\"{o}lder inequality we have
\begin{align}\label{estimate2.2}
    \|u\|_{L^k(\R^n)} & \le \|u\|_{L^{m+1}(\R^n)}^{(m+1)\frac{k-1}{km}}\|u\|_{L^1(\R^n)}^{\frac{m-k+1}{km}} \le C(1+t)^{-\frac{k-1}{k(\alpha-1)}},
\end{align}
where C is a constant depending on $m$,$n, \chi$, $M_0, k$, $\|u_{0\varepsilon}\|_{L^{m+1}(\R^n)}$.

For $m+1<k<\infty$, taking
\begin{align*}
    w=u^{\frac{k+m-1}{2}},~~q=\frac{2(k+\alpha-1)}{k+m-1},~~r=\frac{2(m+1)}{k+m-1},~~ C_0=\frac{2m(k-1)}{\chi M_0(k+m-1)^2}
\end{align*}
in Lemma \ref{2018dcds} for $k>\max(2-\alpha,1)$, we have
\begin{align}\label{estimate2}
    \|u\|_{L^{k+\alpha-1}(\R^n)}^{k+\alpha-1} \le \frac{2m(k-1)}{\chi M_0(k+m-1)^2}\|\nabla u^{\frac{k+m-1}{2}}\|_{L^2(\R^n)}^2 +C(k,M_0)\|u\|_{L^{m+1}(\R^n)}^{(m+1){\frac{k+\alpha-2}{m}}}.
\end{align}
Recalling \er{estimate2.1} and substituting \er{estimate2} into \er{fsk} one has
\begin{align}
    & \frac{d}{dt}\int_{\R^n} u^kdx+\frac{2mk(k-1)}{(k+m-1)^2}\|\nabla u^{\frac{k+m-1}{2}}\|_{L^2(\R^n)}^2 +k\chi\int_{\R^n} udx\int_{\R^n} u^{k+\alpha-1}dx \nonumber \\
    \le & C(k,M_0)\left( 1+t \right)^{-\frac{k+\alpha-2}{\alpha-1}}\label{estimate2.3}.
\end{align}
Therefore, we have
\begin{align}
    \frac{d}{dt}\int_{\R^n} u^kdx +\frac{k\chi m_0}{M_0^{\frac{\alpha-1}{k-1}}}\left( \int_{\R^n} u^kdx \right)^{\frac{k+\alpha-2}{k-1}}
    \le  C(k,M_0)\left( 1+t \right)^{-\frac{k+\alpha-2}{\alpha-1}}.
\end{align}
Denote
\begin{align*}
    y(t)=\int_{\R^n} u^kdx, ~~ f(t)=C(k,M_0)\left( 1+t \right)^{-\frac{k+\alpha-2}{\alpha-1}},~~ a=1+\frac{\alpha-1}{k-1}
\end{align*}
and take $t_0=t/2$ in Lemma \ref{ODE} we derive that for any $t>0$
\begin{align}
    \|u(\cdot,t)\|_{L^k(\R^n)} \le C \left( 1+ t \right)^{-\frac{k-1}{k(\alpha-1)}},
\end{align}
where $C$ is a constant depending on $\|u_{0\varepsilon}\|_{L^k(\R^n)}$, $k, m$, $n$, $\chi$, $m_0, M_0$.

Next, integrating \er{estimate2.3} from $0$ to $T$ in time we obtain that for any $T>0$
\begin{align*}
    &\int_{\R^n} u^k(T) dx +\frac{2mk(k-1)}{(k+m-1)^2}\int_0^T\int_{\R^n} |\nabla u^{\frac{k+m-1}{2}}|^2 dxdt +k\chi\int_0^T\int_{\R^n} udx \int_{\R^n} u^{k+\alpha-1}dxdt \nonumber\\
    \le  & \int_{\R^n} u_{0\varepsilon}^kdx +C(k,M_0).
\end{align*}
Hence we also obtain the following regularities
\begin{align}
   u\in L^{k+\alpha-1}(0,T;L^{k+\alpha-1}(\R^n)),~ \nabla u^{\frac{k+m-1}{2}} \in L^2(0,T;L^2(\R^n)).
\end{align}
Thus completes the proof of this lemma. \quad $\Box$

We are now in a position to begin the study of the supercritical case $\alpha>m+2/n$.
\begin{lemma}\label{lemma3}
(Case of $\alpha>m+2/n$) Let the same assumptions as that in proposition \ref{Prop2} hold, $p_0=\frac{n(\alpha-m)}{2}$. If we assume
\begin{align*}
    \|u_{0\varepsilon}\|_{L^{\frac{n(\alpha-m)}{2}}(\R^n)}<C_{p_0},
\end{align*}
where $C_{p_0}$ is given by \er{cp0}, then for any $t>0$
\begin{align}
     \|u_\varepsilon\|_{L^k(\R^n)} \le C\left( 1+t \right)^{-\frac{k-1}{k(\alpha-1)}},~~\mbox{for~any}~1<k<\infty,
\end{align}
where $C$ is a positive constant depending on $\|u_{0\varepsilon}\|_{L^k(\R^n)}$, $k, \alpha, m,n, \chi, m_0,M_0$ but not on $\varepsilon$.
\end{lemma}
\noindent{\it Proof.} Firstly we plug
\begin{align*}
    w=u^{\frac{k+m-1}{2}},~~ q=\frac{2(k+\alpha-1)}{k+m-1},~~r=\frac{2k'}{k+m-1},~~ C_0=\frac{4m(k-1)}{\chi M_0(k+m-1)^2}
\end{align*}
into Lemma \ref{2018dcds} for any $k'>\frac{n(\alpha-m)}{2}$ and $k>\max\left(\frac{n(\alpha-m)}{2}-(\alpha-1),1\right)$ to obtain
\begin{align}
  \|u\|_{L^{k+\alpha-1}(\R^n)}^{k+\alpha-1}
 \le \frac{4m(k-1)}{\chi M_0(k+m-1)^2} \left\|\nabla u^{\frac{k+m-1}{2}} \right\|_{L^2(\R^n)}^2 +C_{knm}\left\|u^{\frac{k+m-1}{2}}\right\|_{L^{\frac{2k'}{k+m-1}}(\R^n)}^{\frac{2(1-\lambda)(k+\alpha-1)}{k+m-1}\frac{1}{1-\frac{\lambda(k+\alpha-1)}{k+m-1}}},
\end{align}
where
\begin{align*}
   C_{knm}
   = \left( \frac{8mS_n(k-1)}{\chi M_0\lambda q(k+m-1)^2}\right)^{-\frac{\lambda q}{2-\lambda q}}\frac{2-\lambda q}{2},~~
   \lambda=\frac{\frac{k+m-1}{2k'}-\frac{k+m-1}{2(k+\alpha-1)}}{\frac{k+m-1}{2k'}-\frac{n-2}{2n}}.
\end{align*}
Then substituting the above estimates into \er{fsk} we get
\begin{align}\label{r}
     \frac{d}{dt}\int_{\R^n} u^kdx +k\chi m(t)\int_{\R^n} u^{k+\alpha-1}dx \le  k\chi M_0 C_{knm}\|u\|_{L^{k'}(\R^n)}^b,
\end{align}
where $b={\frac{(1-\lambda)(k+\alpha-1)}{1-\frac{\lambda(k+\alpha-1)}{k+m-1}}}$. Because of $1<\frac{n(\alpha-m)}{2}<k'<k+\alpha-1$, by interpolation inequality we compute
\begin{align}\label{bbb}
    \|u\|_{L^{k'}(\R^n)}^b \le \|u\|_{L^{k+\alpha-1}(\R^n)}^{b\theta}\|u\|_{L^{\frac{n(\alpha-m)}{2}}(\R^n)}^{b(1-\theta)},
\end{align}
where $\frac{1}{k'}=\frac{2(1-\theta)}{n(\alpha-m)}+\frac{\theta}{k+\alpha-1}$. Some calculations yield that for any $\frac{n(\alpha-m)}{2}<k'<k+\alpha-1$.
\begin{align*}
    b\theta=k+\alpha-1, ~~ b(1-\theta)=\frac{n(\alpha-m)}{2}\frac{k+\alpha-1-k'}{k'-\frac{n(\alpha-m)}{2}}.
\end{align*}
We choose
$$ k'=\frac{k+\alpha-1+\frac{n(\alpha-m)}{2}}{2} $$
used in \er{bbb} such that it satisfies
\begin{align}\label{gj3.1}
    \|u\|_{L^{k'}(\R^n)}^b \le \|u\|_{L^{k+\alpha-1}(\R^n)}^{k+\alpha-1}\|u\|_{L^{\frac{n(\alpha-m)}{2}}(\R^n)}^{\frac{n(\alpha-m)}{2}}.
\end{align}
Keeping in mind the fact that $m_0 \le m(t)$ and collecting \er{r} and \er{gj3.1} one obtains that for any $k>\max\left(\frac{n(\alpha-m)}{2}-(\alpha-1),1\right)$
\begin{align}\label{gj4.1}
    &\frac{d}{dt}\int_{\R^n} u^kdx +k\chi m_0C_k^{-p_0}\left( C_k^{p_0}-\|u\|_{L^{p_0}(\R^n)}^{p_0} \right)\|u\|_{L^{k+\alpha-1}(\R^n)}^{k+\alpha-1}\le 0,
\end{align}
where $C_k$ is defined as
\begin{align}\label{Bigstar}
    C_k=\left( \frac{m_0}{M_0 C_{knm} }\right)^\frac{1}{p_0}=\left(\frac{m_0(n+2)}{2 M_0}\right)^{\frac{1}{p_0}}\left( \frac{4S_n m(k-1)(n+2)}{n\chi M_0(k+m-1)^2} \right)^{\frac{1}{\alpha-m}}
\end{align}
by the choice of $k'$.

Next we will demonstrate the boundedness of $\|u\|_{L^k(\R^n)}$ starting from the case of $k=p_0.$

\noindent {\it\textbf{Step 1}}~(Decay estimate in $\|u\|_{L^{p_0}(\R^n)}$ )\quad Taking $k=p_0$ in \er{gj4.1} we obtain
\begin{align}\label{**}
    \frac{d}{dt}\int_{\R^n} u^{p_0}dx +p_0\chi m_0 C_{p_0}^{-p_0}\left( C_{p_0}^{p_0}-\|u\|_{L^{p_0}(\R^n)}^{p_0} \right)\|u\|_{L^{p_0+\alpha-1}(\R^n)}^{p_0+\alpha-1}\le 0.
\end{align}
Since we assume
\begin{align}\label{p00}
\|u_{0\varepsilon}\|_{L^{p_0}(\R^n)}<C_{p_0},
\end{align}
by bootstrap arguments on \er{**} we have the following estimate
\begin{align}
    \|u\|_{L^{p_0}(\R^n)}< \|u_{0\varepsilon}\|_{L^{p_0}(\R^n)}<C_{p_0}.
\end{align}
By H\"{o}lder inequality, \er{**} can be rewritten as
\begin{align}
    \frac{d}{dt}\int_{\R^n} u^{p_0}dx + \frac{p_0\chi m_0(C_{p_0}^{p_0}-\|u_{0\varepsilon}\|_{L^{p_0}(\R^n)}^{p_0})}{C_{p_0}^{p_0}M_0^{\frac{\alpha-1}{p_0-1}}} \left( \int_{\R^n} u^{p_0}dx \right)^{\frac{p_0+\alpha-2}{p_0-1}} \le 0.
\end{align}
After some calculations we have
\begin{align}\label{gj4.2}
    \|u\|_{L^{p_0}(\R^n)}
     \le \left( \frac{1}{\frac{p_0\chi m_0(C_{p_0}^{p_0}-\|u_{0\varepsilon}\|_{L^{p_0}(\R^n)}^{p_0})(\alpha-1)}{(p_0-1)C_{p_0}^{p_0}M_0^{\frac{\alpha-1}{p_0-1}}}t+\|u_{0\varepsilon}\|_{L^{p_0}(\R^n)}^{-\frac{p_0(\alpha-1)}{p_0-1}}} \right)^{\frac{p_0-1}{p_0(\alpha-1)}} \le C(1+t)^{-\frac{p_0-1}{p_0(\alpha-1)}},
\end{align}
where $C$ is a constant depending on $\|u_{0\varepsilon}\|_{L^{p_0}(\R^n)}, C_{p_0}$, $p_0$, $\alpha, \chi$, $M_0,m_0$. Integrating \er{**} from $0$ to $T$ in time we obtain that for any $T>0$
\begin{align*}
    \int_{\R^n} u^{p_0}(T)dx +p_0\chi m_0 C_{p_0}^{-p_0}\left( C_{p_0}^{p_0}-\|u_{0\varepsilon}\|_{L^{p_0}(\R^n)}^{p_0} \right)\int_0^T\int_{\R^n} u^{p_0+\alpha-1}dxdt\le \int_{\R^n} u_{0\varepsilon}^{p_0}dx.
\end{align*}
This assures that for any $T>0$
\begin{align}\label{linjie}
    u\in L^\infty(0,T;L^{p_0}(\R^n)),~~u\in L^{p_0+\alpha-1}(0,T;L^{p_0+\alpha-1}(\R^n)).
\end{align}
{\it\textbf{Step 2}}~(Decay estimates in $\|u\|_{L^k(\R^n)}$ for $1<k<\infty$)\quad In this step, we will show the decay properties of $\|u\|_{L^k(\R^n)}$ based on the decay of $\|u\|_{L^{p_0}(\R^n)}$ in time. We divide $k$ into two cases $1<k<p_0$ and $p_0<k<\infty$. \\
\textbf{(1) $1<k<p_0$.} By \er{gj4.2} one has
\begin{align}
    \|u\|_{L^k(\R^n)} &\le \|u\|_{L^{p_0}(\R^n)}^{p_0\frac{k-1}{k(p_0-1)}}\|u\|_{L^1(\R^n)}^{\frac{p_0-k}{k(p_0-1)}} \le C(1+t)^{-\frac{k-1}{k(\alpha-1)}},
\end{align}
where $C$ depends on $C_{p_0},p_0,k,\|u_{0\varepsilon}\|_{L^{p_0}(\R^n)},\alpha,\chi,m_0,M_0$. On the other hand, the inequality
\begin{align*}
   \|u\|_{L^{p_0}(\R^n)} <\|u_{0\varepsilon}\|_{L^{p_0}(\R^n)}<C_{p_0}<C_k
\end{align*}
is seen to hold because of the decreasing of $C_k$ with $k$. Hence integrating \er{gj4.1} from $0$ to $T$ in time we obtain that for any $T>0$
\begin{align*}
    \int_{\R^n} u^k(T)dx +k\chi m_0 C_k^{-p_0} \left( C_k^{p_0}-\|u_{0\varepsilon}\|_{L^{p_0}(\R^n)}^{p_0} \right)\int_0^T\int_{\R^n} u^{k+\alpha-1}dx dt\le \int_{\R^n} u_{0\varepsilon}^k dx
\end{align*}
which gives the following regularities that for any $T>0$
\begin{align}\label{xylinjie}
    u\in L^\infty \left(0,T;L^k(\R^n)\right),~~u\in L^{k+\alpha-1}\left(0,T;L^{k+\alpha-1}(\R^n)\right).
\end{align}
\textbf{(2) $p_0<k<\infty$.} Following \cite{cpz04} we consider the $L^k$ norm of the function $(u-N)+$ with $N>1$
\begin{align}
   & \frac{d}{dt}\int_{\R^n} (u-N)_+^kdx\nonumber \\
   = &-km(k-1)\int_{\mathbb{R}^n} (u-N)_+^{k-2}u^{m-1}|\nabla u|^2dx +k\chi\int_{\R^n} (u-N)_+^{k-1}u^\alpha dx \left(M_0-\int_{\R^n} udx\right)\nonumber\\
    \le & -\frac{4mk(k-1)}{(k+m-1)^2}\|\nabla(u-N)_+^{\frac{k+m-1}{2}}\|_{L^2(\R^n)}^2+k\chi\int_{\R^n} (u-N)_+^{k-1}(u-N+N)^\alpha dx \left(M_0-\int_{\R^n} udx\right)\nonumber\\
   \le & -\frac{4mk(k-1)}{(k+m-1)^2}\|\nabla(u-N)_+^{\frac{k+m-1}{2}}\|_{L^2(\R^n)}^2+k\chi M_02^{\alpha-1}\int_{\R^n} (u-N)_+^{k-1}\left( (u-N)^\alpha+N^\alpha\right)dx\nonumber\\
    =& -\frac{4mk(k-1)}{(k+m-1)^2}\|\nabla(u-N)_+^{\frac{k+m-1}{2}}\|_{L^2(\R^n)}^2+k\chi M_0 2^{\alpha-1}\int_{\R^n} (u-N)_+^{k+\alpha-1}dx\nonumber\\
   & +k\chi M_0 N^\alpha 2^{\alpha-1}\int_{\R^n} (u-N)_+^{k-1}dx.\label{I}
\end{align}
The terms involving $\int_{\R^n}(u-N)_+^{k+\alpha-1}dx$  and $\int_{\R^n} (u-N)_+^{k-1}dx $ can be estimated as follows:
\begin{align}
  \int_{\R^n}(u-N)_+^{k+\alpha-1}dx \le S_n^{-1}\|\nabla (u-N)_+^{\frac{k+m-1}{2}}\|_{L^2(\R^n)}^2 \|(u-N)_+\|_{L^{\frac{n(\alpha-m)}{2}}(\R^n)}^{\alpha-m}
\end{align}
and
\begin{align}\label{u-N}
  \int_{\R^n}(u-N)_+^{k-1}dx
  & =\int_{N<u\le N+1} (u-N)_+^{k-1}dx+\int_{u>N+1} (u-N)_+^{k-1}dx,\\
  \int_{N<u\le N+1}(u-N)_+^{k-1}dx
  & \le \int_{N<u \le N+1} 1dx\le\frac{1}{N}\int_{N<u\le N+1}udx\le\frac{m(t)}{N},\nonumber\\
  \int_{u> N+1}(u-N)_+^{k-1}dx & \le \int_{u>N+1} (u-N)_+^kdx\le\int_{\R^n} (u-N)_+^kdx\nonumber.
\end{align}
Collecting everything together \er{I} becomes
 \begin{align}
  & \frac{d}{dt}\int_{\R^n} (u-N)_+^kdx \nonumber \\
  \le
  & -\frac{4mk(k-1)}{(k+m-1)^2}\|\nabla(u-N)_+^{\frac{k+m-1}{2}}\|_{L^2(\R^n)}^2+k\chi M_0N^{\alpha-1}2^{\alpha-1}m(t)+k\chi M_0N^\alpha 2^{\alpha-1}\int_{\R^n} (u-N)_+^kdx\nonumber\\
  & +k\chi M_02^{\alpha-1}S_n^{-1}\|\nabla(u-N)_+^{\frac{k+m-1}{2}}\|_{L^2(\R^n)}^2\|(u-N)_+\|_{L^{\frac{n(\alpha-m)}{2}}(\R^n)}^{\alpha-m}.\label{N}
\end{align}
Next, let us observe that for any fixed $\max\left(1,\frac{n(\alpha-m)}{2}-(\alpha-1)\right)<k<\infty$ and under the condition \er{p00}, from \er{gj4.1} we may choose a $k_0>p_0$ such that $\|u(\cdot,t)\|_{L^{p_0}(\R^n)}< \|u_{0\varepsilon}\|_{L^{p_0}(\R^n)}\le C_{k_0}<C_{p_0}$ which guarantees $C_{k_0}^{p_0}-\|u\|_{L^{p_0}(\R^n)}^{p_0}>0$, then one has that for any $t>0$
\begin{align}\label{u}
 \|u(\cdot,t)\|_{L^{k_0}(\R^n)}\le\|u_{0\varepsilon}\|_{L^{k_0}(\R^n)}.
\end{align}
Using H\"older inequality it can be estimated that
\begin{align*}
 \|(u-N)_+\|_{L^{\frac{n(\alpha-m)}{2}}(\R^n)}
 &\le\|(u-N)_+\|_{L^{k_0}(\R^n)}\Big(\int_{u(t)\ge N}dx\Big)^{\frac{2}{n(\alpha-m)}-\frac{1}{k_0}}\\
 &\le \|(u-N)_+\|_{L^{k_0}(\R^n)}\Big(\int_{u(t)\ge N}\frac{u}{N}dx\Big)^{\frac{2}{n(\alpha-m)}-\frac{1}{k_0}}\\
 &\le \|u\|_{L^{k_0}(\R^n)}\Big(\frac{m(t)}{N}\Big)^{\frac{2}{n(\alpha-m)}-\frac{1}{k_0}}.
\end{align*}
Hence we may choose $N=N(k)$ sufficiently large such that for any $0<t<\infty$
\begin{align}
  \|(u-N)_+\|_{L^{\frac{n(\alpha-m)}{2}}(\R^n)}^{\alpha-m} \le\|u_{0\varepsilon}\|_{L^{k_0}(\R^n)}^{\alpha-m} \left(\frac{m(t)}{N}\right)^{\frac{2}{n}-\frac{\alpha-m}{k_0}}\le \frac{4m(k-1)S_n}{(k+m-1)^2 \chi M_0 2^{\alpha-1}}.
\end{align}
As a consequence, we infer from (\ref{N}) that
\begin{align*}
   \frac{d}{dt}\int_{\R^n} (u-N)_+^k dx\le k\chi M_0 N^\alpha 2^{\alpha-1} \int_{\R^n} (u-N)_+^k dx+k\chi M_0 N^{\alpha-1} 2^{\alpha-1} m(t).
\end{align*}
Finally, Gronwall inequality follows that for any $0\le t<\infty$
\begin{align}\label{u-N2}
  \int_{\R^n} (u-N)_+^kdx\le \Big(\int_{\R^n} (u_{0\varepsilon}-N)_+^kdx\Big)e^{k\chi M_0 N^\alpha 2^{\alpha-1} t}+\frac{M_0}{N}\left(e^{k\chi M_0N^\alpha 2^{\alpha-1} t}-1\right).
\end{align}
To go further, we claim that the bound on $\int_{\R^n} (u-N)_+^kdx$ is enough to treat $\int_{\R^n} u^kdx$. We decompose $\int_{\R^n} u^k dx$ in short and long range parts
\begin{align*}
    \int_{\R^n} u^kdx=\int_{u\le N} u^kdx+\int_{u>N} u^k dx.
\end{align*}
Then the short range part enjoys good properties for our purpose,
\begin{align*}
    \int_{u\le N} u^kdx\le N^{k-1}m(t)\le N^{k-1}M_0.
\end{align*}
As for the long range part we write
\begin{align*}
  \int_{u>N} u^kdx &=\int_{u>N} (u-N+N)^{k-1}udx\\
  &\le \max(2^{k-2},1)\left(\int_{u>N} (u-N)^{k-1}udx+N^{k-1}M_0\right)\\
  &= \max(2^{k-2},1)\left(N\int_{\R^n} (u-N)_+^{k-1}dx+\int_{\R^n}(u-N)_+^kdx+N^{k-1}M_0\right) \\
  &\le \max(2^{k-2},1)\left(m(t)+N\int_{\R^n} (u-N)_+^k dx+\int_{\R^n} (u-N)_+^k dx+N^{k-1}M_0\right),
\end{align*}
where the last line is derived from (\ref{u-N}). Therefore, the previous inequality (\ref{u-N2}) warrants that for any $T>0$
\begin{align}\label{t<T_k}
    \int_{\R^n} u^k dx \le C(k,M_0,N,e^t), ~\mbox{for}~ 0<t<T.
\end{align}
Now we can claim that $\int_{\R^n} u(\cdot,t)^k dx$ decays in time at infinity. Actually, for $t$ is larger than some $T_k$ one has
\begin{align}\label{define}
    S_n^{-1}\chi M_0\|u\|_{L^{p_0}(\R^n)}^{\alpha-m} \le \frac{2m(k-1)}{(k+m-1)^2}, ~~\mbox{for} ~t>T_k
\end{align}
because of the decay property of $\|u\|_{L^{p_0}(\R^n)}$. Then plugging \er{define} into \er{fsk} yields that for $t>T_k$
\begin{align}\label{Coll}
    & \frac{d}{dt}\int_{\R^n} u^kdx +\frac{2mk(k-1)}{(k+m-1)^2}\|\nabla u^{\frac{k+m-1}{2}}\|_{L^2(\R^n)}^2 +k\chi\int_{\R^n} udx\int_{\R^n} u^{k+\alpha-1} dx \le 0
\end{align}
by Lemma \ref{GNS}. Hence by H\"{o}lder inequality one has
\begin{align*}
    \frac{d}{dt}\int_{\R^n} u^kdx +\frac{k\chi m_0}{M_0^{\frac{\alpha-1}{k-1}}}\left( \|u\|_{L^k(\R^n)}^k \right)^{1+\frac{\alpha-1}{k-1}} \le 0,  ~~\mbox{for} ~t>T_k
\end{align*}
which leads to
\begin{align}\label{t>T_k}
    \|u\|_{L^k(\R^n)} \le C\left(M_0,m_0,k,\|u(\cdot,T_k)\|_{L^k(\R^n)}\right)(1+t-T_k)^{-\frac{k-1}{k(\alpha-1)}},~~\mbox{for} ~t>T_k.
\end{align}
Here $\|u(\cdot,T_k)\|_{L^k(\R^n)}$ is bounded from above by a positive constant $C\left(M_0,k,T_k\right)$ due to \er{t<T_k}. Moreover, taking \er{t<T_k} and \er{t>T_k} into account we conclude that for any $0<T<\infty$
\begin{align}\label{0<t<infty}
    \int_0^T\int_{\R^n} u^{k+\alpha-1} dxdt \le C \left(\|u_{0\varepsilon}\|_{L^k(\R^n)}, k, m_0,T\right),~~\mbox{for}~~ 1<k<\infty.
\end{align}
So taking \er{linjie}, \er{xylinjie} and \er{0<t<infty} together we deduce that for any $T>0$ and any $1<k<\infty$
 \begin{align}\label{sum}
     u\in L^\infty \left(0,T;L^k(\R^n)\right),~~ u\in L^{k+\alpha-1}\left(0,T;L^{k+\alpha-1}(\R^n)\right).
 \end{align}

 Furthermore, integrating \er{fsk} from $0$ to $T$ in time we conclude that for any $T>0$ and $1<k<\infty$
 \begin{align*}
     &\int_{\R^n} u^k(T)dx +\frac{4mk(k-1)}{(k+m-1)^2}\int_0^T\int_{\R^n} |\nabla u^{\frac{m+k-1}{2}}|^2 dx dt +k\chi\int_0^T\int_{\R^n} u dx \int_{\R^n} u^{k+\alpha-1}dx dt \\
     = &\int_{\R^n} u_{0\varepsilon}^k dx+ k\chi M_0 \int_0^T \int_{\R^n} u^{k+\alpha-1} dx dt.
 \end{align*}
As a result, we have that for any $0<T<\infty$ and $1<k<\infty$
\begin{align}
    \nabla u^{\frac{k+m-1}{2}} \in L^2\left(0,T;L^2(\R^n)\right),
\end{align}
by \er{0<t<infty}. Thus ends the proof.~\quad$\Box$

On account of the above arguments, our last task is to give the uniform boundedness of solutions for any $t>0$.
\begin{proposition}\label{Prop3}
(Uniform estimate in $L^\infty$) Let the same assumptions as that in proposition \ref{Prop2} hold. Then there exists a positive constant $C$  depending on $\|u_{0\varepsilon}\|_{L^1\cap L^\infty(\R^n)}$, $\alpha$, $m, n, \chi, m_0, M_0$ but not on $\varepsilon$ such that $u_\varepsilon$ is uniformly bounded for any $t> 0$, i.e.
\begin{align}
    \|u_\varepsilon(\cdot,t) \|_{L^\infty(\R^n)} \le C.
\end{align}
\end{proposition}
\noindent{\it Proof.} Firstly, we denote $q_k=2^k+n\alpha+nm$ and estimate $\int_{\R^n} u^{q_k}dx$. Multiplying \er{nkpp} with $q_k u^{q_k-1}$ we have
\begin{align}\label{guji2}
&\frac{d}{dt}\int_{\R^n} u^{q_k}dx+\frac{4mq_k(q_k-1)}{(q_k+m-1)^2}\int_{\R^n} |\nabla u^{\frac{q_k+m-1}{2}}|^2 dx+q_k\chi\int_{\R^n}udx\int_{\R^n} u^{q_k+\alpha-1} dx\nonumber\\
= &M_0 q_k\chi\int_{\R^n} u^{q_k+\alpha-1}dx.
\end{align}
Let
\begin{align*}
   q=\frac{2(q_k+\alpha-1)}{q_k+m-1},~~r=\frac{2q_{k-1}}{q_k+m-1},
\end{align*}
it's easy to verify that $1<\frac{q}{r}<\frac{2n}{r(n-2)}$ and $\frac{q}{r}<\frac{2}{r}+\frac{2}{n}$, so we take
$w=u^{\frac{q_k+m-1}{2}}$ in Lemma \ref{2018dcds} to obtain
\begin{align}\label{guji3}
\int_{\R^n} u^{q_k+\alpha-1}dx \le C(n)C_0^{-\frac{1}{\delta_1-1}}\Big(\int_{\R^n}u^{q_{k-1}}\Big)^{\gamma_1}+C_0\|\nabla u^{ \frac{q_k+m-1}{2}}\|^2_{L^2(\R^n)},
\end{align}
where
\begin{align*}
     & \delta_1 =\frac{q_k+m-1-(1-2/n)q_{k-1}}{q_k+\alpha-1-q_{k-1}}>1+2/n, \\
     & \gamma_1 =1+\frac{q_k+\alpha-1-q_{k-1}}{q_{k-1}+\frac{n}{2}(m-\alpha)}<2~ ~\mbox{for}~~\alpha \ge 1,~ m>1,
\end{align*}
 $C_0$ is a positive constant to be determined. Substituting (\ref{guji3}) into (\ref{guji2}) yields
\begin{align}
  & \frac{d}{dt}\int_{\R^n} u^{q_k}dx+\left(\frac{4mq_k(q_k-1)}{(q_k+m-1)^2}-C_0 M_0 q_k\chi\right)\int_{\R^n} |\nabla u^{\frac{q_k+m-1}{2}}|^2 dx \nonumber\\
   \le & C(n)q_k\chi C_0^{-\frac{1}{\delta_1-1}}\left(\int_{\R^n} u^{q_{k-1}}dx \right)^{\gamma_1}.
\end{align}
Notice that $\frac{4mq_k(q_k-1)}{(q_k+m-1)^2}>2m$, thus choosing $C_0=\frac{m}{M_0 q_k\chi }$ we have
$$\frac{4mq_k(q_k-1)}{(q_k+m-1)^2}-C_0 M_0 q_k\chi\ge m$$
which follows
\begin{align}\label{guji5}
   \frac{d}{dt}\int_{\R^n} u^{q_k}dx+m\int_{\R^n} |\nabla u^{\frac{q_k+m-1}{2}}|^2 dx \le C(n,m,\chi,M_0)q_k^{\frac{\delta_1}{\delta_1-1}}\Big(\int_{\R^n} u^{q_{k-1}}dx \Big)^{\gamma_1}.
\end{align}
On the other hand, taking
\begin{align*}
  w=u^{\frac{q_k+m-1}{2}},~~ q=\frac{2q_k}{q_k+m-1},~~ r=\frac{2q_{k-1}}{q_k+m-1},~~C_0=m
\end{align*}
in Lemma \ref{2018dcds} gives
\begin{align}\label{guji6}
   \int_{\R^n} u^{q_k}dx \le C(n)m^{-\frac{1}{\delta_2-1}}\left(\int_{\R^n} u^{q_{k-1}}dx \right)^{\gamma_1}+m\int_{\R^n} |\nabla u^{\frac{q_k+m-1}{2}}|^2 dx,
\end{align}
where
\begin{align*}
    \delta_2
    & =\frac{q_k+m-1-(1-2/n)q_{k-1}}{q_k-q_{k-1}}=O(1),\\
    \gamma_2
    & =1+\frac{q_k-q_{k-1}}{q_{k-1}+n(m-1)/2}<2 ~~\mbox{for}~~\alpha \ge 1,~ m>1.
\end{align*}
We may insert (\ref{guji6}) into (\ref{guji5}) and take the fact $\gamma_1 <2$, $\gamma_2<2$ into account to get
\begin{align}
    &\frac{d}{dt}\int_{\R^n} u^{q_k}dx+\int_{\R^n} u^{q_k}dx\nonumber\\
    \le &
    C(n,m,\chi,M_0)q_k^{\frac{\delta_1}{\delta_1-1}}\left(\int_{\R^n} u^{q_{k-1}}dx \right)^{\gamma_1}+C(n,m)\left(\int_{\R^n} u^{q_{k-1}}dx \right)^{\gamma_2} \nonumber\\
   \le & C(n,m,\chi,M_0)q_k^{\frac{\delta_1}{\delta_1-1}}\left(\left(\int_{\R^n} u^{q_{k-1}}dx \right)^{\gamma_1}+\left(\int_{\R^n} u^{q_{k-1}}dx \right)^{\gamma_2}\right) \nonumber\\
     \le & 2C(n,m,\chi,M_0)q_k^{\frac{\delta_1}{\delta_1-1}}\max\left(\left(\int_{\R^n} u^{q_{k-1}}dx\right)^2,1\right) \nonumber\\
    \le & 2C(n,m,\chi,M_0)q_k^{1+n/2}\max\left(\left(\int_{\R^n} u^{q_{k-1}}dx\right)^2,1\right) \nonumber\\
    \le & 2C(n,m,\chi,M_0)2^{k(1+n/2)}(n\alpha+nm+1)^{1+n/2}\max\left(\left(\int_{\R^n} u^{q_{k-1}}dx\right)^2,1\right)\label{guji7}.
\end{align}
Here we have used the fact that $\delta_1>1+2/n$. Let $K_0=\max\Big(1,\|u_{0\varepsilon}\|_{L^1(\R^n)}),\|u_{0\varepsilon}\|_{L^\infty (\R^n)}\Big)$, we have the following inequality for the initial data
$$
\int_{\R^n} u_{0\varepsilon}^{q_k}dx\le [\max(\|u_{0\varepsilon}\|_{L^1(\R^n)}),\|u_{0\varepsilon}\|_{L^\infty(\R^n)})]^{q_k}\le K_0^{q_k}.
$$
Denote
$$y_k(t)=\int_{\R^n} u^{q_k}dx,~~ d_0=1+n/2,~~ \Bar{a}=C(n,m,\chi,M_0)(n\alpha+nm+1)^{d_0},$$
(\ref{guji7}) can be recasted as
$$
y'_k(t)+y_k(t)\le 2\Bar{a}2^{d_0k}\max\left(y_{k-1}^2(t),1\right).
$$
By virtue of Lemma 4.1 of \cite{BL14} one can solve that
\begin{align}\label{guji8}
    \int_{\R^n} u(\cdot,t)^{q_k}dx\le (2\Bar{a})^{2^k-1}2^{d_0(2^{k+1}-k-2)}\max\Bigg(\sup_{t\ge0}\Big(\int_{\R^n}u(\cdot,t)^{q_0}dx\Big)^{2^k},K_0^{q_k}\Bigg).
\end{align}
Recalling $q_k=2^k+n\alpha+nm$ and taking the power $\frac{1}{q_k}$ to both sides of (\ref{guji8}) we have the uniformly boundedness of the solution $u$ by passing to the limit $k\to\infty $
\begin{align}\label{20210701}
   \|u(\cdot,t)\|_{L^\infty(\R^n)}\le 2\Bar{a}2^{2d_0}\max\Big(\sup_{t\ge0}\int_{\R^n} u(\cdot,t)^{q_0}dx,K_0\Big).
\end{align}
Thanks to Proposition \ref{Prop2} by $q_0=1+n\alpha+nm$, it allows us to find
\begin{align*}
    \int_{\R^n} u(\cdot,t)^{1+n\alpha+nm}dx\le C\left(\|u_{0\varepsilon}\|_{L^{1+n\alpha+nm}(\R^n)},m_0,M_0 \right).
\end{align*}
Thus \er{20210701} implies
\begin{align*}
   \|u(\cdot,t)\|_{L^\infty(\R^n)}\le C(K_0).
\end{align*}
This is exactly the anticipated result.\quad\quad$\Box$

Combining the local existence result Proposition \ref{Prop1} and the uniform boundedness Proposition \ref{Prop2}, \ref{Prop3} we close the proof of Theorem \ref{global}.

\section{Proof of the main theorems} \label{sec4}
\def\theequation{4.\arabic{equation}}\makeatother
\setcounter{equation}{0}
\def\thetheorem{4.\arabic{theorem}}\makeatother
\setcounter{theorem}{0}

In this section, we give the proof of Theorem \ref{sub}, \ref{critical} and \ref{super} and show the global existence of a weak solution to \er{nkpp} for the three cases.

\noindent{\it Proof of Theorem \ref{sub}, \ref{critical} and \ref{super}:} By virtue of proposition \ref{Prop2} and proposition \ref{Prop3}, for the initial data satisfies $u_{0\varepsilon}\in L^1\cap L^\infty(\R^n)$, the following basic estimates are obtained that for any $T>0:$
\begin{align}
    & \|u_\varepsilon\|_{L^\infty(0,T;L^1\cap L^q(\R^n))}\le C,~~\mbox{for any}~1<q\le \infty,\label{i}\\
    & \|u_\varepsilon\|_{L^{k+\alpha-1}(0,T;L^{k+\alpha-1}(\R^n))}\le C,~~ \mbox{for any}~ 1<k\le \infty,\label{ii}\\
    &\left\|\nabla u_\varepsilon^{\frac{m+k-1}{2}}\right\|_{L^2(0,T;L^2(\R^n))}\le C,~~ \mbox{for any}~ 1<k<\infty,\label{iii}
\end{align}
where $C$ are constants depending on $\|u_{0\varepsilon}\|_{L^1\cap L^\infty(\R^n)}$, $m$, $\alpha$, $n, \chi, m_0, M_0$ but not on $\varepsilon$. Hence there exists a subsequence $u_\varepsilon$ without relabeling such that for any $T>0$
\begin{align}
    & u_\varepsilon \rightharpoonup u~ \mbox{in}~ L^q(0,T;L^1\cap L^q(\R^n)),~~\mbox{for any}~ 1\le q<\infty,\label{iv}\\
    & u_\varepsilon\stackrel{\ast}{\rightharpoonup} u~~\mbox{in}~ L^\infty(0,T;L^1\cap L^\infty(\R^n)).\nonumber
\end{align}
Furthermore we will show that for any $T>0:$
\begin{align}
   & \nabla u_\varepsilon^m \in L^\infty(0,T;L^2(\R^n)),\label{v}\\
   & (u_\varepsilon^m)_t \in L^2(0,T;L^2(\R^n)).\label{vi}
\end{align}
Multiplying \er{fkppepsilon} with $(u_\varepsilon^m)_t$ and $(u_\varepsilon)_t$ respectively we have
\begin{align}
    & \frac{4m}{(m+1)^2}\int_{\R^n} \left| \left(u_\varepsilon^{\frac{m+1}{2}}\right)_t \right|^2 dx \nonumber\\
    = & -\frac{1}{2}\frac{d}{dt}\int_{\R^n} |\nabla u_\varepsilon^m|^2 dx-\varepsilon\int_{\R^n} \nabla u_\varepsilon \cdot(\nabla u_\varepsilon^m)_t dx +\chi m\int_{\R^n} u_\varepsilon^{\alpha+m-1} (u_\varepsilon)_t dx\left(M_0-\int_{\R^n} u_\varepsilon dx \right) \nonumber\\
    \le & \frac{2m}{(m+1)^2}\int_{\R^n} |(u_\varepsilon^{\frac{m+1}{2}})_t|^2 dx +C(\chi,M_0,m)\int_{\R^n} u_\varepsilon^{2\alpha+m-1}dx,\label{mep}
\end{align}
and
\begin{align}
    &\varepsilon\int_{\R^n} |(u_\varepsilon)_t|^2dx \nonumber\\
    = & -\varepsilon\int_{\R^n} \nabla u_\varepsilon^m \cdot(\nabla u_\varepsilon)_t dx-\frac{\varepsilon^2}{2}\frac{d}{dt}\int_{\R^n} |\nabla u_\varepsilon|^2 dx +\chi\varepsilon\int_{\R^n} u_\varepsilon^\alpha (u_\varepsilon)_t dx\left(M_0-\int_{\R^n} u_\varepsilon dx\right) \nonumber\\
    \le & -\varepsilon\int_{\R^n} \nabla u_\varepsilon^m \cdot(\nabla u_\varepsilon)_t dx-\frac{\varepsilon^2}{2}\frac{d}{dt}\int_{\R^n} |\nabla u_\varepsilon|^2 dx +\frac{\varepsilon}{2}\int_{\R^n} |(u_\varepsilon)_t|^2 dx +C(\chi,M_0,\varepsilon)\int_{\R^n} u_\varepsilon^{2\alpha} dx.\label{ep}
\end{align}
Combining \er{mep} with \er{ep} one has
\begin{align}
    & \frac{2m}{(m+1)^2} \int_{\R^n} |u_\varepsilon^{\frac{m+1}{2}}|^2 dx +\frac{1}{2} \frac{d}{dt}\int_{\R^n} |\nabla u_\varepsilon^m|^2 dx \nonumber\\
    & +\frac{\varepsilon}{2}\int_{\R^n} |(u_\varepsilon)_t|^2 dx +\frac{\varepsilon^2}{2}\frac{d}{dt}\int_{\R^n} |\nabla u_\varepsilon|^2 dx +\varepsilon\frac{4m}{(m+1)^2}\frac{d}{dt}\int_{\R^n} |\nabla u_\varepsilon^{\frac{m+1}{2}}|^2 dx \nonumber\\
    \le & C(\chi,M_0,m)\left( \int_{\R^n} u_\varepsilon^{2\alpha+m-1} dx +\int_{\R^n} u_\varepsilon^{2\alpha} dx \right).\label{zong1}
\end{align}
Integrating \er{zong1} with respect to time follows that for any $T>0$
\begin{align}
    & \frac{2m}{(m+1)^2}\int_0^T \int_{\R^n} \left| \left(u_\varepsilon^{\frac{m+1}{2}}\right)_t\right|^2 dx dt +\frac{1}{2}\displaystyle\sup_{0<t<T}\int_{\R^n} |\nabla u_\varepsilon^m|^2 dx \nonumber\\
    & +\frac{\varepsilon}{2}\int_0^T \int_{\R^n} |(u_\varepsilon)_t|^2 dx dt +\frac{\varepsilon^2}{2}\displaystyle\sup_{0<t<T} \int_{\R^n
    } |\nabla u_\varepsilon|^2 dx +\varepsilon\frac{4m}{(m+1)^2}\displaystyle\sup_{0<t<T} \int_{\R^n} |\nabla u_\varepsilon^{\frac{m+1}{2}}|^2 dx \nonumber\\
    \le & C(\chi,M_0,m)\int_0^T \int_{\R^n} u_\varepsilon^{2\alpha+m} dxdt + \int_0^T \int_{\R^n} u_\varepsilon^{2\alpha}dx dt +\frac{1}{2}\int_{\R^n} \left|\nabla u_{0\varepsilon}^m\right|^2 dx \nonumber\\
    & +\frac{\varepsilon^2}{2} \int_{\R^n} |\nabla u_{0\varepsilon}|^2 dx +\frac{4m\varepsilon}{(m+1)^2}\int_{\R^n} |\nabla u_{0\varepsilon}^{\frac{m+1}{2}}|^2 dx.\label{zong2}
\end{align}
By \er{i} and \er{ii} we see that for $\alpha \ge 1$, $m> 1$, there exists a positive constant $C$ which is independent of $\varepsilon$ such that
\begin{align*}
    &\int_0^T \int_{\R^n} |(u_\varepsilon^m)_t|^2 dx dt +\displaystyle\sup_{0<t<T} \int_{\R^n} |\nabla u_\varepsilon^m|^2 dx\\
    \le & \frac{4m^2}{(m+1)^2} \|u_\varepsilon\|_{L^\infty(Q_T)}^{m-1} \int_0^T \int_{\R^n} \left|(u_\varepsilon^{\frac{m+1}{2}})_t\right|^2 dxdt+\displaystyle\sup_{0<t<T} \int_{\R^n} |\nabla u_\varepsilon^m|^2dx\le C.
\end{align*}
It can be seen that $u_\varepsilon^m \in L^\infty(0,T;H^1(\R^n))\cap H^1(0,T;L^2(\R^n))$ which proves \er{v} and \er{vi}. Consequently, there exists a subsequence $u_\varepsilon$ without relabeling such that
\begin{align}\label{ys1}
    u_\varepsilon^m \to\xi ~~\mbox{in}~~ C\left(0,T;L_{loc}^2(\R^n)\right)
\end{align}
which directly gives that for any bounded domain $\Omega$
\begin{align}\label{ys2}
    u_\varepsilon \to \xi^{\frac{1}{m}} ~~\mbox{in}~~ \Omega,~t\in(0,T).
\end{align}
On the other hand, recalling \er{iv} and using the dominated convergence theorem leads to
\begin{align}
    u_\varepsilon\to u ~~\mbox{in}~~ L^q(0,T;L^1\cap L^q(\Omega))~~ \mbox{for any }~~ 1\le q<\infty,
\end{align}
and thus from \er{ys2} we find
\begin{align}
    u_\varepsilon \to \xi^{\frac{1}{m}}=u~~ \mbox{in}~~\Omega, ~t\in(0,T).
\end{align}
By \er{ys1} we arrive at
\begin{align}
    u_\varepsilon^m \to u^m~~ \mbox{in}~~ C\left(0,T;L^2(\Omega)\right).
\end{align}
Hence from \er{v} we conclude
\begin{align}\label{ast}
    \nabla u_\varepsilon^m \stackrel{\ast}{\rightharpoonup} \nabla u^m ~~\mbox{in}~~ L^\infty(0,T;L^2(\R^n)).
\end{align}
We can also obtain
\begin{align}
    u_\varepsilon \to u~~\mbox{in}~~ C\left(0,T;L^{2m}(\Omega)\right)
\end{align}
by virtue of the fact $|u_\varepsilon-u|^m\le |u_\varepsilon^m-u^m|$ for $m>1$. Because of
\begin{align*}
    |u_\varepsilon-u|^k\le (2\|u_\varepsilon\|_{L^\infty(0,T;L^\infty(\Omega))})^{k-2m}|u_\varepsilon-u|^{2m}~~ \mbox{for any }~~ k\ge 2m,
\end{align*}
we can go further to obtain that for any $2<k<\infty$
\begin{align}\label{ys3}
    u_\varepsilon \to u ~~\mbox{in}~~ C\left(0,T;L^k(\Omega)\right).
\end{align}

In addition, we need to prove
\begin{align}\label{ys4}
    u_\varepsilon \to u ~~\mbox{in}~~ L^q(0,T;L^1({\R^n}), ~~\mbox{for any }~~ 1<q<\infty.
\end{align}
Here we will apply the second moment estimate to establish the uniform integrability of $u_\varepsilon$ at far field.
From \er{fkppepsilon} we find
\begin{align*}
    \frac{d}{dt}\int_{\R^n} |x|^2 u_\varepsilon(\cdot,t)dx
    & =2n\int_{\R^n} u_\varepsilon^m dx +2n\varepsilon\int_{\R^n} u_\varepsilon dx+\chi\int_{\R^n} u_\varepsilon^\alpha |x|^2dx\left(M_0-\int_{\R^n} u_\varepsilon dx \right) \nonumber\\
    & \le C\left(\|u_\varepsilon\|_{L^1\cap L^\infty(\R^n)}\right)+\chi\|u_\varepsilon\|_{L^\infty(Q_T)}^{\alpha-1} M_0\int_{\R^n} u_\varepsilon |x|^2 dx.
\end{align*}
Then from \er{i} and by Gronwall inequality one gets
\begin{align}
    \int_{\R^n} u_\varepsilon |x|^2dx \le \left( \int_{\R^n} u_{0\varepsilon}|x|^2 dx +C_1 \right)e^{C_2t}<C,
\end{align}
where $C_1$, $C_2$ are constants depending on $\|u_{0\varepsilon}\|_{L^1\cap L^\infty(\R^n)}$, $\chi$, $M_0, \alpha,m,n$. Then we compute
\begin{align}\label{lim1}
    \int_0^T \|u_\varepsilon\|_{L^1(|x|>R)}^q dt \le \int_0^T \frac{\left(\int_{\R^n} u_\varepsilon |x|^2 dx\right)^q}{R^{2q}} dt \to 0 ~~\mbox{as}~~ R \to \infty
\end{align}
for any $1<q<\infty$, and the weak semi-continuity of $L^q(0,T;L^1(|x|>R))$ yields
\begin{align}\label{lim2}
    \int_0^T \|u\|_{L^1(|x|>R)}^q dt \le \displaystyle\liminf_{\varepsilon \to 0} \int_0^T \|u_\varepsilon\|_{L^1(|x|>R)}^q dt \to 0 ~~\mbox{as}~~ R \to \infty.
\end{align}
Therefore, the following inequality is derived that for any $1<q<\infty$, as $R\to \infty$, $\varepsilon \to0$,
\begin{align*}
     & \int_0^T \|u_\varepsilon-u\|_{L^1(\R^n)}^q dt =\int_0^T \left( \|u_\varepsilon-u\|_{L^1(|x|>R)}+\|u_\varepsilon-u\|_{L^1(|x|\le R)} \right)^q dt\\
    \le  & C(q)
     \left(\int_0^T
     \|u_\varepsilon\|_{L^1(|x|>R)}^q dt +\int_0^T \|u\|_{L^1(|x|>R)}^q dt +\int_0^T \|u_\varepsilon-u\|_{L^1(|x|\le R)}^q dt \right) \to 0.
\end{align*}
In the last inequality, the first term goes to zero due to \er{lim1}, the second term is given by \er{lim2} and \er{ys3} yields the third term, thus one proves \er{ys4}.

Now integrating \er{fkppepsilon} with respect to $x$, $t$ we get the weak formulation for $u_\varepsilon$
\begin{align}
    & \int_0^T \int_{\R^n} \nabla u_\varepsilon^m \cdot\nabla \varphi dxdt +\varepsilon\int_0^T \int_{\R^n} \nabla u_\varepsilon \cdot \nabla \varphi dxdt -\chi\int_0^T\int_{\R^n} u_\varepsilon^\alpha\varphi \left( M_0-\int_{\R^n} u_\varepsilon dx \right) dxdt \nonumber\\
    = & \int_{\R^n} u_{0\varepsilon}(x) \varphi(x,0) dx +\int_0^T \int_{\R^n} u_\varepsilon \varphi_t dxdt \label{zong3}
\end{align}
for any continuously differentiable function $\varphi$ with compact support in $\R^n \times [0,T)$. Thanks to \er{ast}, \er{ys3} and \er{ys4}, passing to the limit $\varepsilon \to 0$ in \er{zong3} we obtain
\begin{align}
   & \int_0^T \int_{\R^n} \nabla u^m \cdot \nabla \varphi dxdt -\chi\int_0^T \int_{\R^n} u^\alpha\varphi \left( M_0-\int_{\R^n} u dx \right) dxdt \nonumber\\
    = & \int_{\R^n} u_0(x)\varphi(x,0)dx +\int_0^T \int_{\R^n} u\varphi_t dxdt,
\end{align}
which finishes the proof of our main results. \quad $\Box$

\section{Numerical results} \label{sec5}
\def\theequation{5.\arabic{equation}}\makeatother
\setcounter{equation}{0}
\def\thetheorem{5.\arabic{theorem}}\makeatother
\setcounter{theorem}{0}

Throughout this section, we numerically consider the radial solutions satisfying
\begin{align}\label{radialPDE}
\left\{
  \begin{array}{ll}
    u_t=(u^m)_{rr}+\frac{n-1}{r}(u^m)_r+\chi u^\alpha \left( M_0-n \beta_n \int_0^\infty u(r,t)r^{n-1} dr \right), \quad  0<r<\infty, ~t>0, \\
    u'(0,t)=0,\quad u \to 0 \mbox{~~as~~} r \to \infty, \quad t>0, \\
    u(r,0)=u_0(r) \ge 0,
  \end{array}
\right.
\end{align}
where $n \ge 3$ and $\beta_n=\frac{\pi^{n/2}}{\Gamma(n/2+1)}$. Here the initial mass $m_0=n \beta_n \int_0^\infty u(r,t)r^{n-1} dr$ is assumed to be $m_0<M_0$ such that $M_0-m(t)$ remains non-negative for all $0<t<\infty$. A series of numerical experiments of the PDE \er{radialPDE} with different forms of initial data are presented in this section. The numerical results illustrate the theoretical predications of the earlier sections, as well as explore other issues that lie beyond the scope of the analysis. Numerical simulations are carried out using semi-implicit finite difference scheme for the diffusion term and linearized method specialized for the nonlinear reaction term. Here $\alpha$ is divided into three cases: $1\le \alpha<m+2/n$, $\alpha=m+2/n$ and $\alpha>m+2/n.$

\subsection{Global existence for the subcritical case $1 \le \alpha<m+2/n$}

\begin{figure}[htbp]
\begin{minipage}{0.5\linewidth}
\vspace{3pt}
\centerline{\includegraphics[width=8cm]{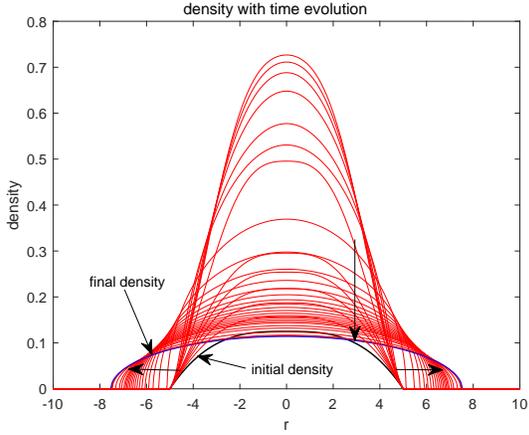}}
\centerline{(a) density with time evolution}
\end{minipage}
\begin{minipage}{0.5\linewidth}
\vspace{3pt}
\centerline{\includegraphics[width=8cm]{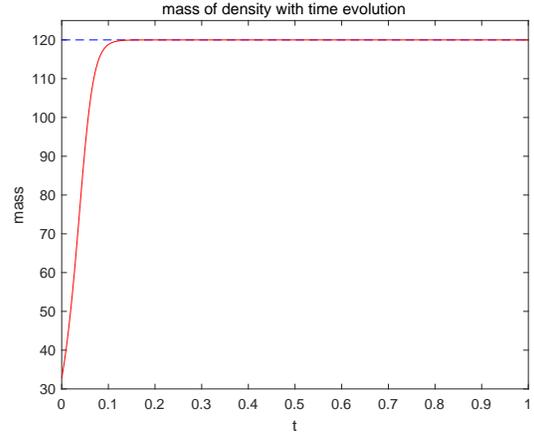}}
\centerline{(b) mass with time evolution}
\end{minipage}

\begin{minipage}{0.5\linewidth}
\vspace{3pt}
\centerline{\includegraphics[width=8cm]{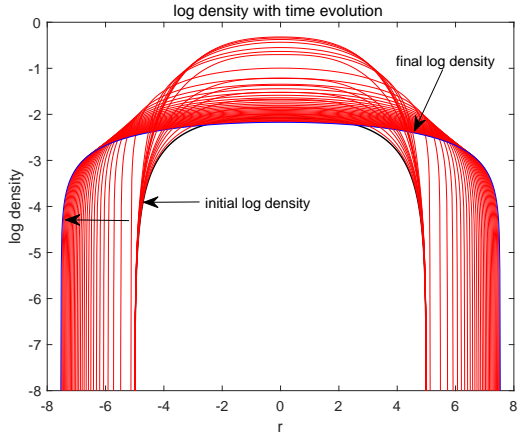}}
\centerline{(c) log density with time evolution}
\end{minipage}
\begin{minipage}{0.5\linewidth}
\vspace{3pt}
\centerline{\includegraphics[width=8cm]{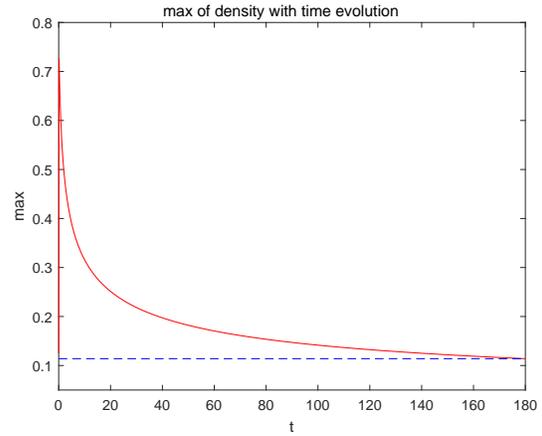}}
\centerline{(d) maximum with time evolution}
\end{minipage}
\caption{Numerical simulations of the global existence with compactly supported initial data for $n=3$: $m_0=32.72413808$, $M_0=120$, the density converges to the unique compactly supported steady solution with mass $M_0$.} \label{fig0}
\end{figure}

\begin{figure}[htbp]
\begin{minipage}{0.5\linewidth}
\vspace{3pt}
\centerline{\includegraphics[width=8cm]{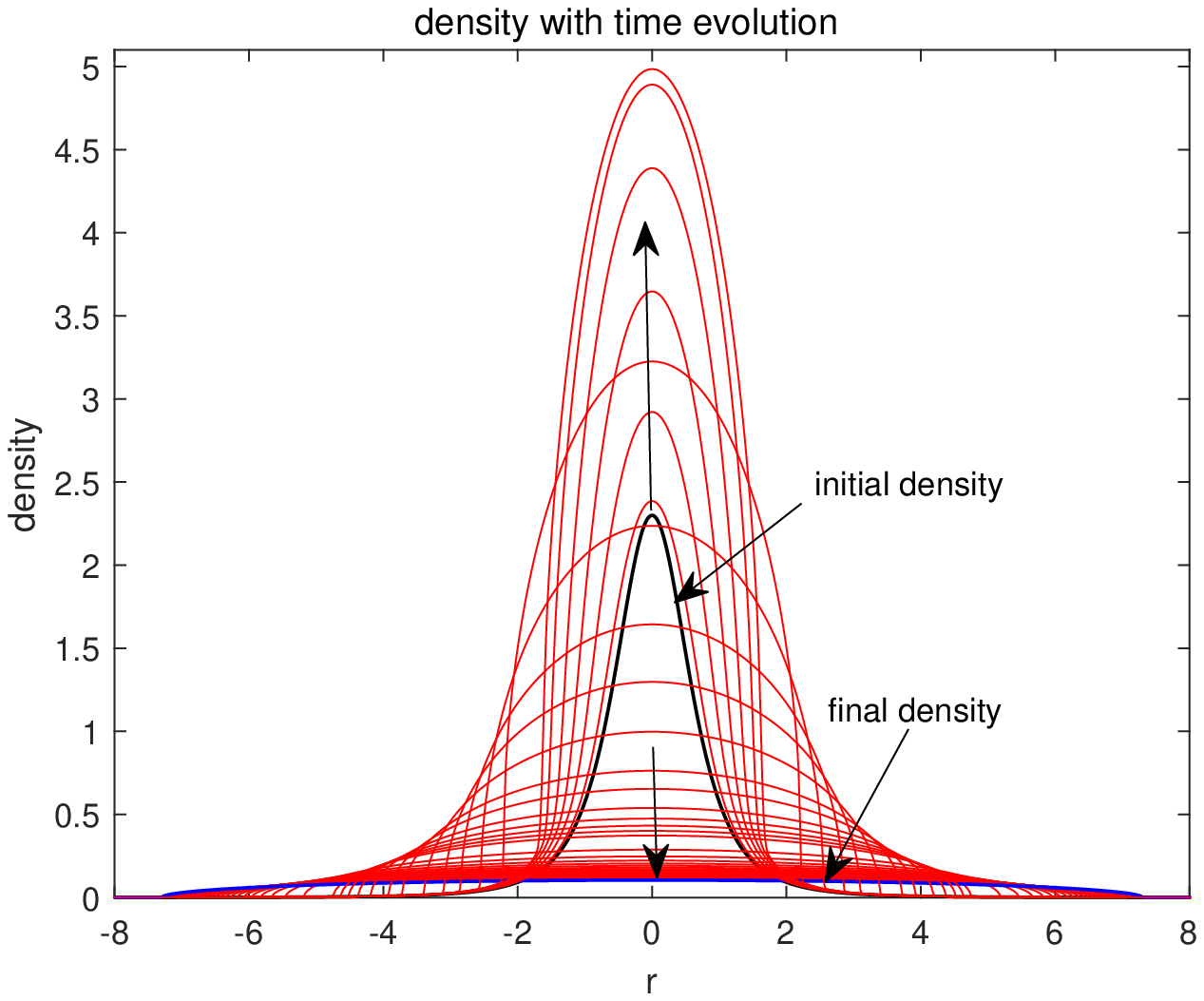}}
\centerline{(a) density with time evolution}
\end{minipage}
\begin{minipage}{0.5\linewidth}
\vspace{3pt}
\centerline{\includegraphics[width=8cm]{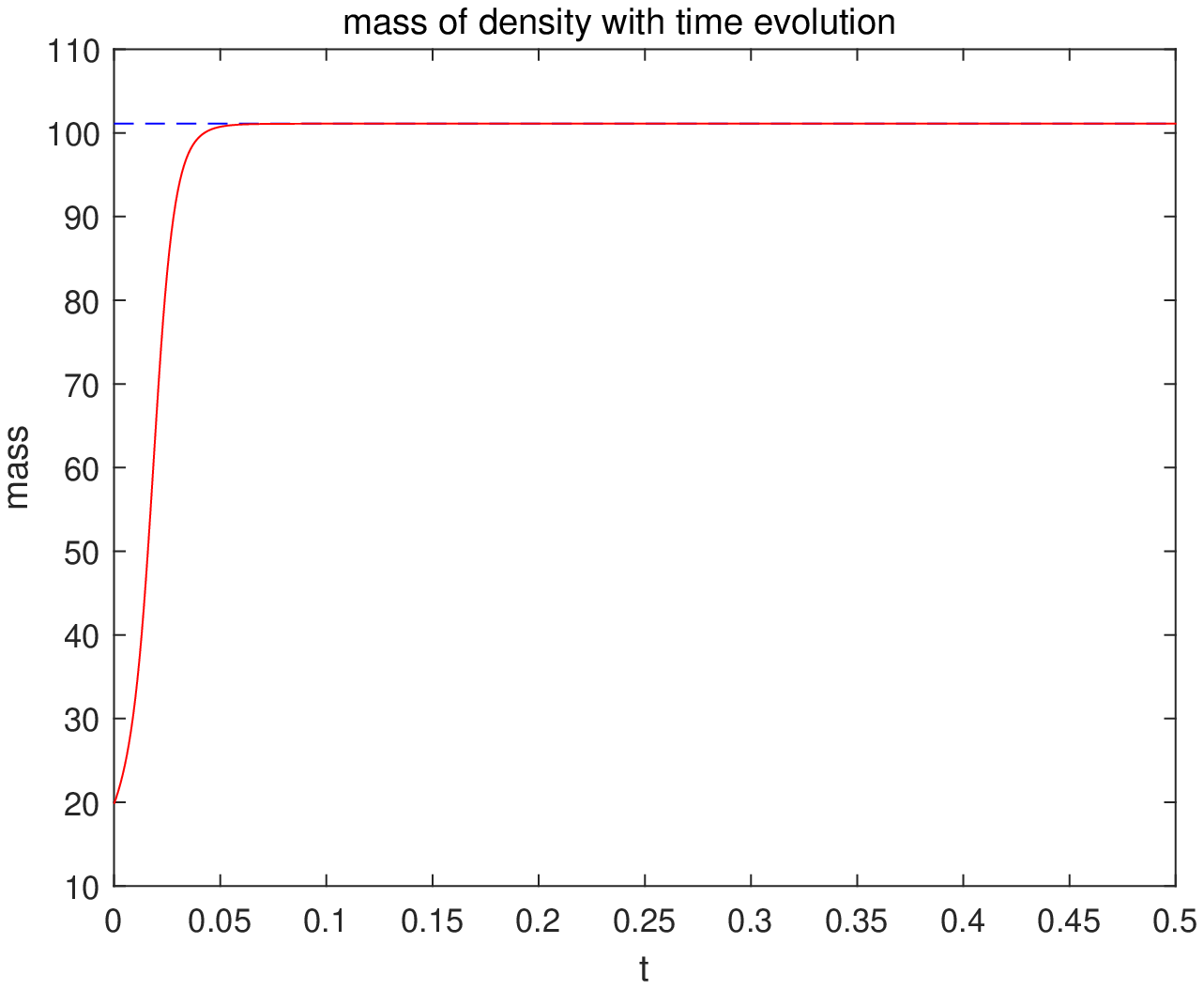}}
\centerline{(b) mass with time evolution}
\end{minipage}

\begin{minipage}{0.5\linewidth}
\vspace{3pt}
\centerline{\includegraphics[width=8cm]{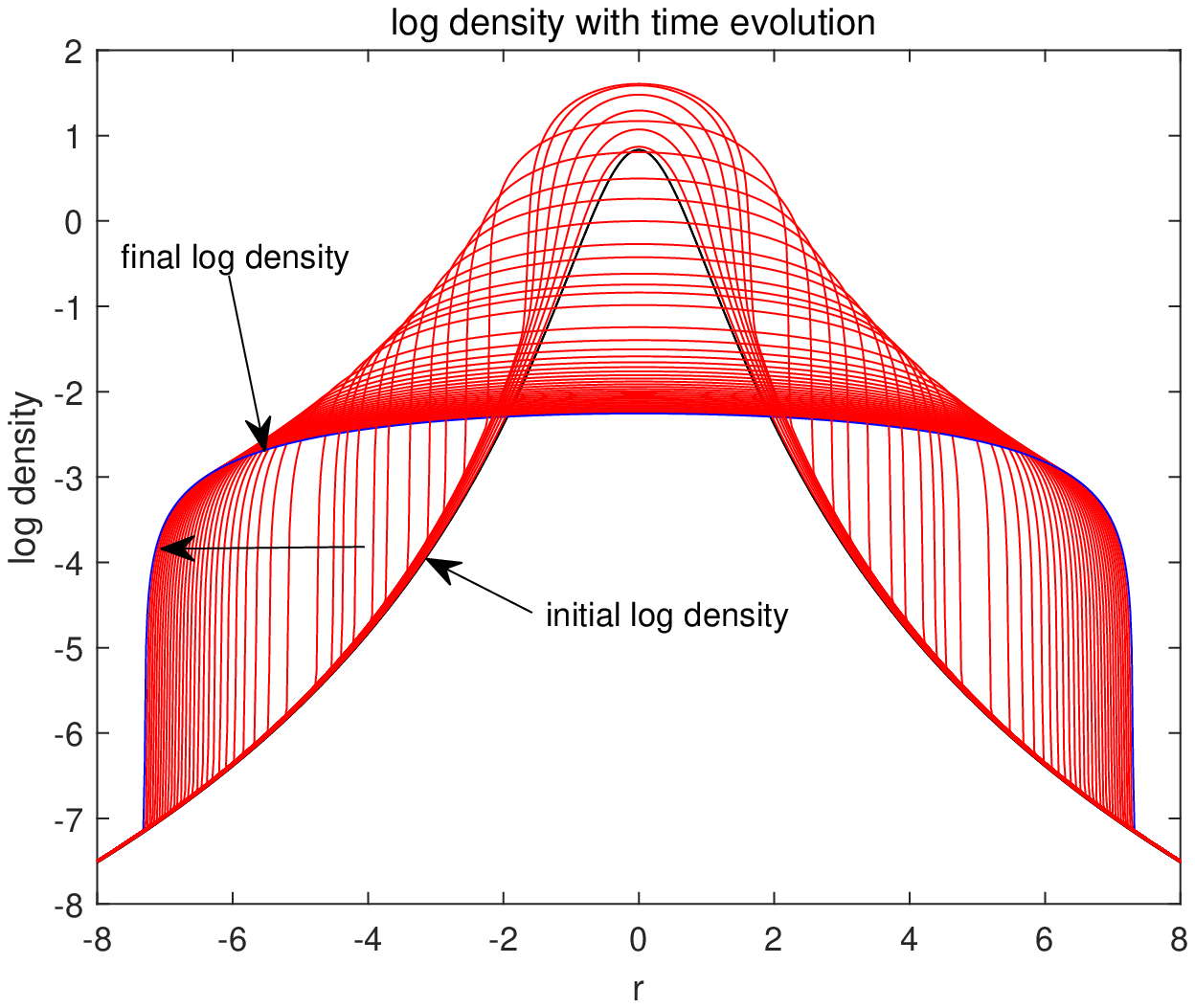}}
\centerline{(c) log density with time evolution}
\end{minipage}
\begin{minipage}{0.5\linewidth}
\vspace{3pt}
\centerline{\includegraphics[width=8cm]{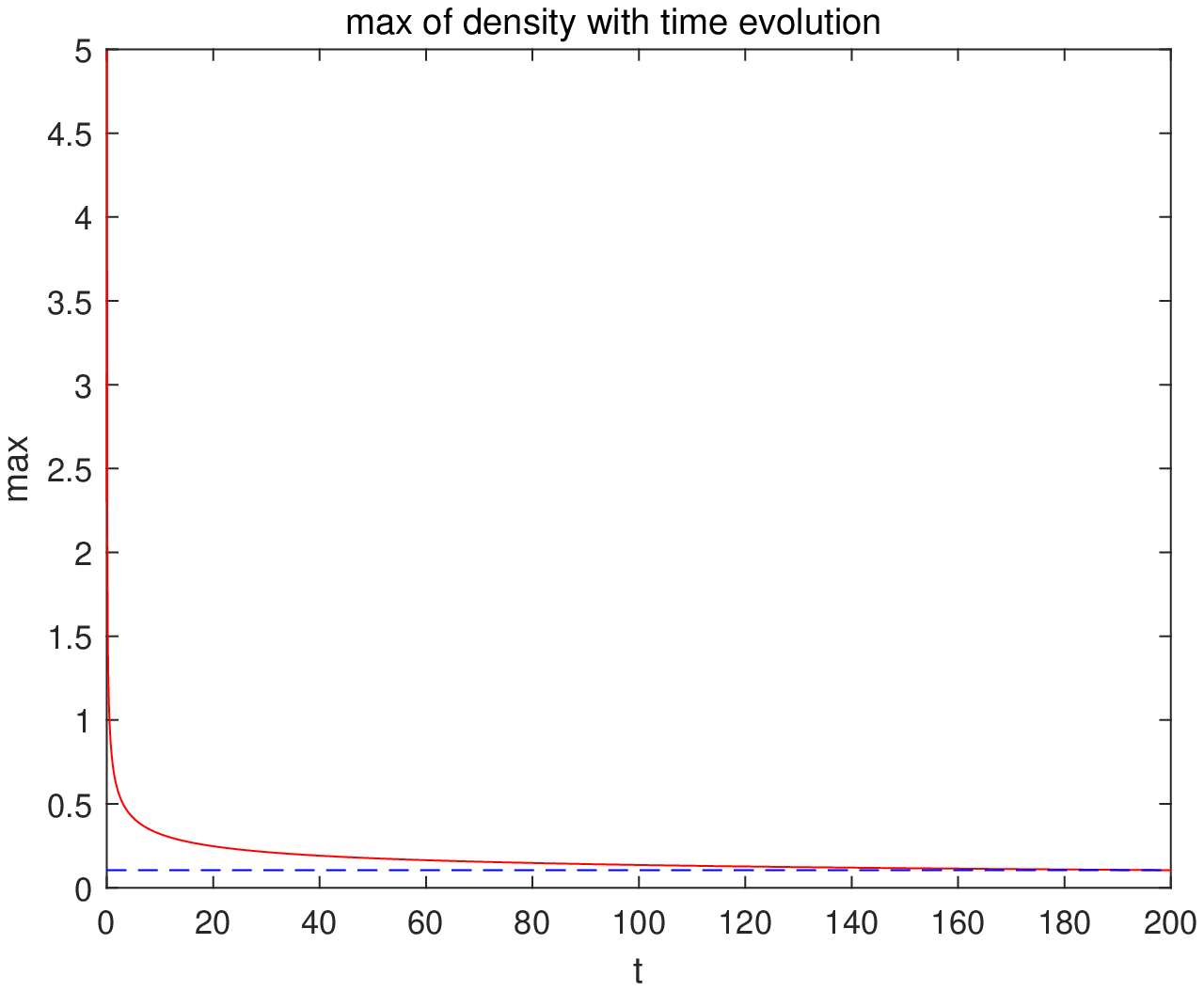}}
\centerline{(d) maximum with time evolution}
\end{minipage}
\caption{Numerical simulations of the global existence with non-compactly supported initial data for $n=3$: $m_0=19.82455437$, $M_0=101.10522729$, the density converges to the unique compactly supported steady solution with mass $M_0$. } \label{fig00}
\end{figure}

In this subsection, we solve the problem starting from non-negative initial data of forms: compactly supported, non-compactly supported. As is expected in Theorem \ref{sub}, for any given initial mass $m_0$ and $M_0,$ the solution converges to the compact supported steady profile with mass $M_0$, see Figure \ref{fig0}(a)(b) and Figure \ref{fig00}(a)(b). This is more evident in Figure \ref{fig0}(c)(d) and Figure \ref{fig00}(c)(d) where time-profiles of the solution are plotted on a log-scale graph, and the maximum of the solution tends to that of the stationary solution.

\subsection{Finite time blow-up for the critical case $\alpha=m+2/n$}

In Theorem \ref{critical}, we have proved that the solution will exist globally under the condition of $M_0 \le M_\ast$ where $M_\ast$ is defined as in \er{Mstar}. For $M_0> M_\ast$, we solve \er{radialPDE} starting from initial data with different masses $m_0$ and mass capacities $M_0$ to illustrate the influence of $M_0,m_0$ on the dynamics. Numerical experiments show that for any given initial mass $m_0$, there exists a unique critical value $M_c$ such that for $M_0<M_c$, all the solutions will converge to the unique steady solution with mass $M_0$. While for $M_0>M_c,$ all the solutions will blow up in finite time. We can infer from Table \ref{tab1} that higher dimensions require larger $M_c$ for finite time blow-up.

\begin{table}
\caption{Critical threshold $M_c$ on $M_0$ separating finite time blow-up and global existence}\label{tab1}
 \centering
\begin{tabular}{ccccc}
\hline
  \mbox{Dimension} & $n=3$ & $n=4$ & $n=5$ & $n=6$ \\
  \hline
  $M_c$ & 46.20083432 & 58.10808013 & 69.81554367 & 82.83702312 \\
  \hline
\end{tabular}
\\[1mm]
The initial mass $m_0=19.82455437$.
\end{table}

\begin{figure}[htbp]
\begin{minipage}{0.33\linewidth}
\vspace{3pt}
\centerline{\includegraphics[width=5.9cm]{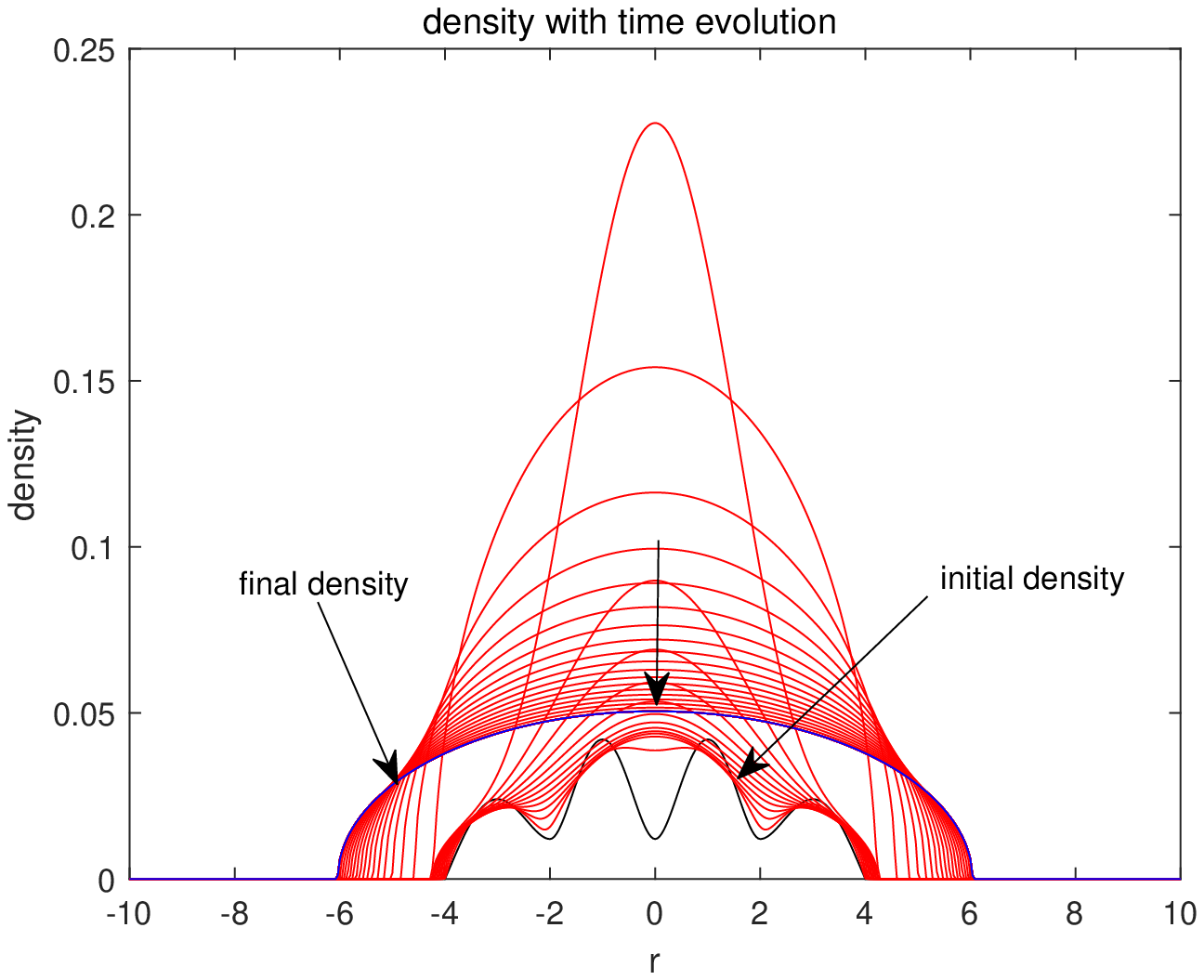}}
\centerline{(a1) convergence of density}
\end{minipage}
\begin{minipage}{0.33\linewidth}
\vspace{3pt}
\centerline{\includegraphics[width=5.9cm]{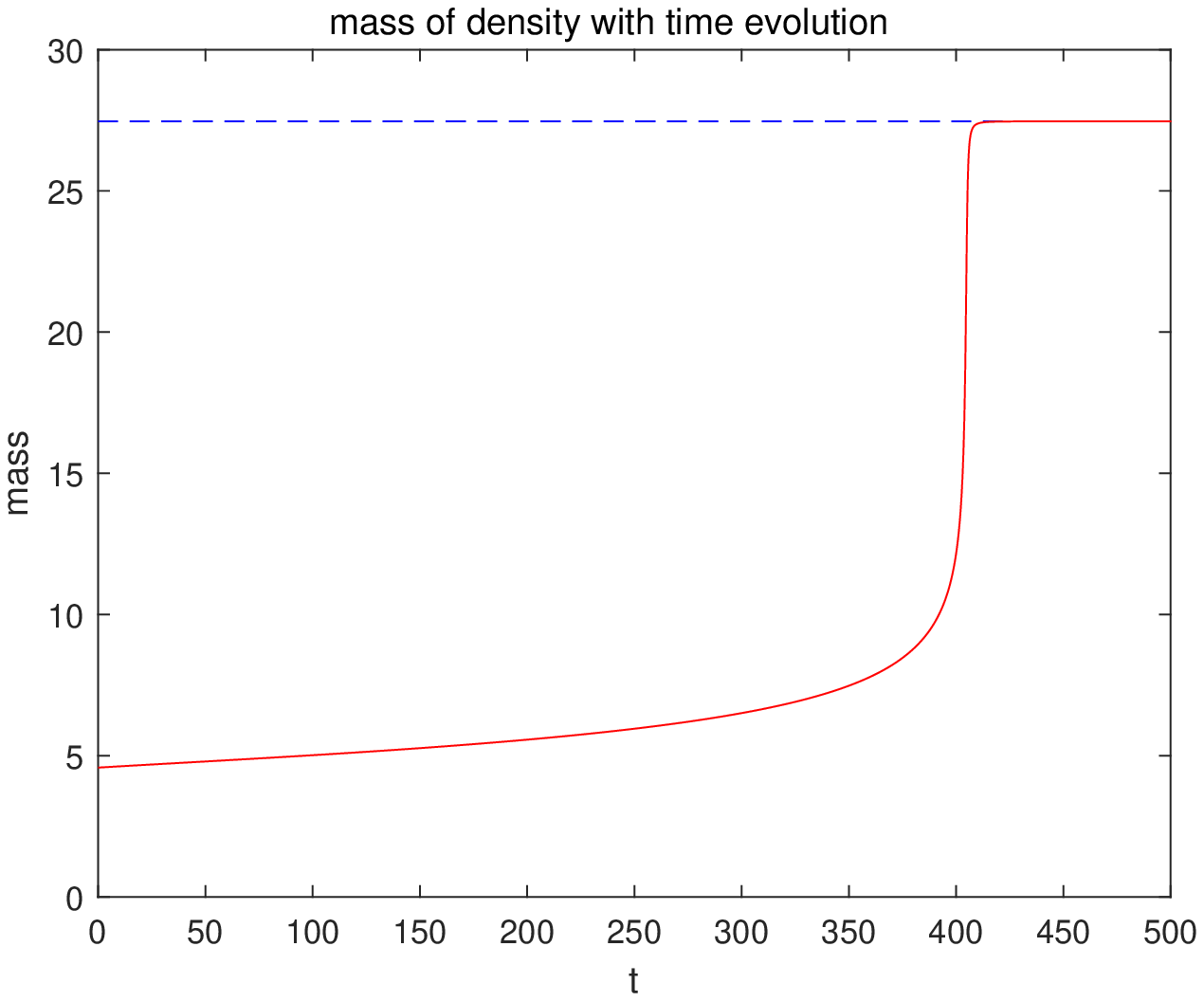}}
\centerline{(a2) convergence of mass}
\end{minipage}
\begin{minipage}{0.33\linewidth}
\vspace{3pt}
\centerline{\includegraphics[width=5.9cm]{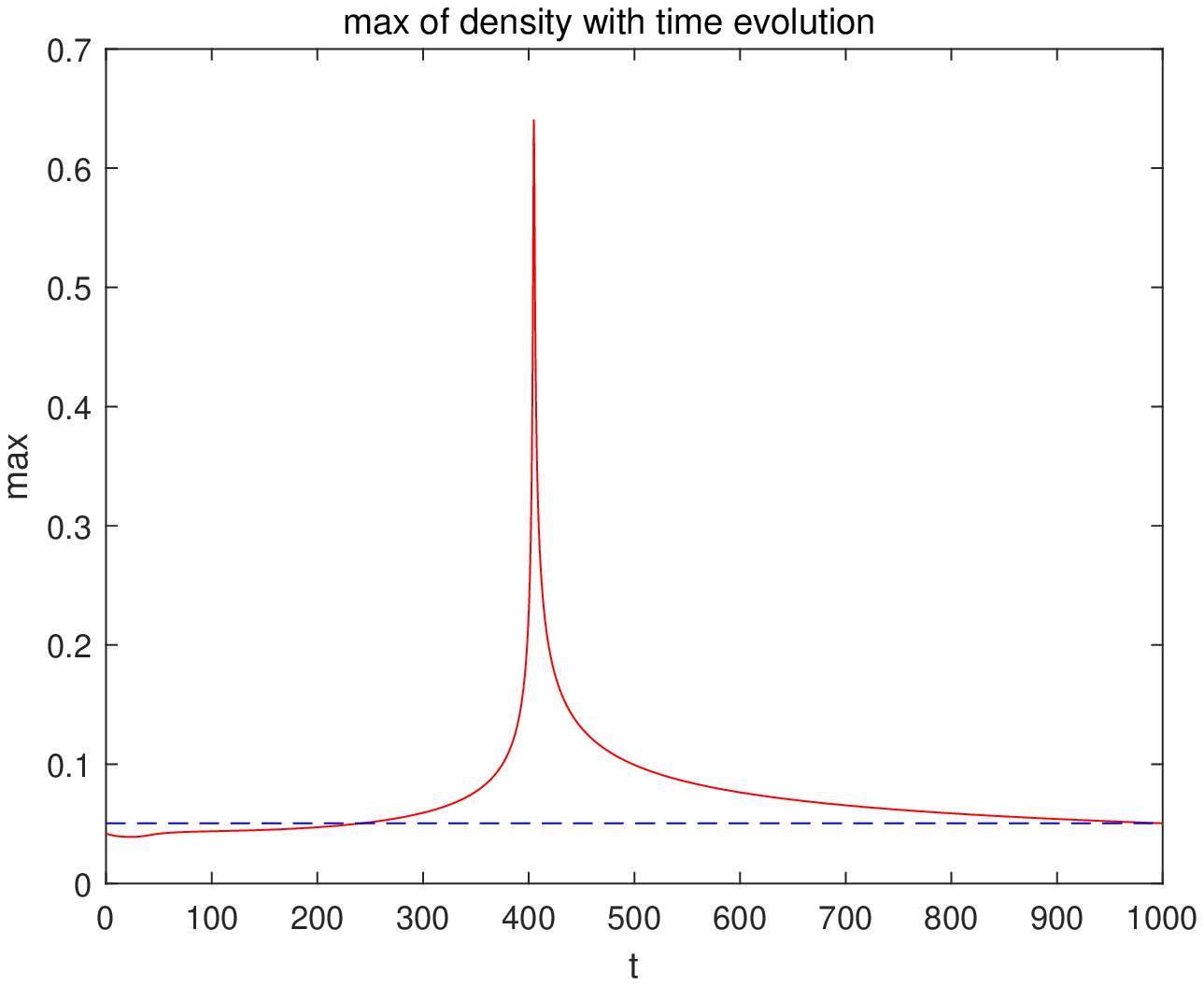}}
\centerline{(a3) evolution of maximum}
\end{minipage}

\begin{minipage}{0.33\linewidth}
\vspace{3pt}
\centerline{\includegraphics[width=5.9cm]{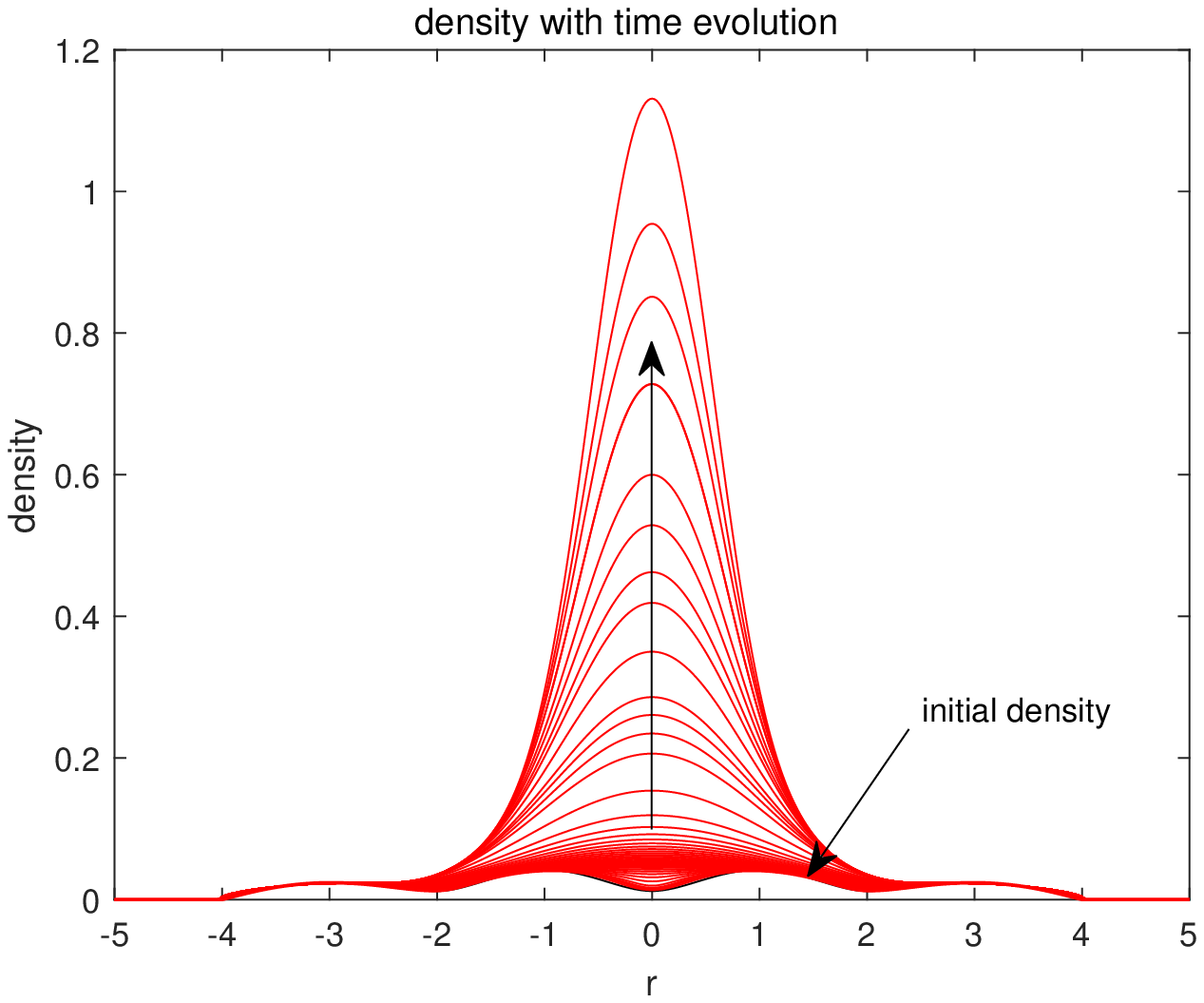}}
\centerline{(b1) blow-up of density}
\end{minipage}
\begin{minipage}{0.33\linewidth}
\vspace{3pt}
\centerline{\includegraphics[width=5.9cm]{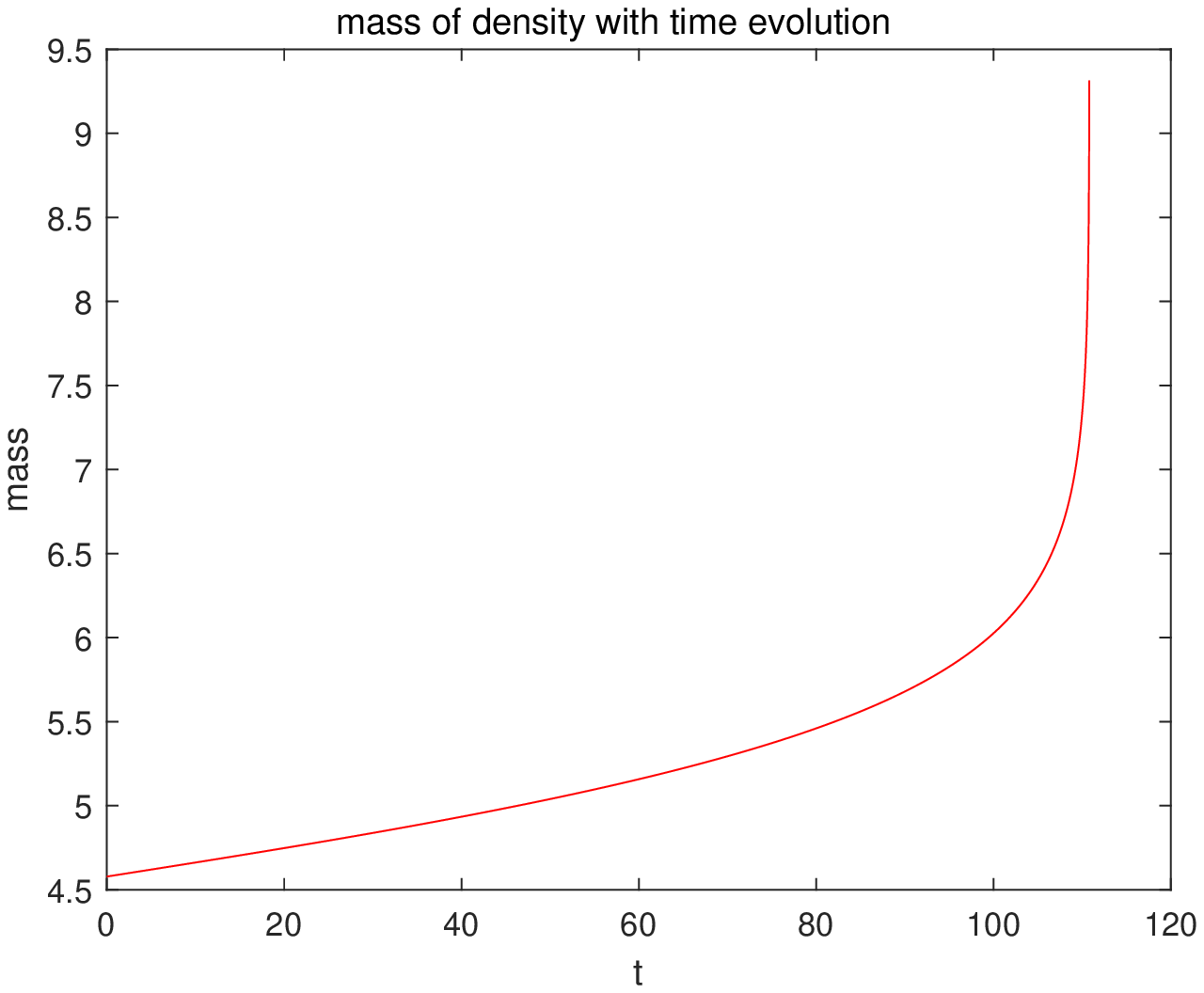}}
\centerline{(b2) evolution of mass}
\end{minipage}
\begin{minipage}{0.33\linewidth}
\vspace{3pt}
\centerline{\includegraphics[width=5.9cm]{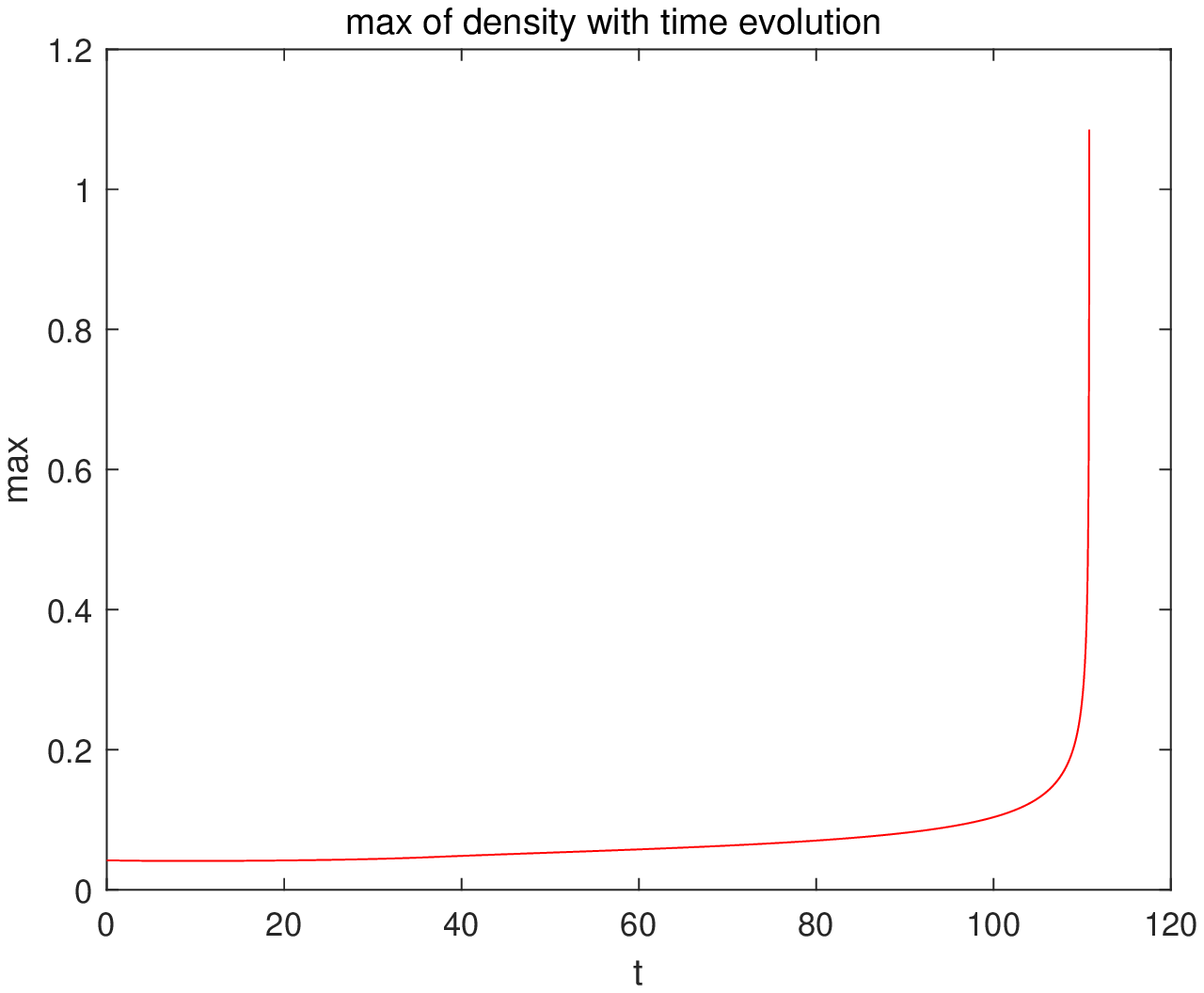}}
\centerline{(b3) evolution of maximum}
\end{minipage}
\caption{The initial mass $m_0=4.57840705$, the global existence (up) and finite time blow-up (down): (a1-a3) $M_0=27.47044232$, the solution will converge to the unique steady solution with mass $M_0$, (b1-b3) $M_0=45.78407053$, at $T_b=110.823$s, the solution forms a local maximum that evolves to produce blow-up and the mass increases to $9.31417986$. } \label{fig1}
\end{figure}

\begin{figure}[htbp]
\begin{minipage}{0.5\linewidth}
\vspace{3pt}
(a) \centerline{\includegraphics[width=8cm]{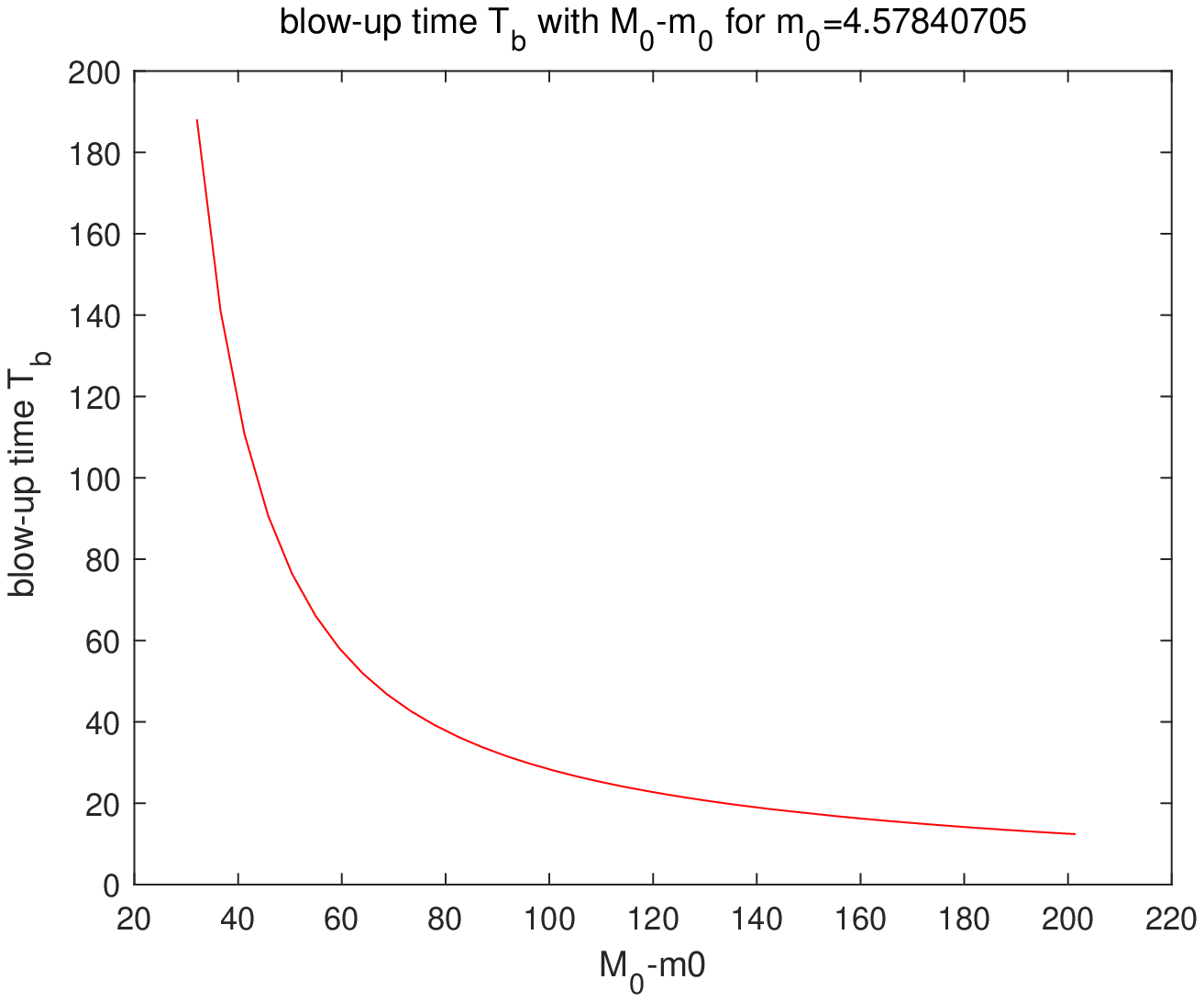}}
\end{minipage}
\begin{minipage}{0.5\linewidth}
\vspace{3pt}
(b) \centerline{\includegraphics[width=8cm]{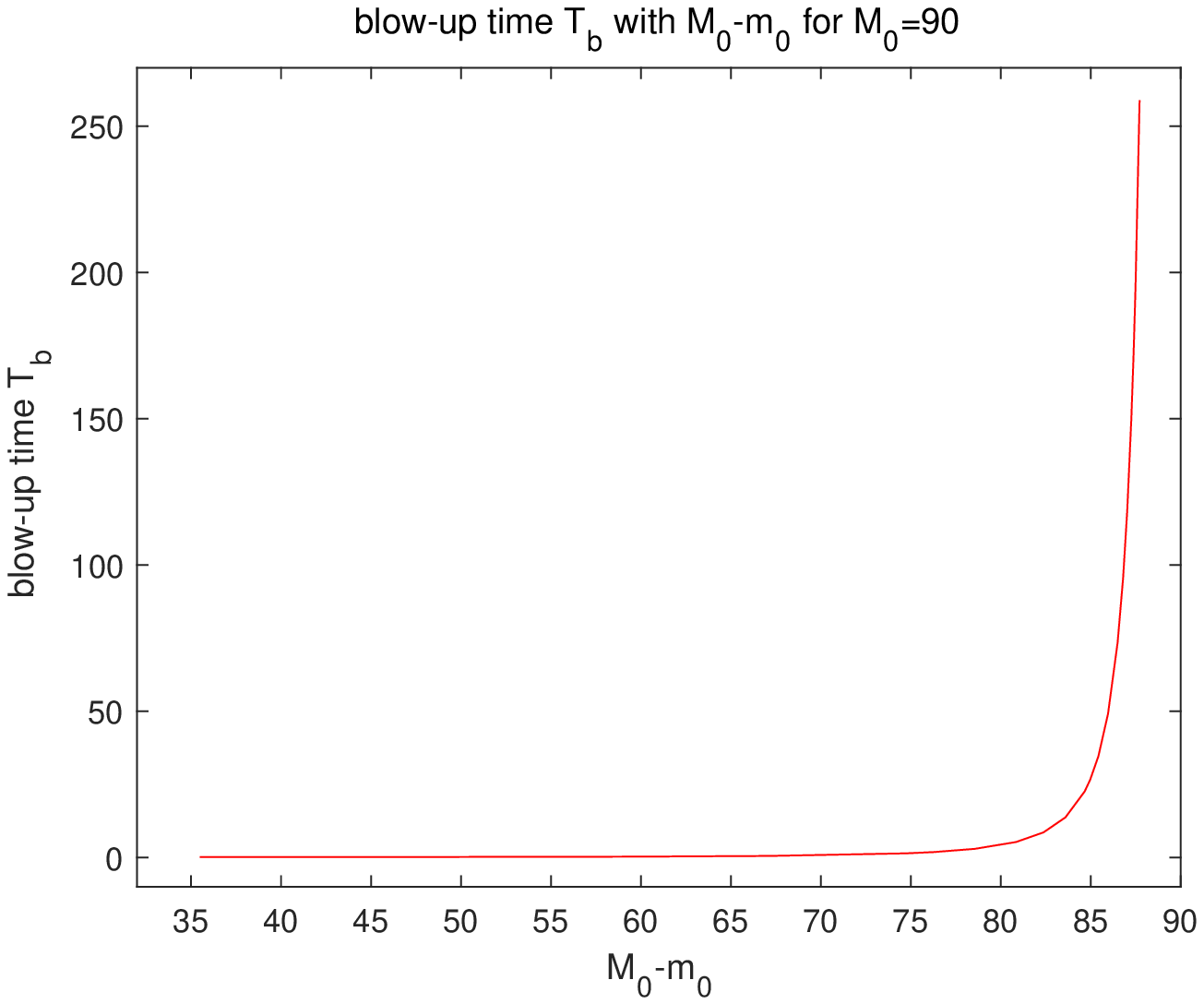}}
\end{minipage}
\caption{Relationship between the blow-up time $T_b$ and $M_0-m_0$: (a) given $m_0=4.57840705$, blow-up time $T_b$ corresponding to different $M_0$ in the range $M_0>34.36599205$, (b) given $M_0=90$, blow-up time $T_b$ for $m_0<54.53645869$. } \label{fig2}
\end{figure}

We take $n=3$ as an example, beginning with a multi-bump initial data with mass $m_0=4.57840705$, for $m_0<M_\ast<M_0<34.36599205$, by viewing the simulations of the maximum and mass of the solution, we observe that the bumps firstly move towards the center and the solution quickly increases to the maximum, then it decreases to spread outwards and finally converges to the compactly supported steady profile with mass $M_0$, see Figure \ref{fig1}(a1-a3). While for $M_0>34.36599205$, the coefficient $M_0-m_0$ of the growth term is large to prevent the solution from spreading, the bumps approach to zero and increase to form a local maximum that determines the position of a blow-up singularity, see Figure \ref{fig1}(b1). This is further verified by the time evolution of the mass and maximum of the solution (see Figure \ref{fig1}(b2-b3)). In addition, numerical simulations in the blow-up profile illustrate that the larger $M_0,$ the faster the solution blows up (see Figure \ref{fig2}(a)). On the other hand, for any given $M_0,$ there is a critical value $m_c$ such that the solution will blow up in finite time for $m_0<m_c$ (or $M_0-m_0>M_0-m_c$), and the blow-up time $T_b$ becomes longer for smaller $m_0,$ see Figure \ref{fig2}(b).

\begin{figure}[htbp]
\begin{minipage}{0.33\linewidth}
\vspace{3pt}
\centerline{\includegraphics[width=5.9cm]{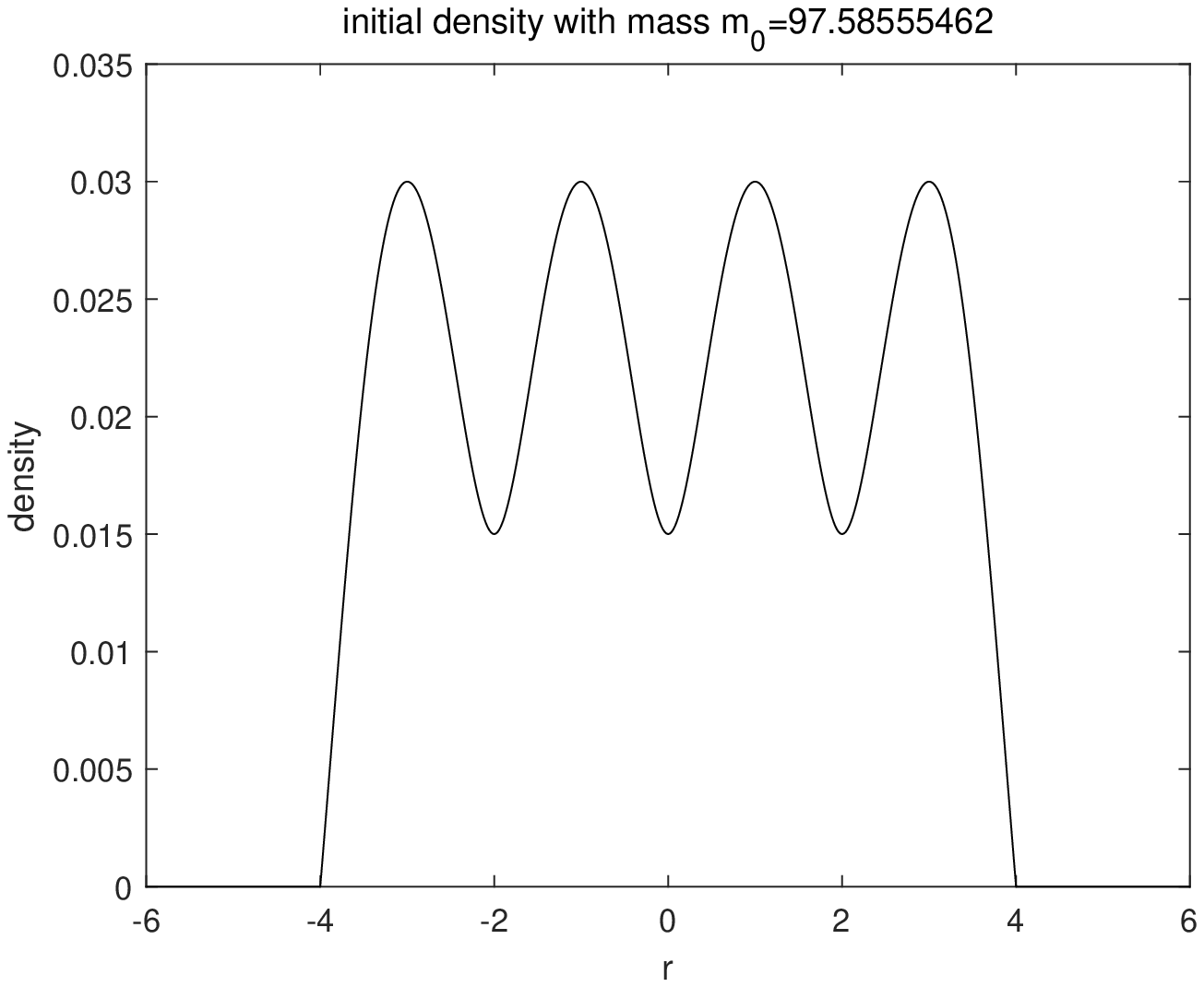}}
\centerline{(a) the initial density}
\end{minipage}
\begin{minipage}{0.33\linewidth}
\vspace{3pt}
\centerline{\includegraphics[width=5.9cm]{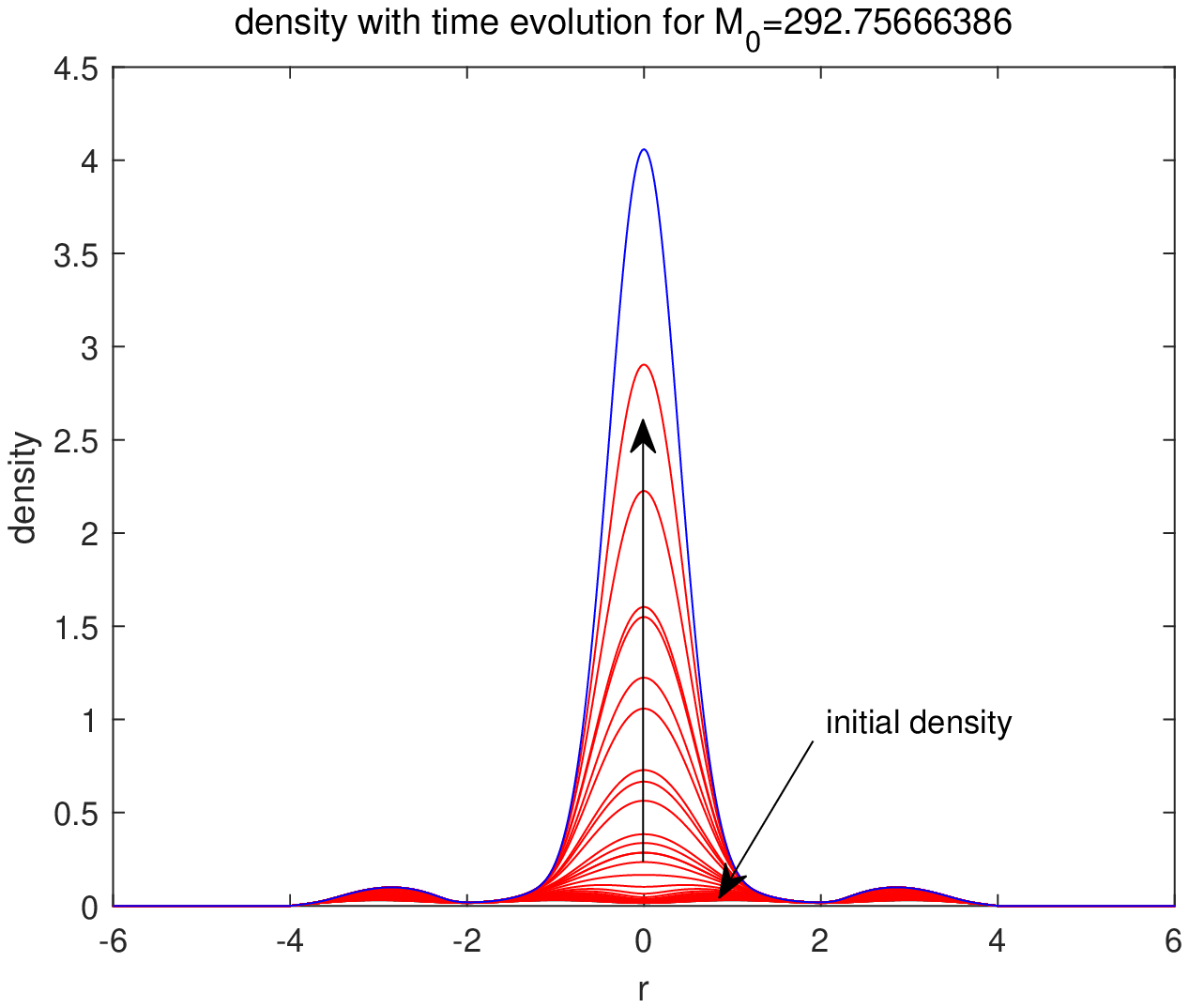}}
\centerline{(b) one singularity}
\end{minipage}
\begin{minipage}{0.33\linewidth}
\vspace{3pt}
\centerline{\includegraphics[width=5.9cm]{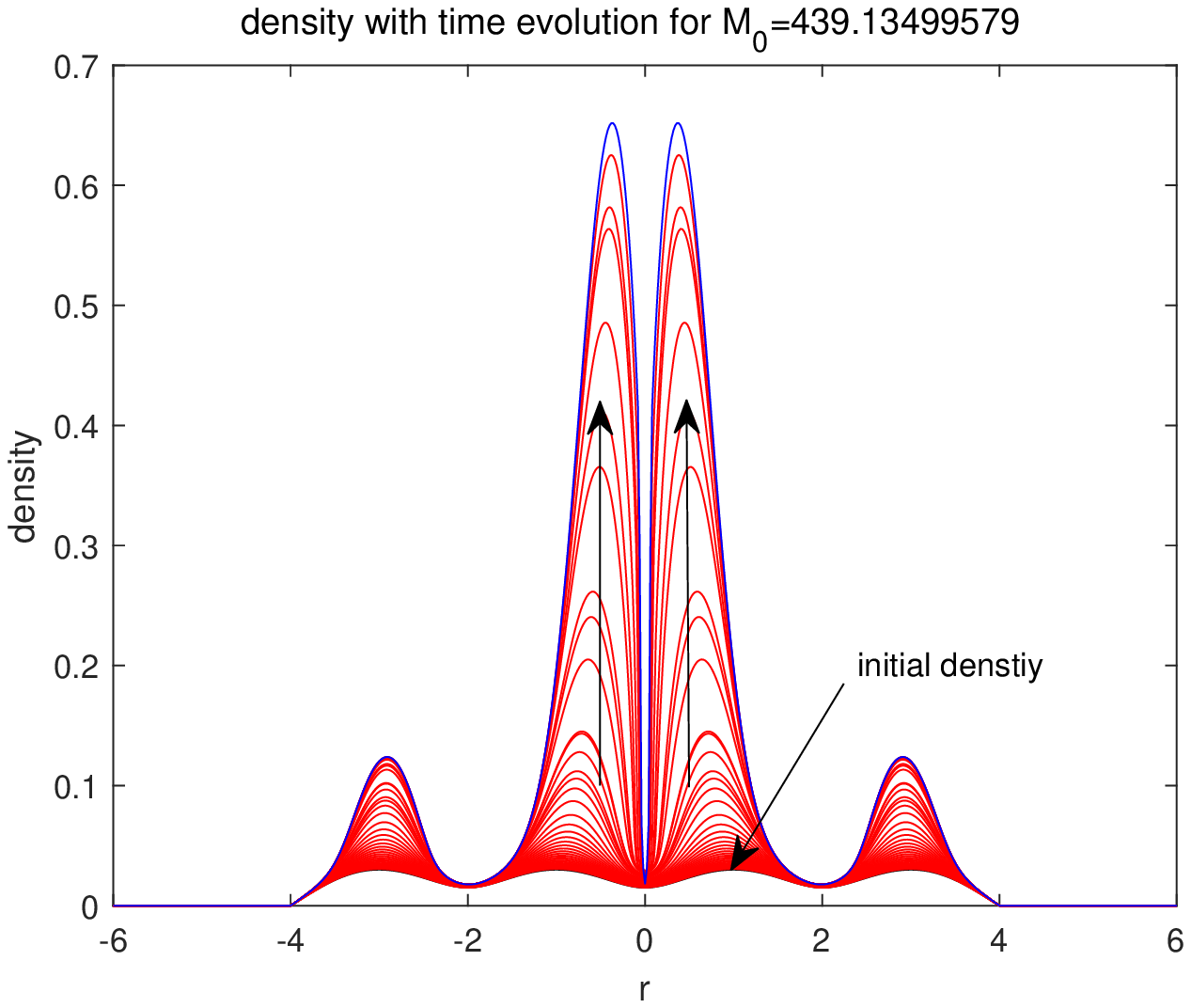}}
\centerline{(c) two singularities}
\end{minipage}
\caption{Different configurations for different $M_0$: (a) the initial density with $m_0=97.58555462$, (b) $M_0=292.75666386$, formation of one local maximum at the blow-up point, (c) $M_0=439.13499579$, formation of two local maxima that determine the position of two blow-up singularities. } \label{fig3}
\end{figure}

In the process of simulations, we found that $M_0,m_0$ directly affect the blow-up point and the behavior of the blow-up solution. We begin with a multi-bump initial data (see Figure \ref{fig3}(a)) with mass $m_0=97.58555462$ in dimension $n=5$. For $M_0>283.97396395$, the solution will blow up at finite time, otherwise it will exist globally. It is interesting to note that for $M_0=292.75666386$, the bumps merge to become a single bump and only one singularity rather than four, ultimately occurs, see Figure \ref{fig3}(b). While for larger $M_0=439.13499579$, it's investigated that two singularities appear near the center, see Figure \ref{fig3}(c).

\begin{figure}[htbp]
\begin{minipage}{0.5\linewidth}
\vspace{3pt}
\centerline{\includegraphics[width=8cm]{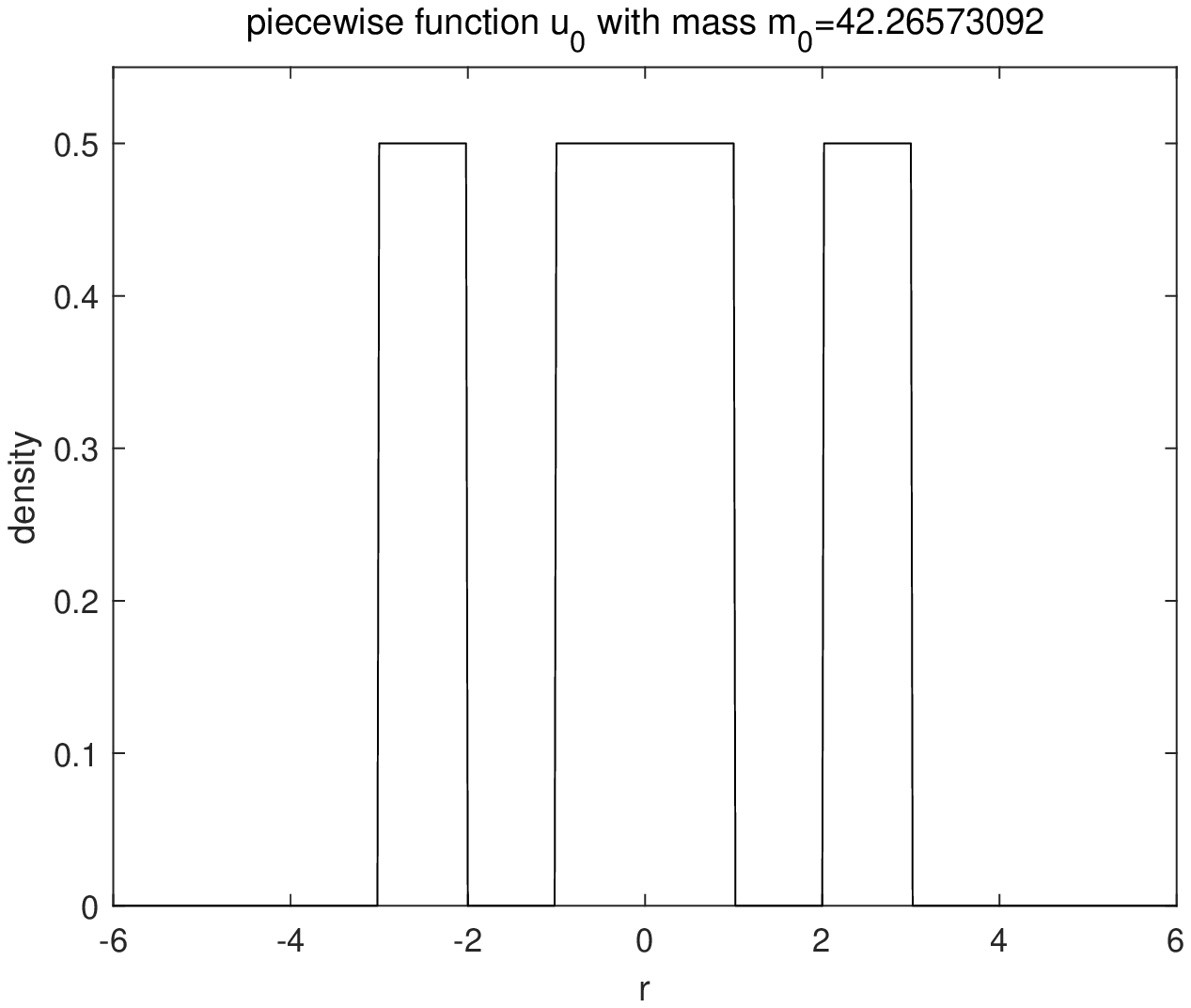}}
\centerline{(a1) piecewise function $u_0$}
\end{minipage}
\begin{minipage}{0.5\linewidth}
\vspace{3pt}
\centerline{\includegraphics[width=8cm]{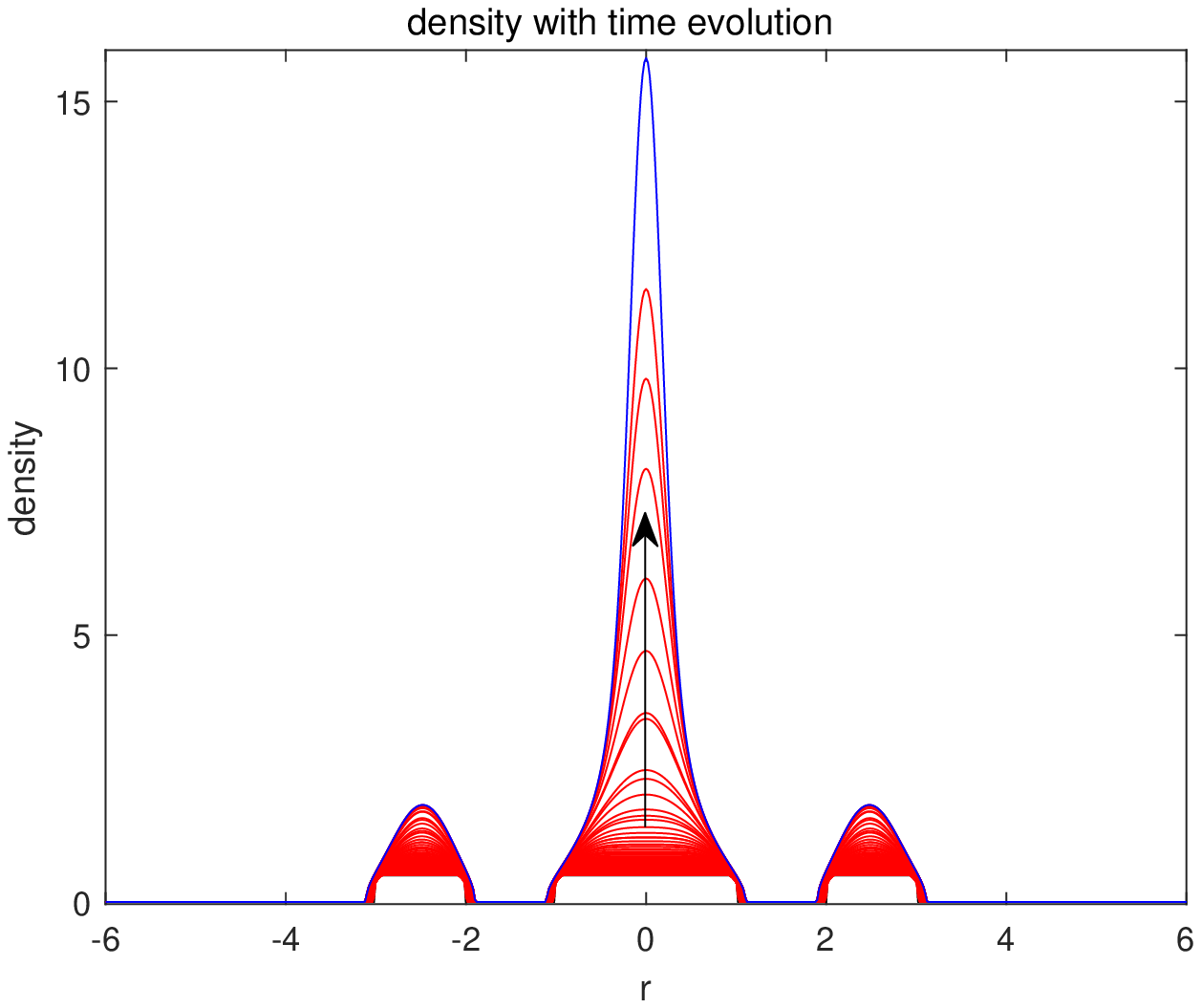}}
\centerline{(a2) formation of one singularity}
\end{minipage}

\begin{minipage}{0.5\linewidth}
\vspace{3pt}
\centerline{\includegraphics[width=8cm]{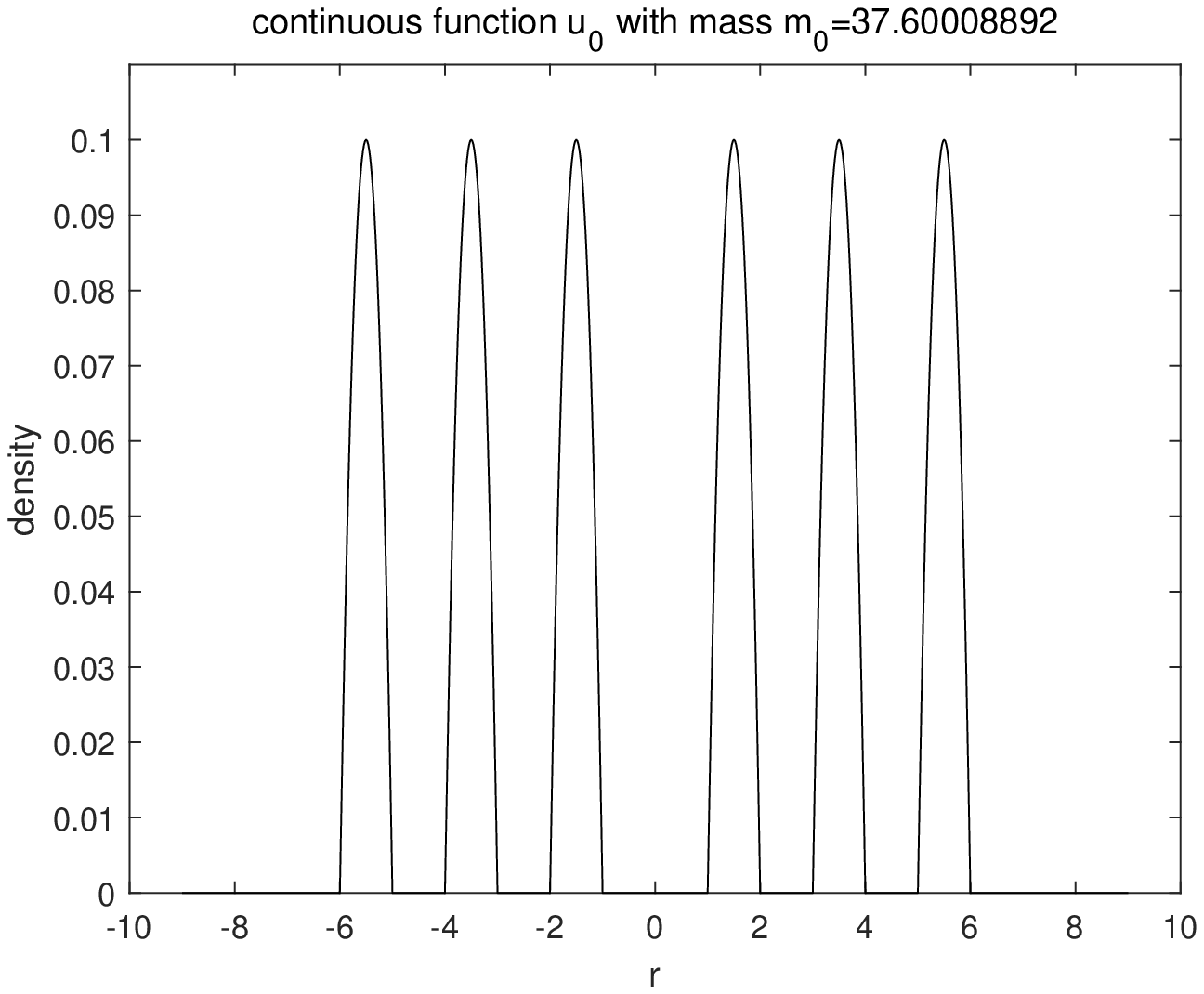}}
\centerline{(b1) the initial data with six supported bumps}
\end{minipage}
\begin{minipage}{0.5\linewidth}
\vspace{3pt}
\centerline{\includegraphics[width=8cm]{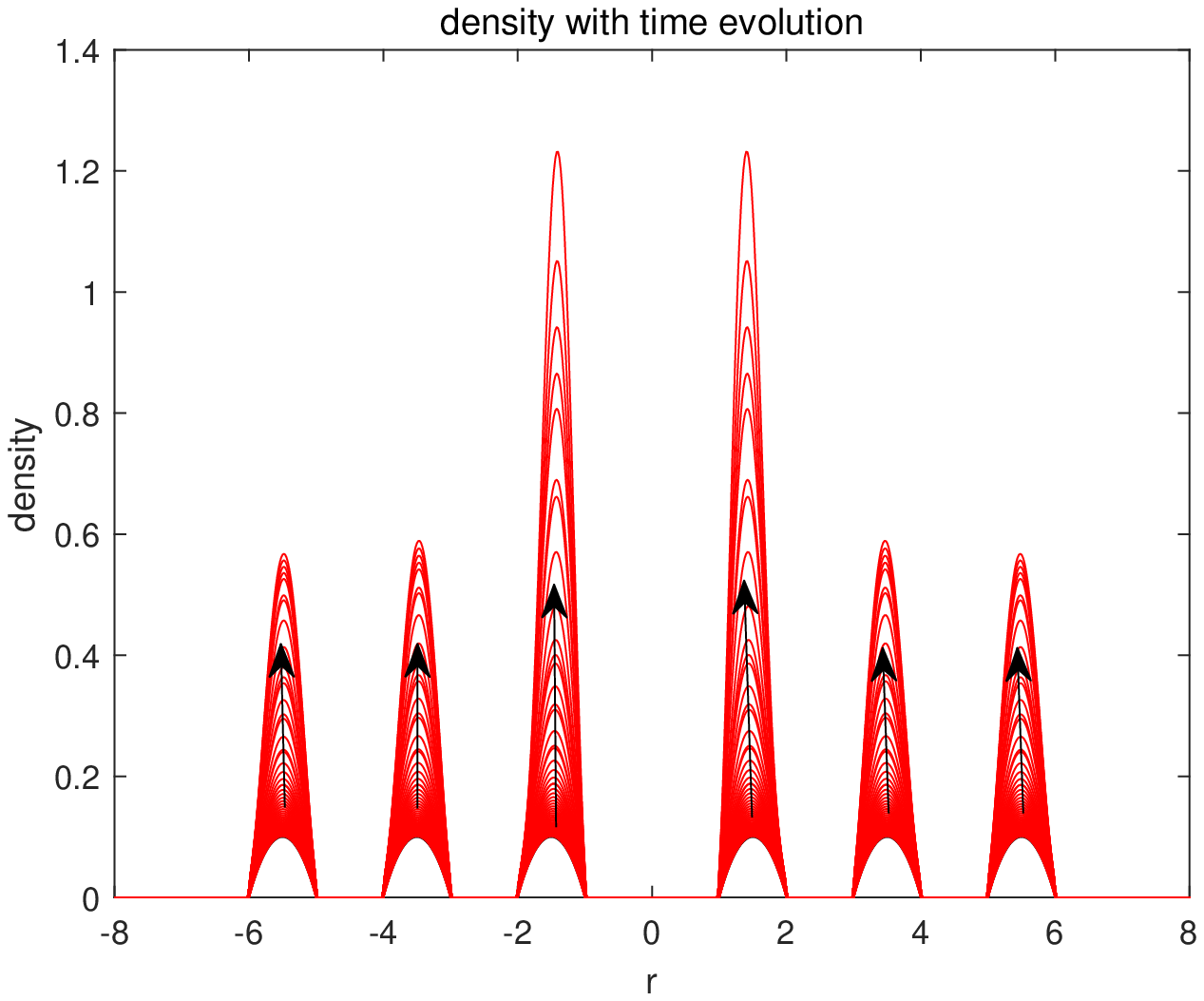}}
\centerline{(b2) formation of two singularities}
\end{minipage}
\caption{Time profiles for solutions approaching blow-up with different forms of initial data for $n=3$: (a) the piecewise function evolves to one local maximum at zero, (b) the six disjointly-supported bumps evolve independently to form two local maxima next to zero. } \label{fig4}
\end{figure}

An important property of degenerate diffusion equations like \er{radialPDE} is the nonlinear superposition principle for disjointly-supported solutions \cite{witel04}. That is, as long as their respective regions of support do not overlap, any combination of non-negative solutions can be pasted together in the domain to yield another configuration. For times when there is no overlap of domains, each bump evolves independently of the others, examples of this are shown in Figure \ref{fig4}(a2) where the central bump and the two compactly-supported bumps away from zero develop independently and finally the central bump evolves to one singularity at zero, and Figure \ref{fig4}(b2) where six supported bumps develop respectively and ultimately two singularities next to zero are formed.

\subsection{Finite time blow-up and infinite-time convergence for the supercritical case $\alpha>m+2/n$}

Theorem \ref{super} tells us that the solution of \er{radialPDE} will exist globally for $\|u_0\|_{L^{p_0}(\R^n)}<C_{p_0}$. In this subsection, we will numerically explore the global existence and finite time blow-up of solutions under the assumption of $\|u_0\|_{L^{p_0}(\R^n)} \ge C_{p_0}$.

\begin{figure}[htbp]
\begin{minipage}{0.33\linewidth}
\vspace{3pt}
\centerline{\includegraphics[width=5.9cm]{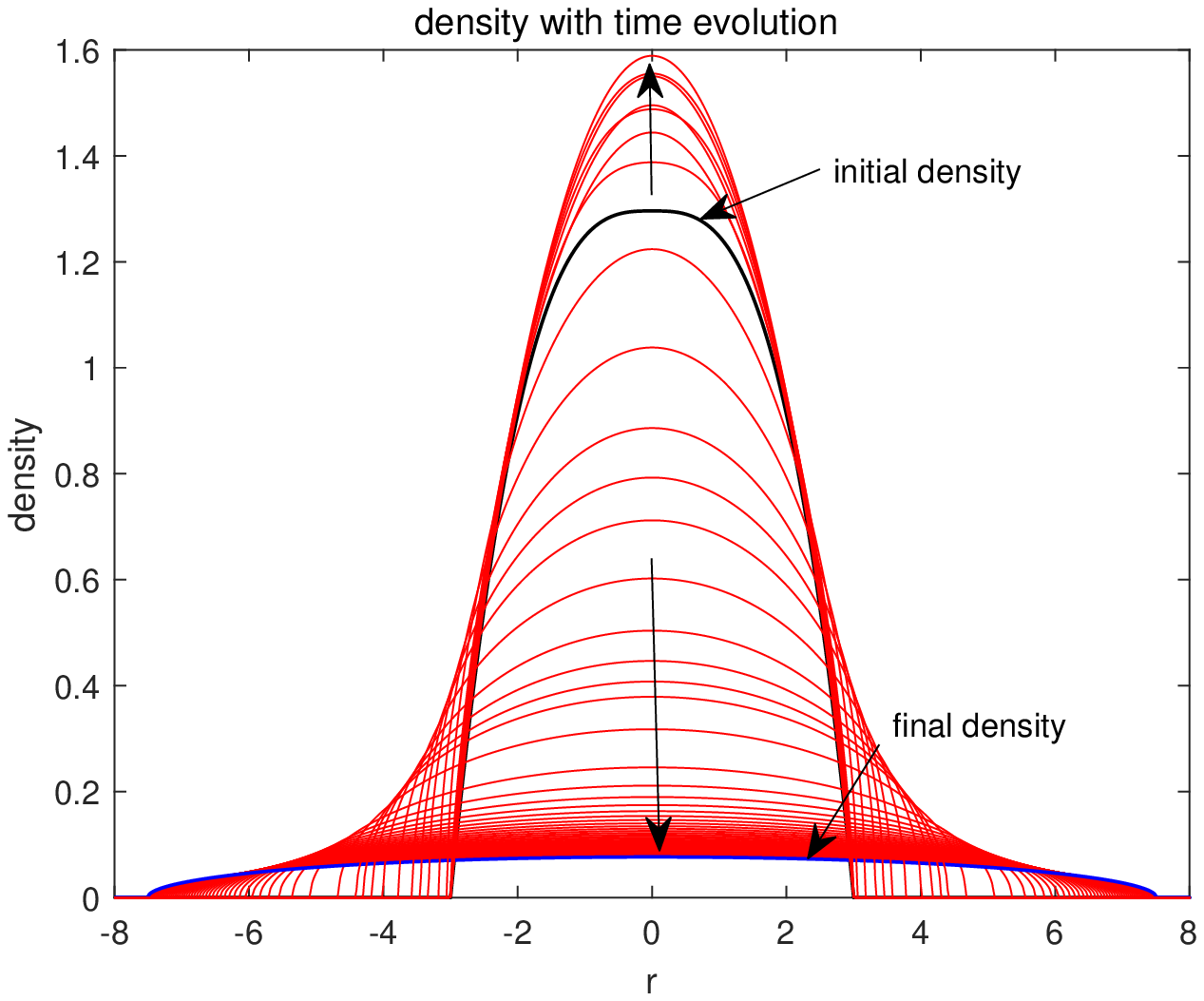}}
\centerline{(a1) density with time evolution}
\end{minipage}
\begin{minipage}{0.33\linewidth}
\vspace{3pt}
\centerline{\includegraphics[width=5.9cm]{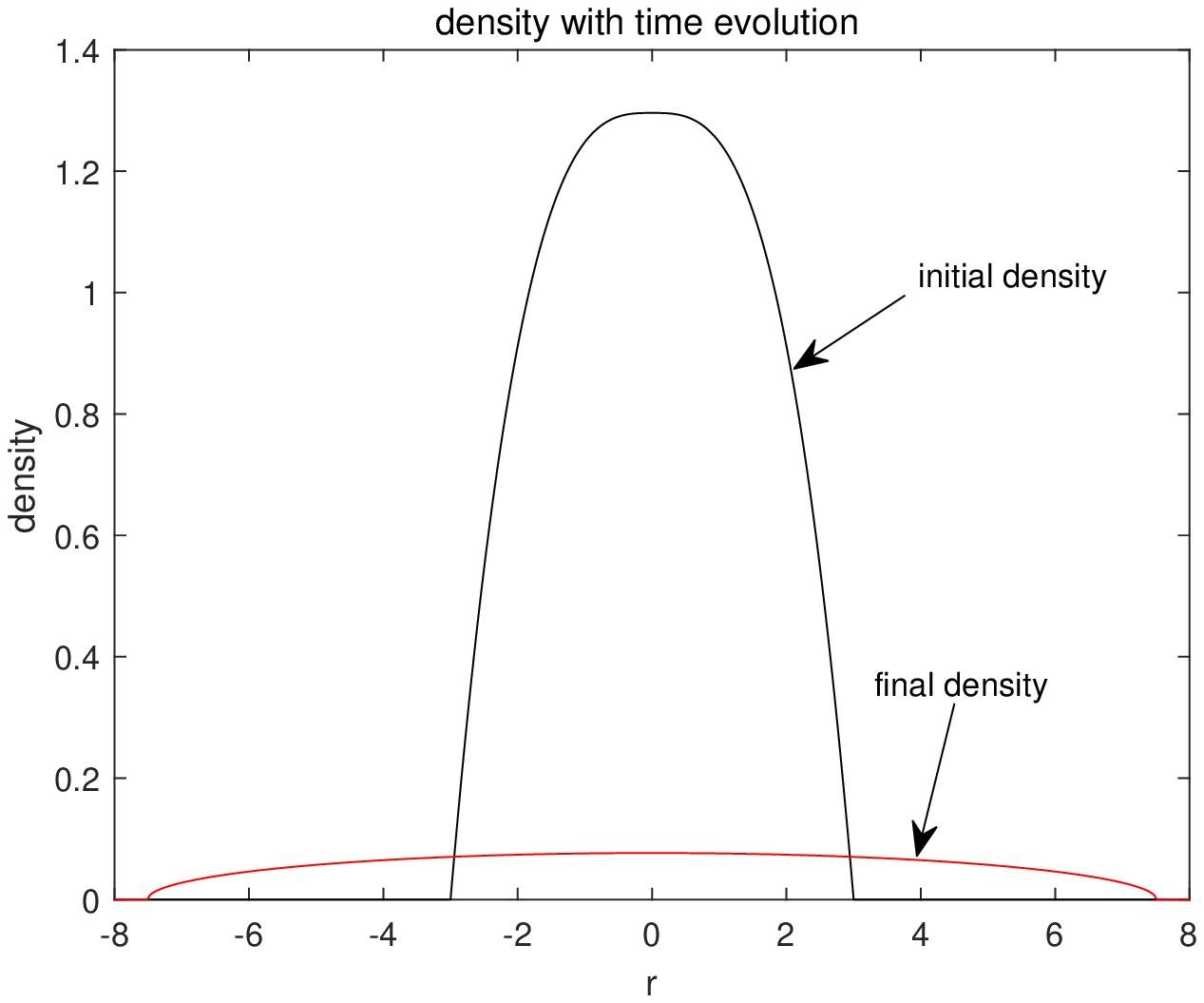}}
\centerline{(a2) initial and final densities}
\end{minipage}
\begin{minipage}{0.33\linewidth}
\vspace{3pt}
\centerline{\includegraphics[width=5.9cm]{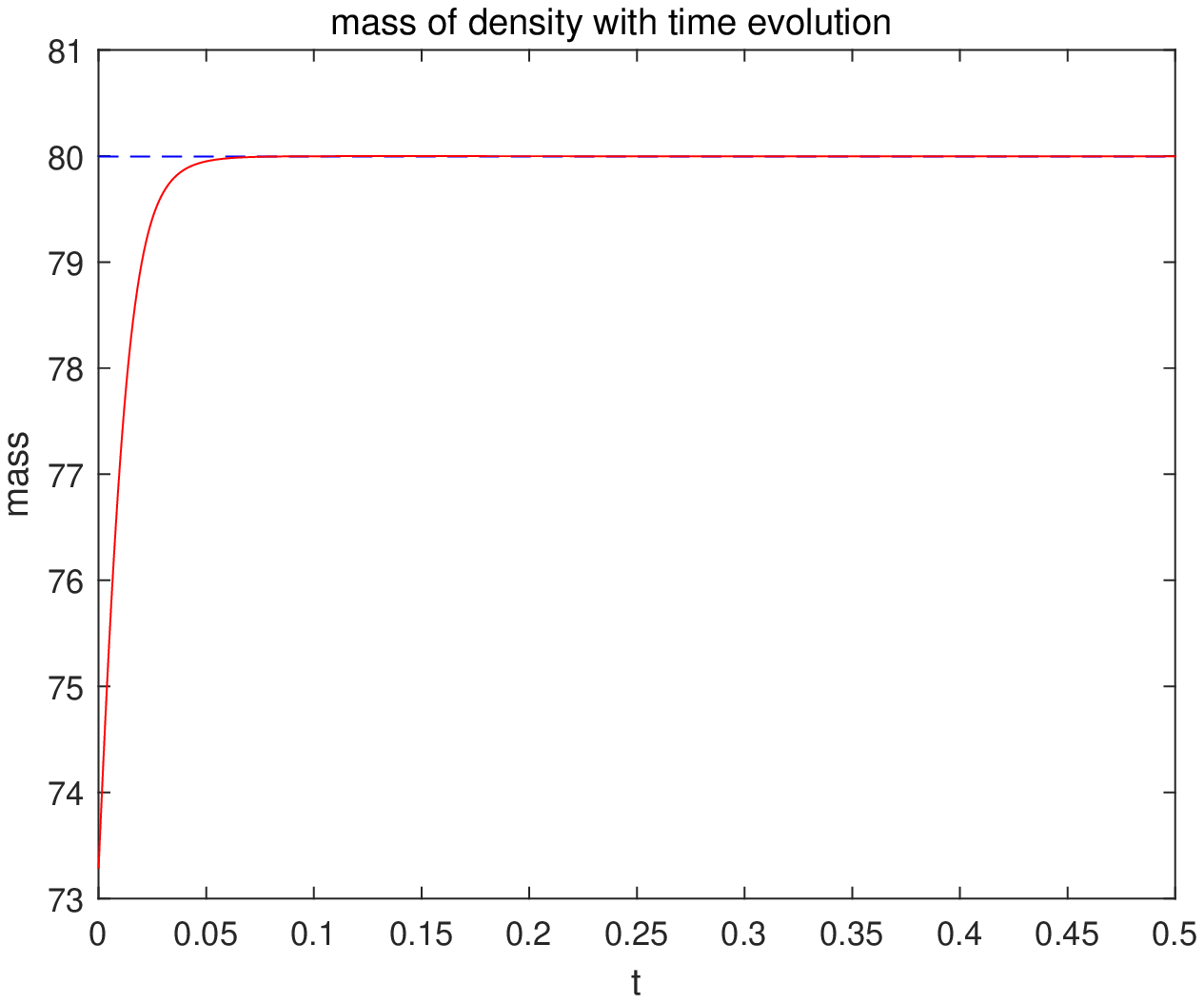}}
\centerline{(a3) mass with time evolution}
\end{minipage}

\begin{minipage}{0.33\linewidth}
\vspace{3pt}
\centerline{\includegraphics[width=5.9cm]{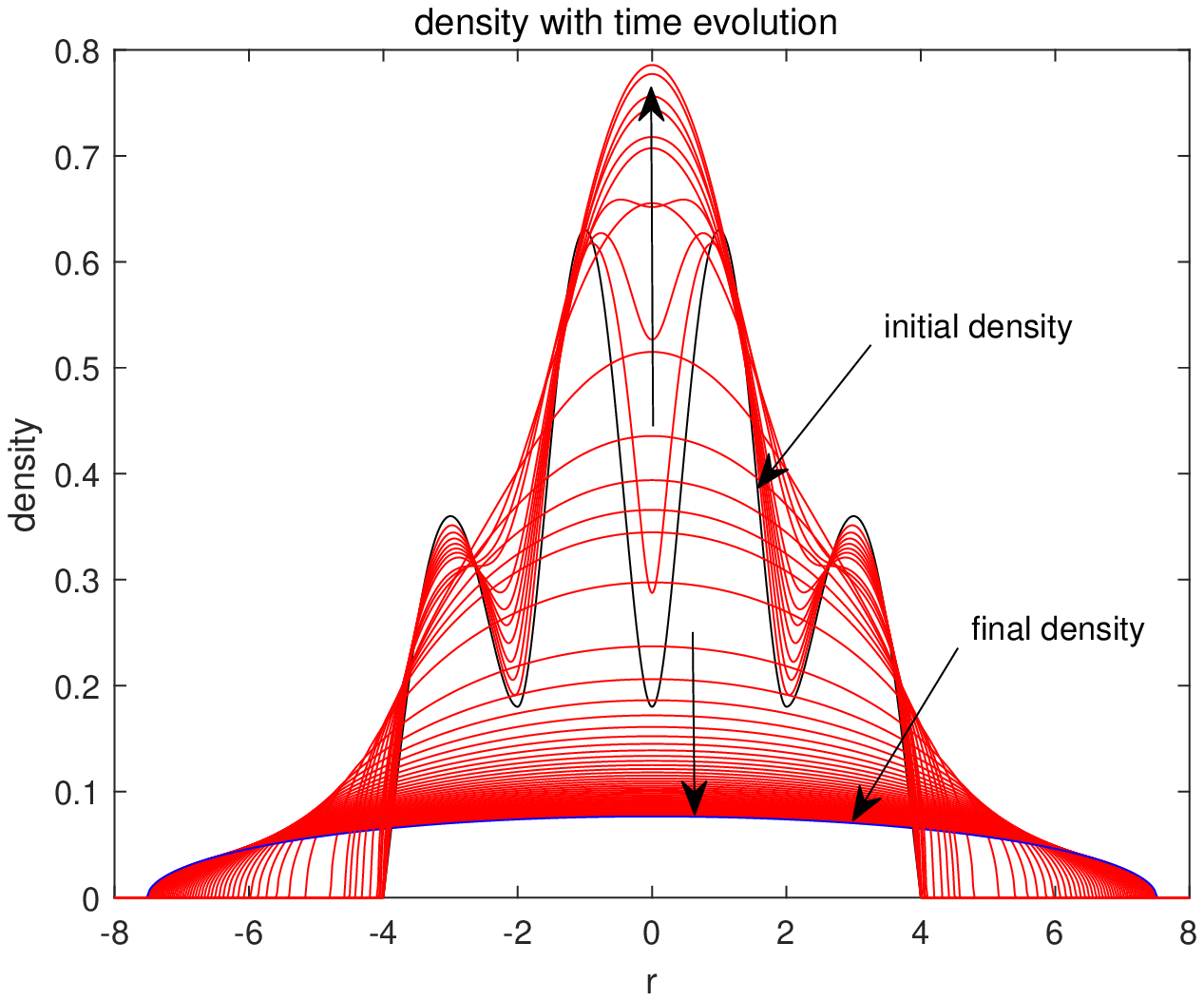}}
\centerline{(b1) density with time evolution}
\end{minipage}
\begin{minipage}{0.33\linewidth}
\vspace{3pt}
\centerline{\includegraphics[width=5.9cm]{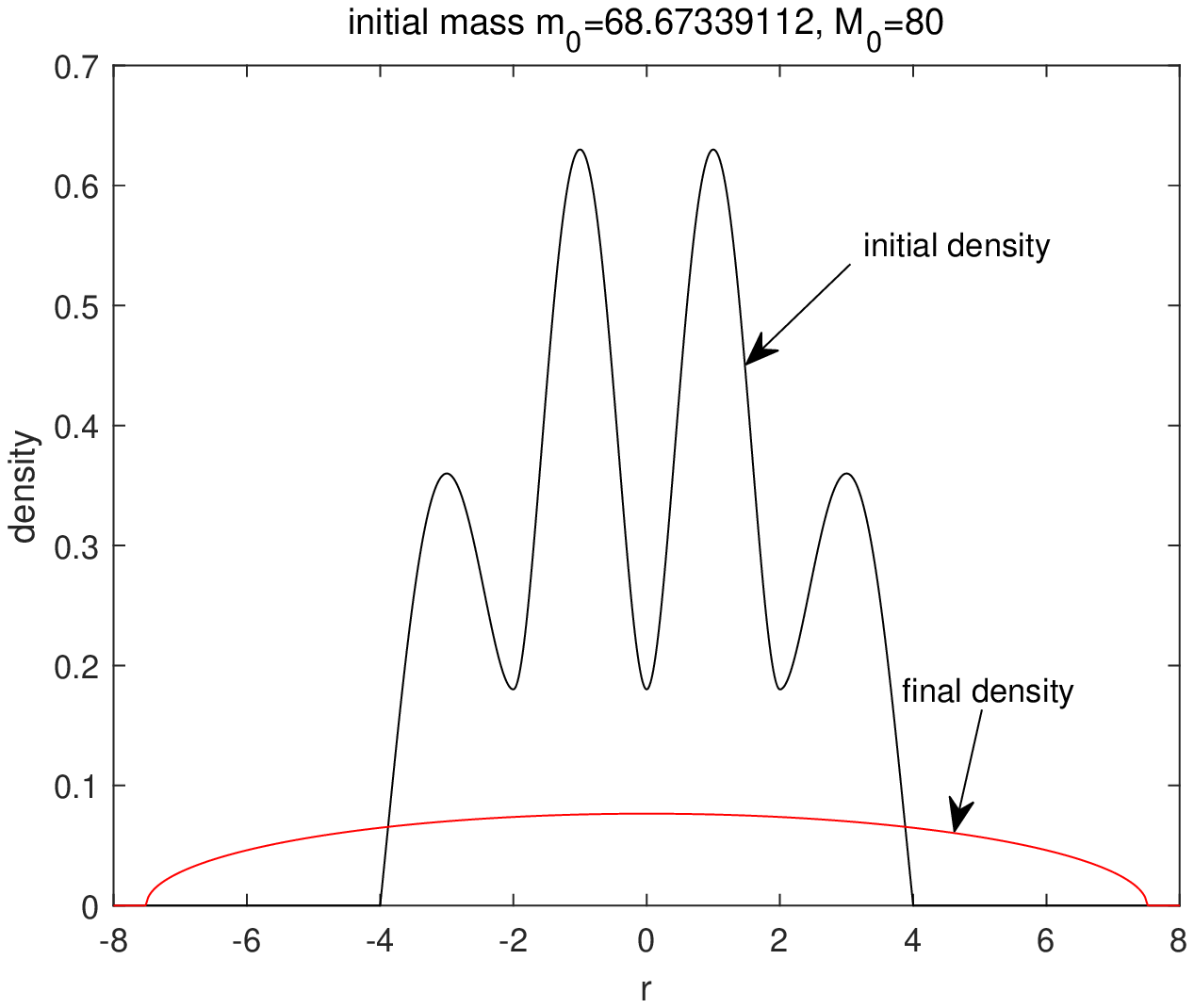}}
\centerline{(b2) initial and final densities}
\end{minipage}
\begin{minipage}{0.33\linewidth}
\vspace{3pt}
\centerline{\includegraphics[width=5.9cm]{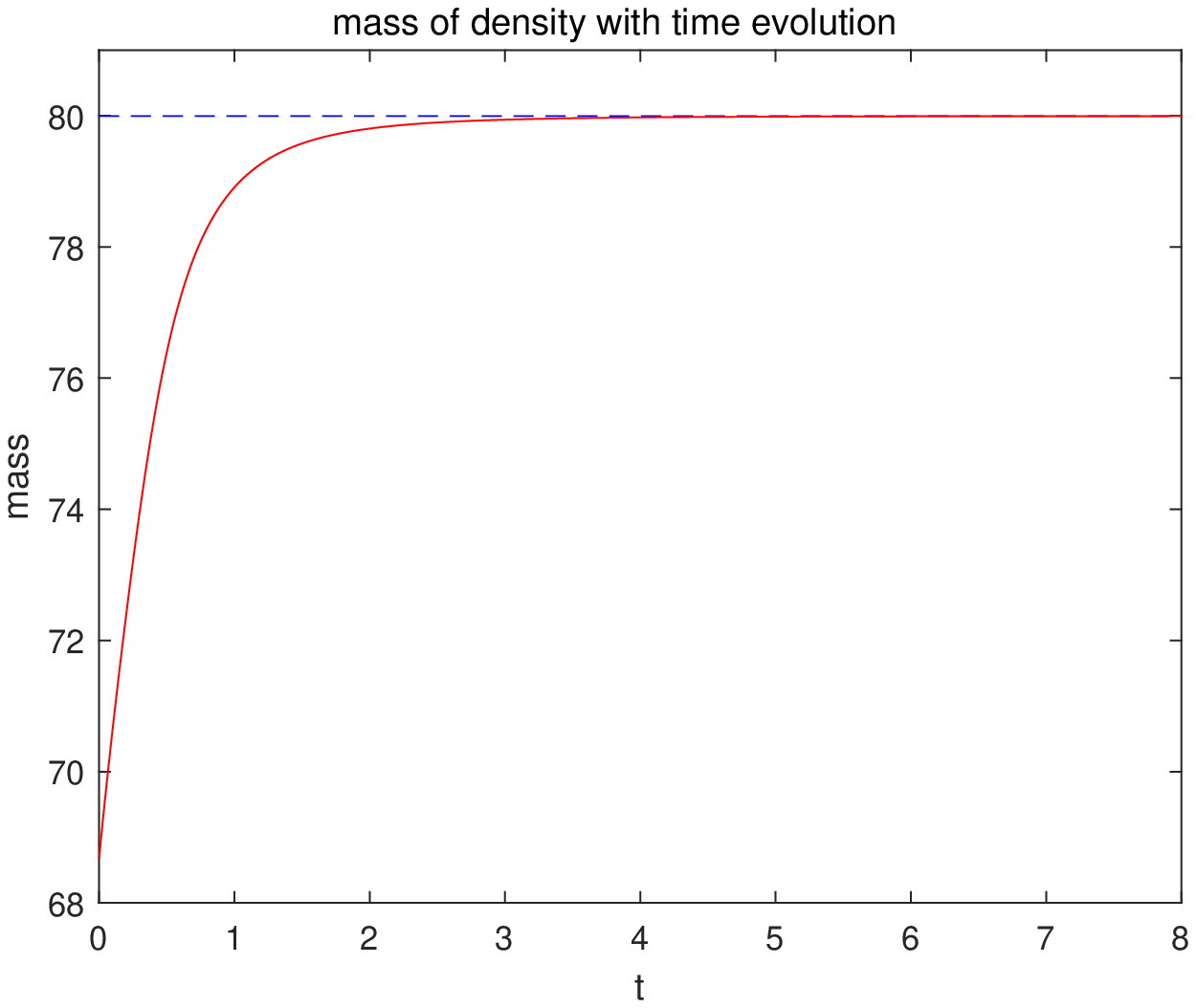}}
\centerline{(b3) mass with time evolution}
\end{minipage}
\caption{Convergence to the steady profile: (a1-a3) $M_0=80$, $m_0=73.28218769$, the initial convex function evolves to the unique compactly supported stationary solution with mass $80$, (b1-b3) $M_0=80$, $m_0=68.67339112$, the multi-bump initial data tends to the same steady state with mass $80$.} \label{fig5}
\end{figure}

\begin{figure}[htbp]
\begin{minipage}{0.33\linewidth}
\vspace{3pt}
\centerline{\includegraphics[width=5.9cm]{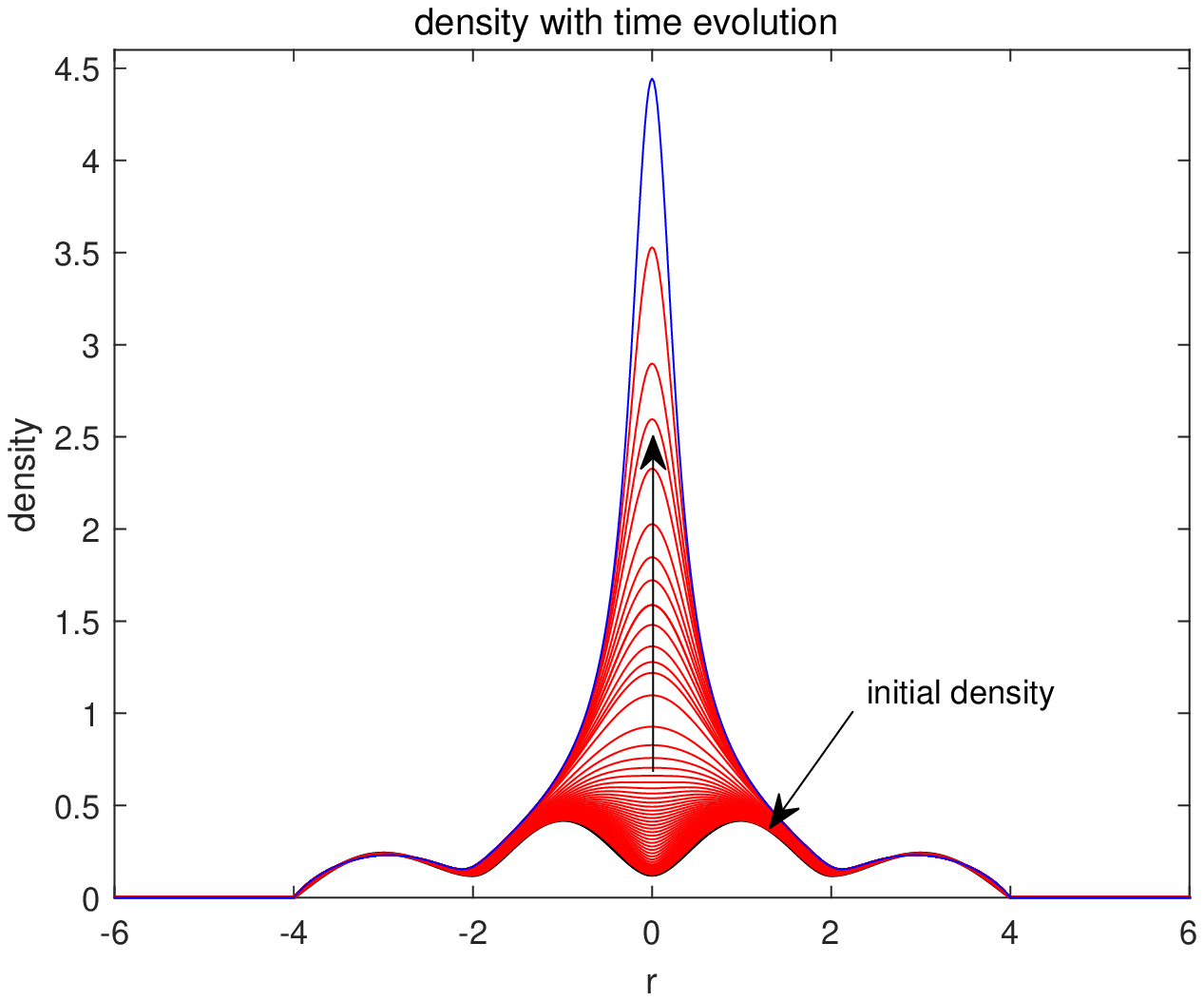}}
\centerline{(a) density with time evolution}
\end{minipage}
\begin{minipage}{0.33\linewidth}
\vspace{3pt}
\centerline{\includegraphics[width=5.9cm]{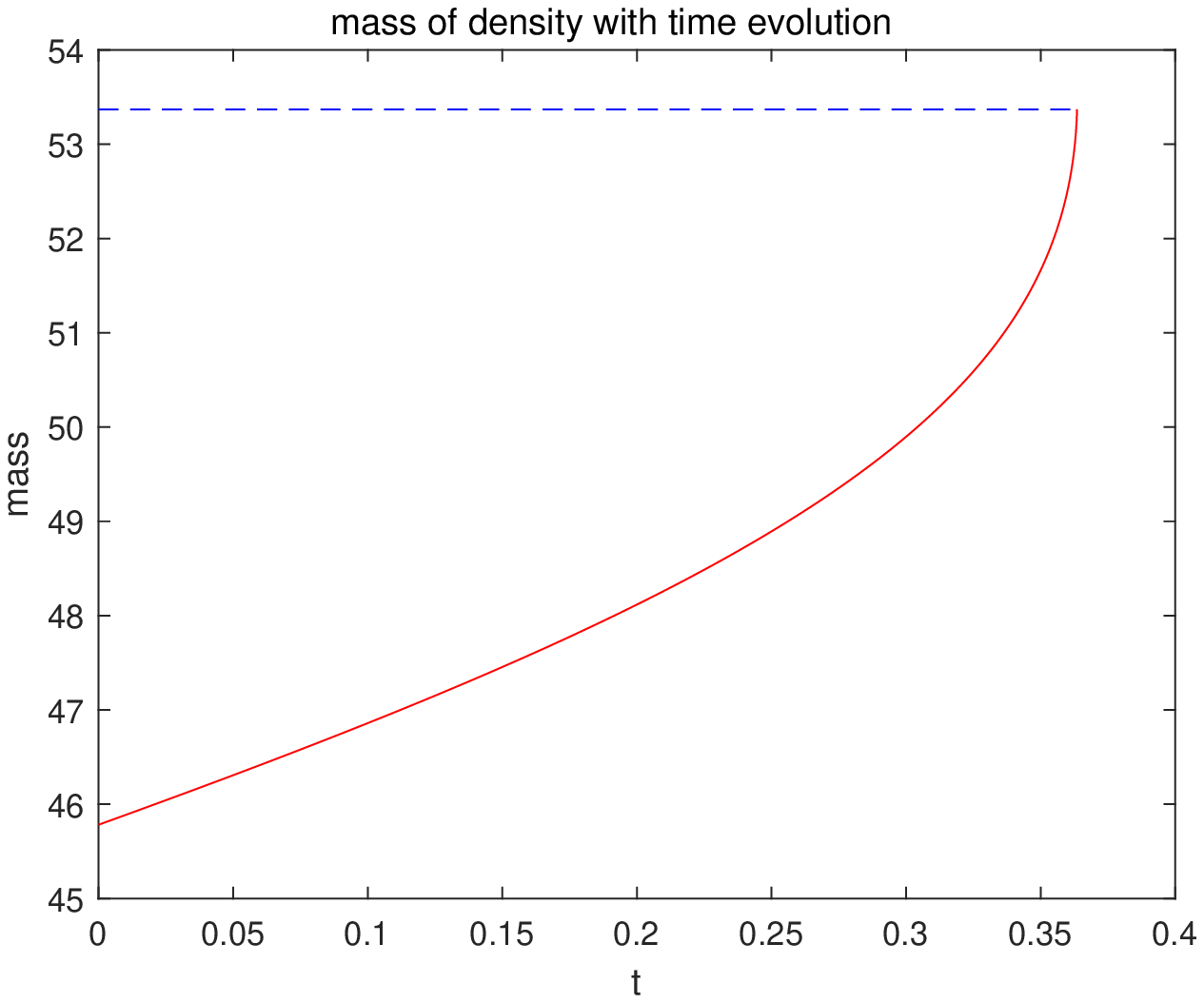}}
\centerline{(b) mass with time evolution}
\end{minipage}
\begin{minipage}{0.33\linewidth}
\vspace{3pt}
\centerline{\includegraphics[width=5.9cm]{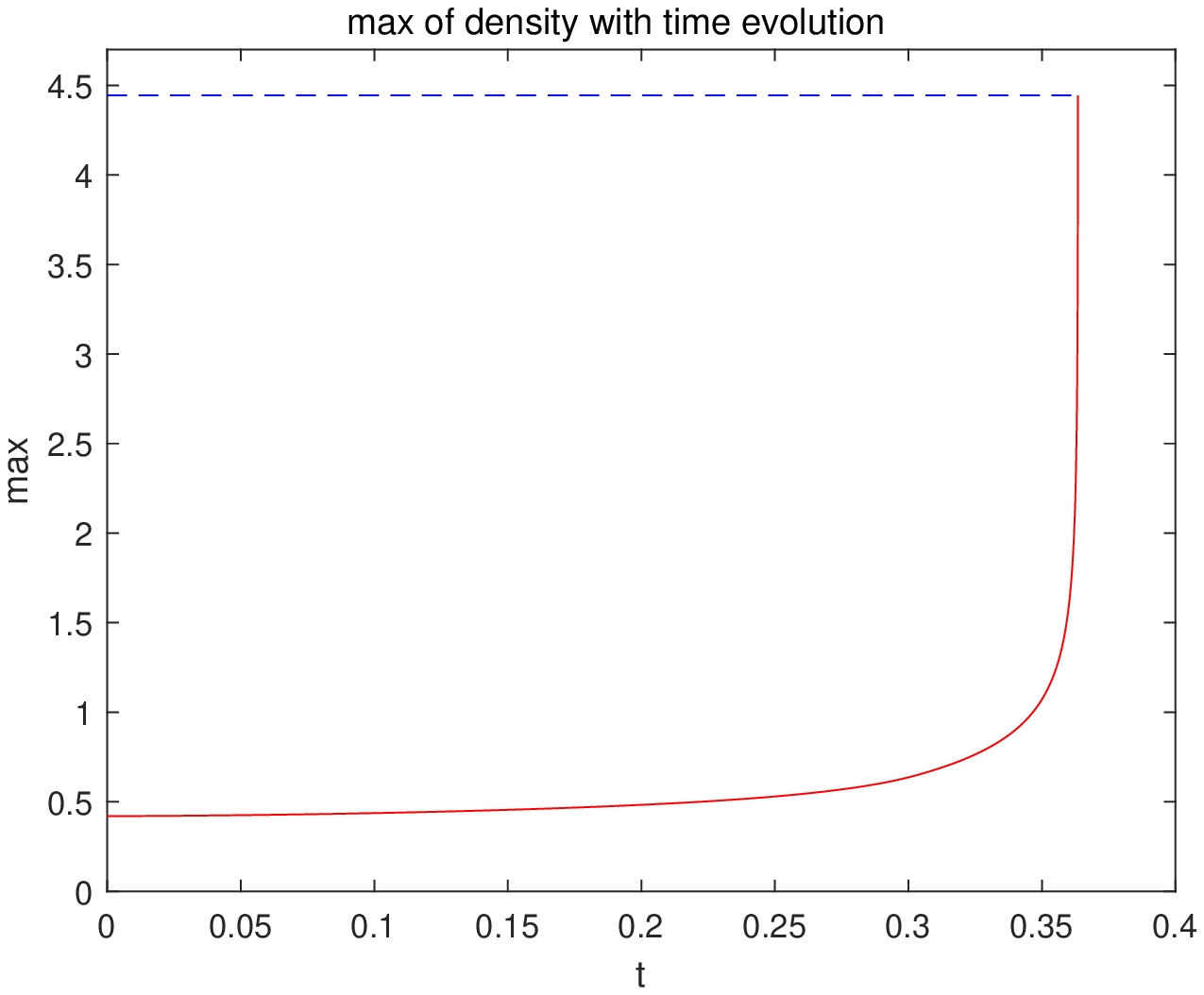}}
\centerline{(c) maximum with time evolution}
\end{minipage}
\caption{The finite time blow-up for $M_0=80$, $m_0=45.78226075$: (a) the multi-bump initial data forms one local maximum that produces blow-up in finite time, (b) the mass increases to $53.36984539$ at time $T_b=0.3635$s, (c) the maximum increases dramatically.} \label{fig6}
\end{figure}

Without loss of generality, we consider the case $n=3.$ Given $M_0=80,$ for $m_0>61.18862434$, that's $\|u_0\|_{L^{p_0}(\R^n)}>11.38458401>C_{p_0}=0.58396290$, the solution will exist globally and converge to the unique compactly supported stationary solution which is only decided by $M_0$ but not on $m_0,$ see Figure \ref{fig5}(a2)(b2). On the contrary, for $m_0<61.18862434$ and $C_{p_0}=0.58396290<\|u_0\|_{L^{p_0}(\R^n)}<11.38458401$, the solution will blow up in finite time, see Figure \ref{fig6} for the formation of a singularity.

\begin{figure}[htbp]
\begin{minipage}{0.5\linewidth}
\vspace{3pt}
(a) \centerline{\includegraphics[width=8cm]{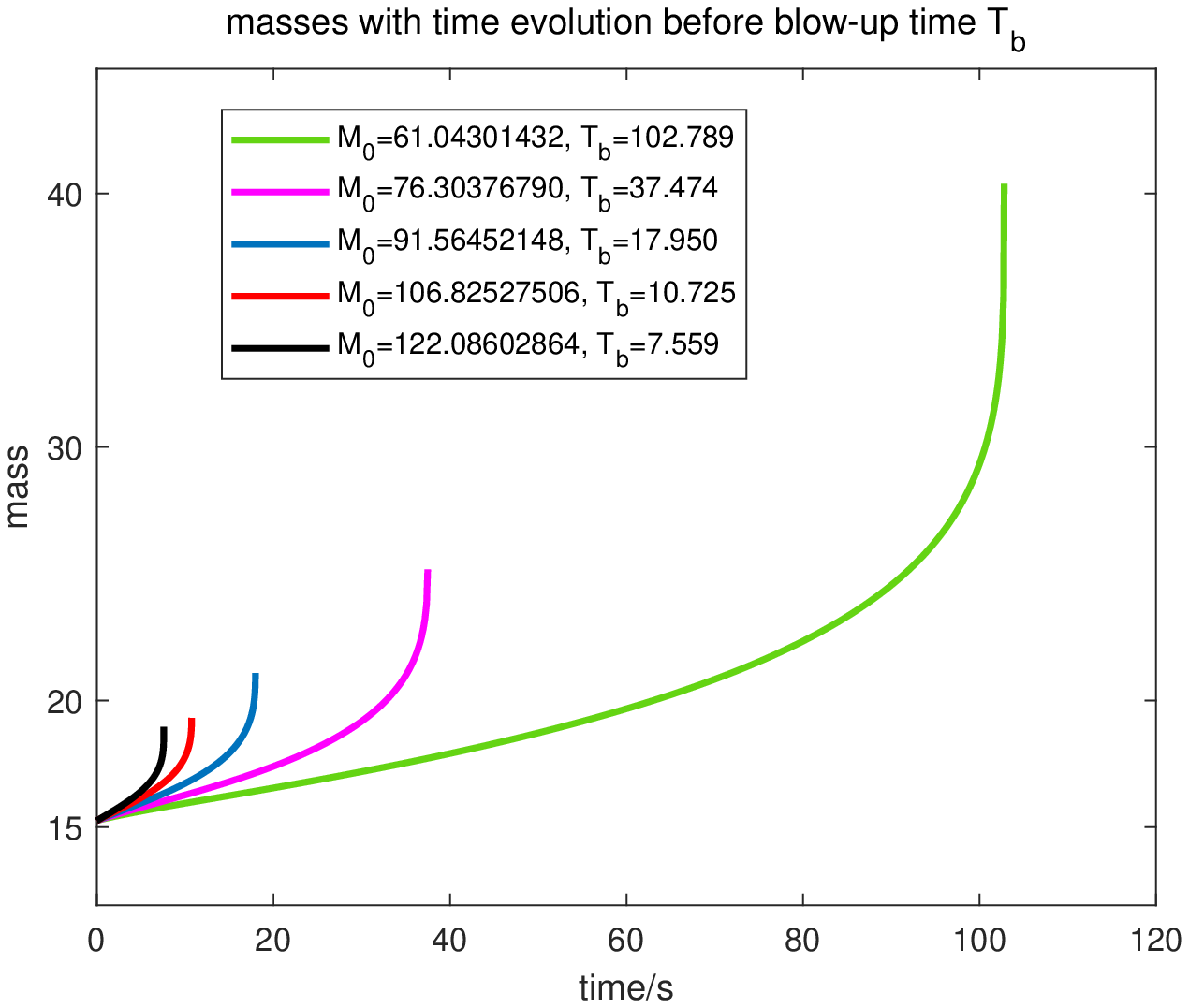}}
\end{minipage}
\begin{minipage}{0.5\linewidth}
\vspace{3pt}
(b) \centerline{\includegraphics[width=8cm]{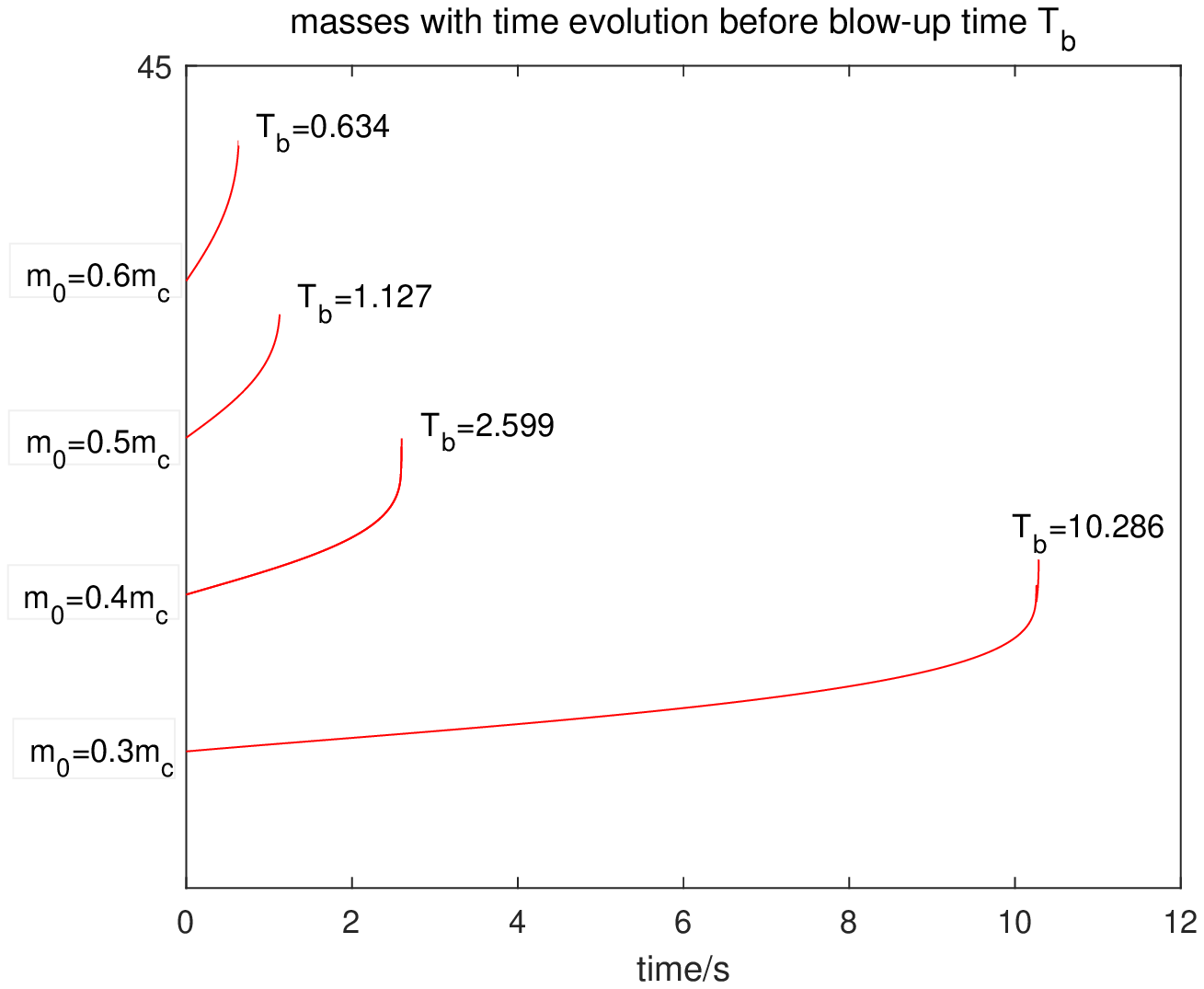}}
\end{minipage}
\caption{ The influence of $M_0,m_0$ on the blow-up time $T_b$: (a) $m_0=15.26075358$, $M_0>58.53590504$, five branches of masses with time evolution before blow-up time $T_b$, (b) $M_0=80,$ $m_0<m_c=61.18862433$, four branches of masses with time evolution before blow-up time $T_b$.} \label{fig7}
\end{figure}

Furthermore, in order to explore the relationship among $M_0,m_0$ and the blow-up time $T_b$, some numerical experiments by choosing different $M_0,m_0$ are carried out. It is shown that for any given $m_0$ and its corresponding $M_c$, larger $M_0$ leads to shorter blow-up time when $M_0>M_c$, see Figure \ref{fig7}(a). Similarly, for any given $M_0$, there exists a critical threshold $m_c$ such that for $m_0<m_c$, the blow-up time is longer when $m_0$ is smaller, see Figure \ref{fig7}(b).

\section{Conclusions} \label{sec6}
\def\theequation{6.\arabic{equation}}\makeatother
\setcounter{equation}{0}
\def\thetheorem{6.\arabic{theorem}}\makeatother
\setcounter{theorem}{0}

In this paper, we have identified a critical exponent $\alpha=m+2/n$ for  equation \er{nkpp} that separates blow-up solutions from those that exist globally, and we have described the structure in three cases. Firstly, for the subcritical case $1\le \alpha<m+2/n$ (the diffusion dominates for large population), the global existence and uniformly boundedness of a weak solution to \er{nkpp} are obtained. Then for the critical case $\alpha=m+2/n$, there exists an upper bound $M_\ast$ such that the solution exists globally for $m_0<M_0<M_\ast.$  Thirdly, for the supercritical case $\alpha>m+2/n$ where the logistic growth term dominates for large population, there exists a universal constant $C_{p_0}$ such that the solution will exist globally with initial data satisfying $\|u_0\|_{L^{p_0}(\R^n)}<C_{p_0}$. In addition, for the critical case $\alpha=m+2/n$ and the supercritical case $\alpha>m+2/n$, numerical simulations illustrate that for any given initial mass $m_0,$ there exists a critical threshold $M_c$ to separate finite time blow-up and global existence. Precisely, at the time that the population grows, the competitive effect of the logistic growth term becomes more influential such that finite time aggregation occurs when $M_0>M_c$ (the position and number of the blow-up points will change for different $M_0$). Inversely, for $m_0<M_0<M_c,$ $M_0$ is too small to prevent the population from spreading around although the growth term is dominant at the beginning (see Figure \ref{fig1}(a3)), finally the density converges to the compactly supported steady profile. This suggests us to further describe that how is $M_c$ selected for any given $m_0$, which is a challenging question in our future research.

\end{document}